\documentclass[12pt]{amsart}
\usepackage{amssymb, amscd, tikz}

\newcommand{\field}[1]{\mathbb{#1}}
\newcommand{\CC}{\field{C}}
\newcommand{\FF}{\field{F}}
\newcommand{\NN}{\field{N}}
\newcommand{\RR}{\field{R}}
\newcommand{\TT}{\field{T}}
\newcommand{\ZZ}{\field{Z}}

\newcommand{\Aa}{\mathcal A}

\newcommand{\Ff}{\mathcal F}
\newcommand{\Gg}{\mathcal G}

\newcommand{\Kk}{\mathcal K}
\newcommand{\Ll}{\mathcal L}

\newcommand{\Oo}{\mathcal O}
\newcommand{\Tt}{\mathcal T}

\newcommand{\NO}[1]{\operatorname{\mathcal{NO}}_{\!#1}}

\newcommand{\fvcl}[1]{\mathcal{P}_{\rm fin}^{\vee}(#1)}
\newcommand{\reduced}{r}
\newcommand{\normal}{n}
\newcommand{\Aut}{\operatorname{Aut}}

\newcommand{\CNP}{\operatorname{CNP}}
\newcommand{\id}{\operatorname{id}}

\newcommand{\lsp}{\operatorname{span}}
\newcommand{\clsp}{\operatorname{\overline{span\!}\,\,}}
\newcommand{\op}{\operatorname{op}}
\newcommand{\CNPgaug}{\nu}
\newcommand{\redgaug}{\nu^{\normal}}

\newcommand{\iotale}{\tilde{\iota}\hskip1pt}

\newcommand{\intfrm}[1]{\Pi{#1}}

\newcommand{\Tc}{\Tt_{\rm cov}}

\newcommand{\supp}{\operatorname{supp}}

\providecommand{\abs}[1]{\lvert#1\rvert}
\providecommand{\norm}[1]{\lVert#1\rVert}

\theoremstyle{plain}
\newtheorem{theorem}{Theorem}[section]
\newtheorem*{theorem*}{Theorem}
\newtheorem*{prop*}{Proposition}
\newtheorem{cor}[theorem]{Corollary}
\newtheorem{lemma}[theorem]{Lemma}
\newtheorem{prop}[theorem]{Proposition}

\theoremstyle{remark}
\newtheorem{rmk}[theorem]{Remark}

\newtheorem{example}[theorem]{Example}

\theoremstyle{definition}
\newtheorem{dfn}[theorem]{Definition}

\newtheorem{notation}[theorem]{Notation}

\numberwithin{equation}{section}

\begin{document}

\title[Co-universal algebras of product systems]
{Co-universal algebras associated to product systems, and gauge-invariant uniqueness theorems}
\author[Toke M. Carlsen]{Toke M. Carlsen}
\address{Department of Mathematics, Norwegian University of Science and Technology, N-7034 Trondheim, Norway.}
\email{tokemeie@math.ntnu.no}
\author[Nadia S. Larsen]{Nadia S. Larsen}
\address{Department of Mathematics, University of Oslo, PO BOX 1053 Blindern, N-0316 Oslo, Norway.}
\email{nadiasl@math.uio.no}
\author[Aidan Sims]{Aidan Sims}
\address{School of Mathematics and Applied Statistics, University of Wollongong, NSW 2522, Australia.}
\email{asims@uow.edu.au}
\author[Sean T. Vittadello]{Sean T. Vittadello}
\address{School of Mathematics and Applied Statistics, University of Wollongong, NSW 2522, Australia.}
\email{seanv@uow.edu.au}

\thanks{This research was supported by the Australian Research Council, The Research Council of
Norway and The Danish Natural Science Research Council. Part of this work was completed while the first-named three
authors were visiting the Fields Institute.}

\begin{abstract}
Let $X$ be a product system over a quasi-lattice ordered group.
Under mild hypotheses, we associate to $X$ a $C^*$-algebra
which is co-universal for injective Nica covariant
Toeplitz representations of $X$ which preserve the gauge coaction. Under
appropriate amenability criteria, this co-universal
$C^*$-algebra coincides with the Cuntz-Nica-Pimsner algebra
introduced by Sims and Yeend. We prove two key uniqueness
theorems, and indicate how to use our theorems to realise a
number of reduced crossed products as instances of our
co-universal algebras. In each case, it is an easy corollary
that the Cuntz-Nica-Pimsner algebra is isomorphic to the
corresponding full crossed product.
\end{abstract}

%\date{June 23, 2009}

\maketitle

\section{Introduction}

In the late 1970s and early 1980s, Cuntz and Krieger introduced
a class of simple $C^*$-algebras generated by partial
isometries, now known as the Cuntz-Krieger algebras
\cite{Cuntz1977, CK1980}. In 1997, Pimsner developed a
far-reaching generalisation of Cuntz and Krieger's construction
by associating to each right-Hilbert $A$--$A$ bimodule $X$
(also known as a $C^*$-correspondence over $A$) two
$C^*$-algebras $\Tt_X$ and $\Oo_X$ (see \cite{Pimsner1997}). The algebras $\Oo_X$ are
direct generalisations of the Cuntz-Krieger algebras, and are
now known as Cuntz-Pimsner algebras while the algebras $\Tt_X$
are generalisations of their Toeplitz extensions.

Pimsner showed that his construction also generalises crossed
products by $\ZZ$.
%Given a homomorphism $\phi \colon A \to B$, there
%is a standard way to view $B$ as a right-Hilbert $A$--$B$
%bimodule, denoted ${_\phi}B$.
If $\alpha$ is an automorphism of
a $C^*$-algebra $A$, there is a standard way to view $A$ as a right-Hilbert $A$--$A$
 bimodule, denoted  $X = {_\alpha}A$, and the Cuntz-Pimsner algebra $\Oo_X$ is isomorphic
to the crossed product $A \times_\alpha \ZZ$.
%This is striking
%because isomorphism classes of Hilbert bimodules can be
%regarded as the morphisms of a category in which the objects
%are $C^*$-algebras. The assignment $\phi \mapsto {_\phi}A$
%indicates how to relate ordinary $C^*$-homomorphisms to this
%category.
In general, it makes sense to think of right-Hilbert
$A$--$A$ bimodules as generalised endomorphisms of $A$ so that
$\Tt_{X}$ is like a crossed product of $A$ by $\NN$ and $\Oo_X$
is then like a crossed product by $\ZZ$.
% which suggests a kindof dilation of $X$.

In an impressive array of papers \cite{ka1,ka2,ka3} Katsura,
drawing on insight from graph algebras contained in \cite{FMR},
not only expanded Pimsner's theory of $C^*$-algebras associated
with Hilbert bimodules beyond the case of isometric left
actions, but also unified the theory of graph algebras and
homeomorphism $C^*$-algebras under the term topological graphs,
and proved (among other things) uniqueness theorems for his
algebras. For $\Oo_X$ his result says that a Cuntz-Pimsner
representation of the bimodule $X$ generates an isomorphic copy
of $\Oo_X$ precisely when the representation is injective and
admits a gauge action. This type of result, due in genesis to
an Huef and Raeburn \cite{HR1997}, is now commonly known as
``the gauge-invariant uniqueness theorem," and has proven to be
a powerful tool for studying analogues of Cuntz-Krieger
algebras.
%The results in \cite{ka2} generalise earlier results for graph
%algebras, see e.g.\cite{FMR}.

In another direction that also generalises Pimsner's work, Fowler introduced and studied
$C^*$-algebras associated to product systems of Hilbert bimodules \cite{F99}.
Product systems over $(0,\infty)$ of Hilbert spaces were introduced
and studied by Arveson \cite{Arveson1989}, and the concept was
later generalised to other semigroups by Dinh \cite{Dinh1991}
and Fowler \cite{Fowler1999}. Just as a right-Hilbert
$A$--$A$ bimodule can be thought of as a generalised
endomorphism, a product system over a semigroup $P$ of
right-Hilbert $A$--$A$ bimodules can be regarded as an action
of $P$ on $A$ by generalised endomorphisms. To study
$C^*$-algebras associated to such objects, Fowler followed the
lead of Nica \cite{N} who had developed Toeplitz-type algebras
for nonabelian semigroups
%. Nica's work suggested that a
%reasonable theory could be established for semigroups
$P$ which
embed in a group $G$ in such a way as to induce what he called
a quasi-lattice order. Based on Nica's formulation, Fowler
associated $C^*$-algebras $\Tc(X)$ to what he called compactly
aligned product systems $X$ over quasi-lattice ordered groups
$(G,P)$ and established that $\Tc(X)$ had much of the structure
of a twisted crossed product of the coefficient algebra $A$ by
the semigroup $P$, with the ``twist" coming from $X$.

Fowler also associated to each product system a generalised
Cuntz-Pimsner algebra $\Oo_X$. However, $\Oo_X$ need not in
general be a quotient of $\Tc(X)$. Moreover, the canonical
homomorphism from $A$ to $\Oo_X$ is, in general, not injective,
so there is little hope of a gauge-invariant uniqueness
theorem. Recently, however, Sims and Yeend \cite{SY} introduced
what they call the Cuntz-Nica-Pimsner algebra of a product
system $X$ over a quasi-lattice ordered group $(G,P)$. This
algebra is a quotient of $\Tc(X)$, and Sims and Yeend
established that under relatively mild hypotheses the canonical
representation of the product system $X$ on the
Cuntz-Nica-Pimsner algebra $\NO{X}$ is isometric.
%By relating their algebras to Crisp and Laca's
%boundary quotient algebras \cite{CL}, they showed that there is
%a good case for regarding $\NO{X}$ as a crossed product of $A$
%by a generalised action of $G$, suggesting a kind of dilation
%of the product system $X$ to $G$.
However, they were unable to
establish a gauge-invariant uniqueness theorem for $\NO{X}$.

The initial purpose of the research presented in this article
was to understand and describe the fixed-point algebra, called
the \emph{core}, for the canonical coaction of $G$ on $\NO{X}$,
and use this analysis to establish the missing gauge-invariant
uniqueness theorem. Since $\NO{X}$ is defined for pairs $(G,
P)$ in which $G$ need not be abelian, Katsura's gauge action of
$\TT$, equivalently seen as a coaction of $\ZZ$, must be
replaced by a coaction of $G$.

We analyse the core in Section~\ref{sec:core} and we
subsequently prove a gauge-invariant uniqueness theorem in
Corollary~\ref{cor:guit}. This result is quite far-reaching in
itself: in particular, it enables us to recover isomorphisms of
various full and reduced crossed products --- ordinary or
partial --- in the presence of amenability (see
Section~\ref{sec:examples}). More importantly, for the class of
topological higher-rank graphs introduced by Yeend \cite{y} in
his generalisation of Katsura's topological graphs from
\cite{ka1}, the result is new, and its proof follows a rather
different path than earlier proofs in other contexts involving
product systems over $\NN^k$.

However, we do not proceed from the analysis of the core to the
gauge-invariant uniqueness theorem in the usual fashion, and
the main thrust of our results in the later sections of the
paper deals with what happens when $\NO{X}$ does not satisfy a
gauge-invariant uniqueness theorem and with the intriguing
properties of the quotient $\NO{X}^\reduced$ which does.

To discuss the key new idea we introduce in this paper, we
first observe that amenability considerations imply that
$\NO{X}$ will not, in general, satisfy a gauge-invariant
uniqueness theorem. Specifically, suppose that $(G,P)$ is a
quasi-lattice ordered group such that $G$ is not amenable, but
every finite subset of $P$ has a supremum under the
quasi-lattice order (finite-type Artin groups provide examples
of this situation). Define a product system over $P$ by letting
$X_{p} = \CC$ for each $p$. Then $\NO{X}$ is isomorphic to the group
$C^*$-algebra $C^*(G)$, and the quotient map from $C^*(G)$ to
$C^*_r(G)$ preserves the gauge coaction and is injective on the
coefficient algebra but is not injective. It is the
\emph{reduced} group $C^*$-algebra which satisfies a
gauge-invariant uniqueness theorem, but this algebra lacks a
universal property to induce homomorphisms to which this
theorem may be applied.

Since it is the gauge-invariant uniqueness property we are
interested in, we seek an analogue of $C^*_r(G)$ to the context
of $C^*$-algebras associated to product systems. We desire an
``intrinsic'' definition of our $C^*$-algebra which, like the
universal properties of other generalisations of Cuntz-Krieger
algebras, gives us an effective tool for analysis. To this end,
our $C^*$-algebra $\NO{X}^\reduced$ is described in terms of a
co-universal property. Specifically, in
Section~\ref{sec:couniversal} we prove that for $X$ in a large
class of product systems there exists a unique $C^*$-algebra
$\NO{X}^\reduced$ which: (1)~is generated by an injective Nica
covariant Toeplitz representation of $X$; (2)~carries a
coaction of $G$ compatible with the canonical gauge coaction on
Fowler's $\Tc(X)$; and~(3) has the property that given any
other $C^*$-algebra $B$ generated by an injective Nica
covariant representation $\psi$ of $X$ and carrying a coaction
$\beta$ of $G$ compatible with the gauge coaction, there is a
canonical homomorphism $\phi \colon B \to \NO{X}^\reduced$. We
also establish that this homomorphism $\phi$ is injective if
and only if $\psi$ is Cuntz-Nica-Pimsner covariant and $\beta$
is normal.

We identify amenability hypotheses which imply that the
canonical coaction on $\NO{X}$ is normal. It then follows from
our main theorem that $\NO{X}$ and $\NO{X}^\reduced$ coincide
under the same amenability hypotheses. From this we obtain a
gauge-invariant uniqueness property of the usual form for
$\NO{X}$, see Corollary~\ref{cor:guit}.

The basic idea of a co-universal property of a $C^*$-algebra
has appeared before, notably in Exel's work on Fell bundles
\cite{Exel} which we use in our analysis, in Katsura's work
\cite{ka2} on $C^*$-algebras associated to a single bimodule
which our work generalises, and in the work of Laca and
Crisp-Laca on Toeplitz algebras and their boundary quotients
\cite{CL, L}. However, this article is, to our knowledge, the
first time that the co-universal property has been used as the
defining property of a $C^*$-algebra. In
Section~\ref{sec:examples} we make extensive
--- and to our knowledge quite novel --- use of the defining
co-universal property of $\NO{X}^\reduced$ to prove that in
various special cases $\NO{X}^\reduced$ is isomorphic to
appropriate reduced crossed-product like $C^*$-algebras.
%Specifically, we show that: (1) for product systems with
%one-dimensional fibres, we obtain Crisp and Laca's boundary
%quotient algebras as instances of $\NO{X}^\reduced$; and (2)
%for product systems obtained from actions $\alpha$ of $G$ on a
%$C^*$-algebra $A$, if every finite subset of $P$ has a supremum
%in the quasi-lattice order, then we obtain the reduced crossed
%product $A \times_{\alpha,r} G$ as an instance of
%$\NO{X}^\reduced$. In both cases, it follows easily that
%$\NO{X}$ coincides with the corresponding full crossed product.
In particular, we feel justified in regarding the algebras
$\NO{X}^\reduced$ and $\NO{X}$ as reduced- and full twisted
crossed products of the algebra $A$ by a generalised partial
action of the group $G$.

We wish to emphasise the power and utility of the co-universal
property of $\NO{X}^\reduced$. In particular, the co-universal
property involves only the defining relations of the
Nica-Toeplitz algebra and not the Cuntz-Pimsner covariance
condition introduced in \cite{SY}. Since this last is a very
technical relation, and difficult to check in practice, it is a
significant advantage of our approach that we do not need to
check it in any of our applications. In each case, we instead
check that the algebra $A$ which we wish to compare with
$\NO{X}$ is generated by an injective Nica covariant Toeplitz
representation of $X$, use the co-universal property to obtain
a surjective homomorphism $\phi \colon A \to \NO{X}^\reduced$,
and then use properties of $A$ to prove that $\phi$ is
injective. Particularly interesting is that when $A$ is some
sort of reduced crossed product, we then obtain, almost for
free, isomorphism of the corresponding full crossed product
with $\NO{X}$. The point is that proving directly, for example,
the isomorphism between $\NO{X}$ and the full crossed product
associated with Crisp and Laca's boundary quotient algebra
would require using the universal properties in both
directions, and hence checking both the Cuntz-Pimsner relation,
and Crisp and Laca's elementary relations associated with the
essential spectrum of the quasi-lattice ordered group. The
effort could not be reduced by application of a gauge-invariant
uniqueness theorem in either direction because such a theorem
only applies when the full and reduced $C^*$-algebras coincide.

\vskip1ex

The results of the paper are organised in three main sections
following a preliminaries section. In Section~\ref{sec:core},
we analyse the fixed-point algebra in $\Tc(X)$ and establish,
for a large class of product systems, that any representation
of $\NO{X}$ which is injective on the coefficient algebra is
injective on the whole fixed-point algebra. This answers a
question of Sims and Yeend \cite{SY}. In
Section~\ref{sec:couniversal}, we define $\NO{X}^\reduced$ and
prove our uniqueness theorems for it. Using Exel's results, we
also establish a gauge-invariant uniqueness theorem for
$\NO{X}$ under appropriate amenability hypotheses. Finally, in
Section~\ref{sec:examples} we use our theorems, most notably
the co-universal property, to establish for each of a variety
of reduced crossed-product algebras $A$ an isomorphism of $A$
with $\NO{X}^\reduced$ for an appropriate product system $X$.
We also prove in Section~\ref{sec:examples} that Katsura's
construction of a Hilbert bimodule from a topological graph
yields, for each compactly aligned topological higher-rank
graph $\Lambda$ in the sense of Yeend, a compactly aligned
product system $X$ over $\NN^k$ of Hilbert bimodules, and that
for this $X$, $\NO{X}$ is isomorphic to the $C^*$-algebra of
Yeend's boundary-path groupoid of $\Lambda$ (see~\cite{y}). We
have included an appendix detailing how and when coactions
descend to quotients, and when the resulting coaction is
normal. These results are surely known, but were difficult to
locate in the literature, at least in the specific context of
full coactions with which we deal in this paper.

Towards the late stages of completing this paper, the second
named author learned of the possible connection between our
work and that of Arveson in \cite{Arveson2008}. It seems that
the existence of a co-universal algebra for our systems could
probably be derived from Arveson's results. Since Arveson's
algebras are not obtained constructively we believe that our
explicit construction and identification of the co-universal
algebra is of independent interest. We thank Hangfeng Li for
pointing us to \cite{Arveson2008}.

\vskip1ex

\noindent\textbf{Acknowledgements.} We thank each of Narutaka
Ozawa and Iain Raeburn for helpful conversations about
coactions. We thank Marcelo Laca for helpful discussions about
quasi-lattice ordered groups and boundary quotient algebras.
The first three authors acknowledge both the financial support
and the stimulating atmosphere of the Fields Institute.

\section{Preliminaries}

For a discrete group $G$ we write $g\mapsto i_G(g)$ for the
canonical inclusion of $G$ as unitaries in the full group $C^*$-algebra
$C^*(G)$. We let
$\lambda_G$ denote the left regular representation of $G$.
We denote by $A$ a $C^*$-algebra. An unadorned tensor product
of $C^*$-algebras will denote the minimal tensor product.

Much of what follows is a summary of \cite[Section~2]{SY}. We
refer the reader to \cite{Black, Lan, TFB} for more detail.

\subsection{Hilbert bimodules}

A right-Hilbert $A$-module is a
complex vector space $X$ endowed with a right $A$-module
structure and an $A$-valued $A$-sesquilinear form $\langle
\cdot\,,\,\cdot\rangle_A$ ($^*$-linear in the first variable)
such that $X$ is complete in the norm $\|x\|_A := \|\langle x,
x \rangle_A \|^{1/2}$. A map $T \colon X \to X$ is said to be
\emph{adjointable} if there is a map $T^* \colon X \to X$ such
that $\langle Tx, y\rangle_A = \langle x, T^*y\rangle_A$ for
all $x,y \in X$. Every adjointable operator on $X$ is
norm-bounded and linear, and the adjoint $T^*$ is unique. The
collection $\Ll(X)$ of adjointable operators on $X$ endowed
with the operator norm is a $C^*$-algebra. The ideal of
generalised compact operators $\Kk(X) \lhd \Ll(X)$ is the
closed span of the operators $x \otimes y^* \colon z \mapsto x
\cdot\langle y, z \rangle_A$ where $x$ and $y$ range over $X$.

A \emph{right-Hilbert $A$--$A$ bimodule} is a right-Hilbert $A$
module $X$ endowed with a left action of $A$ by adjointable
operators, which we formalise as a homomorphism $\phi \colon A \to
\Ll(X)$. Each $C^*$-algebra $A$ is a right-Hilbert $A$--$A$
bimodule  $_AA_A$ with actions given by multiplication and inner product
given by $(a,b) \mapsto a^*b$. The homomorphism that takes $a
\in A$ to left-multiplication by $a$ on $_AA_A$ is an isomorphism
of $A$ onto $\Kk(_AA_A)$.

% The balanced tensor product $X \otimes_A Y$ of right-Hilbert
% $A$--$A$ bimodules $X$ and $Y$ (with left actions $\phi$ and $\rho$) is the completion of the
% vector space spanned by
% elements $x \otimes_A y$ where $x \in X$ and $y \in Y$ subject
% to the relation $x \cdot a \otimes_A y = x \otimes_A \rho(a)y$, in the norm determined by
% the inner-product
% \[
% \langle x_1 \odot y_1, x_2 \odot y_2 \rangle_A = \langle y_1,
% \langle x_1, x_2\rangle_A \cdot y_2\rangle_A.
% \]
% For $S \in \Ll(X)$ and $T \in \Ll(Y)$, the formula $(S \otimes
% T)(x \otimes_A y) = Sx \otimes Ty$ determines an adjointable
% operator on $X \otimes_A Y$. In particular, there is a left action of $A$
% on $X \otimes_A Y$ implemented by the homomorphism $a
%  \mapsto \phi(a) \otimes I_{\mathcal{L}(Y)}$.

The balanced tensor product $X \otimes_A Y$ of right-Hilbert
$A$--$A$ bimodules $X$ and $Y$ (with left actions $\phi$ and $\rho$) is the completion of the
vector space spanned by
elements $x \otimes_A y$ where $x \in X$ and $y \in Y$ subject
to the relation $x \cdot a \otimes_A y = x \otimes_A \rho(a)y$, in the norm determined by
the inner-product
\[
\langle x_1 \otimes_A y_1, x_2 \otimes_A y_2 \rangle_A = \langle y_1,
\langle x_1, x_2\rangle_A \cdot y_2\rangle_A.
\]
There is a right action of $A$ on $X\otimes_A Y$ given by
$(x\otimes_A y)\cdot a=x\otimes_A(y\cdot a)$. With this $X\otimes_A Y$ is
a right-Hilbert $A$-module. For $S \in \Ll(X)$, the formula $(S
\otimes 1_{\Ll(Y)})(x \otimes_A y) = Sx \otimes y$ determines an
adjointable operator on $X \otimes_A Y$. In particular, there
is a left action of $A$ on $X \otimes_A Y$ implemented by the
homomorphism $a
 \mapsto \phi(a) \otimes 1_{\mathcal{L}(Y)}$. With this, $X\otimes_A Y$ is a right-Hilbert
$A$--$A$ bimodule.

\subsection{Semigroups and product systems of Hilbert bimodules}\label{subsec_prod_syst}

A product system over a unital, discrete semigroup $P$ consists
of a semigroup $X$ equipped with a semigroup homomorphism $d
\colon X \to P$ such that $X_p := d^{-1}(p)$ is a right-Hilbert
$A$--$A$ bimodule for each $p\in P$, $X_e = {_A}A_A$, and the
multiplication on $X$ implements isomorphisms $X_p \otimes_A
X_q \cong X_{pq}$ for $p,q \in P \setminus \{e\}$ and the right
and left actions of $X_e = A$ on each $X_p$.
For $p \in P$, we denote the homomorphism of $A$ to
$\Ll(X_p)$ which implements the left action by $\phi_p$. We automatically
 have $\phi_{pq}(a)(xy) =
(\phi_p(a)x)y$ for all $x \in X_p$, $y \in X_q$ and $a \in A$.

Given $p, q \in P$ with $p \not= e$, there is a homomorphism
$\iota^{pq}_p \colon \Ll(X_p) \to \Ll(X_{pq})$ characterised by
\begin{equation}\label{eq:iotapq def}
\iota^{pq}_p(S)(xy) = (Sx)y\text{ for all $x \in X_p$, $y \in
X_{q}$ and $S \in \Ll(X_p)$.}
\end{equation}
Identifying $\Kk(X_e)$ with $A$ as
above, one can also define $\iota^q_e \colon \Kk(X_e)\to
\Ll(X_{q})$ simply by letting $\iota^q_e=\phi_q$ for all $q$,
see \cite[\S 2.2]{SY}.

We will be interested in semigroups arising in quasi-lattice
ordered groups in the sense of Nica \cite{N}. Given a discrete
group $G$ and a subsemigroup $P$ of $G$ such that $P \cap
P^{-1} = \{e\}$, we say that $(G,P)$ is a \emph{quasi-lattice
ordered group} if, under the partial order $g \le h \iff
g^{-1}h \in P$, any two elements $p, q$ in $G$ with a common
upper bound in $P$ have a least common upper bound $p \vee q$
in $P$ (it follows from \cite[Lemma 7]{CL2002} that this
definition is equivalent to Nica's original definition from
\cite{N}, which Fowler also uses in \cite{F99}, and to the
definition Crisp and Laca use in \cite{CL2002} and \cite{CL}).
We write $p \vee q = \infty$ to indicate that $p,q \in G$ have
no common upper bound in $P$, and we write $p \vee q < \infty$
otherwise. As is standard, see \cite{N}, if $p \vee q < \infty$
for all $p,q \in P$, we say that $P$ is \emph{directed}.

Given a quasi-lattice ordered group $(G,P)$, a product system
$X$ over $P$ is called \emph{compactly aligned} (as in
\cite[Definition 5.7]{F99}) if $\iota^{p \vee q}_p(S) \iota^{p
\vee q}_q(T) \in \Kk(X_{p\vee q})$ whenever $S \in \Kk(X_p)$
and $T \in \Kk(X_q)$, and $p \vee q < \infty$. An explanation
is in order here: Fowler only defines compactly aligned in the
case that each $X_p$ is essential as a left $A$-module.
However, since we use $\iota^{p \vee q}_p(S)$ and not
$S\otimes_A 1_{p^{-1}(p\vee q)}$ as in \cite[Definition
5.7]{F99}, and since these make sense also when $p=e$, we can
work with compactly aligned product systems of not necessarily
essential bimodules.

\subsection{Representations of product systems}\label{subsection:reps_prod_syst}
Given a  product system $X$ over $P$, a Toeplitz representation
$\psi$ of $X$ in a $C^*$-algebra $B$ is a map $\psi : X \to B$
such that:
\begin{enumerate}
\item  for each $p$, $\psi_p:=\psi\vert_{X_p} \colon X_p \to B$
    is linear and $\psi_e$ is a homomorphism;
\item $\psi$ takes multiplication in $X$ to multiplication
    in $B$; and
\item $\psi_e(\langle x,y\rangle^p_A) = \psi_p(x)^*
    \psi_p(y)$ for all $x,y \in X_p$ (where $\langle
    x,y\rangle^p_A$ denotes the $A$-valued inner product on
    $X_p$).
\end{enumerate}
In particular, each $\psi_p$ is a Toeplitz representation of
$X_p$ in $B$, see \cite{F99}. A Toeplitz representation $\psi$
of $X$ is \emph{injective} provided that the homomorphism
$\psi_e:X_e \to B$ is injective. Note that property (3) then
implies that $\psi_p$ is an isometry for each $p\in P$. In this
paper, we will frequently drop the word \emph{Toeplitz} and
refer to a map $\psi$ as above simply as a representation of
$X$.

Given a Toeplitz representation $\psi$ of a product system $X$,
there are $*$-homomorphisms $\psi^{(p)} \colon \Kk(X_p)\to B$
such that $\psi^{(p)}(x \otimes y^*)=\psi_p(x) \psi_p(y)^*$ for
all $x,y\in X_p$ (see for example \cite{Pimsner1997}).
Proposition~2.8 of \cite{F99} shows that there is a universal
$C^*$-algebra $\Tt_X$ generated by a universal Toeplitz
representation $i$ of $X$.

Now suppose that $(G,P)$ is a quasi-lattice ordered group and
$X$ is a compactly aligned product system over $P$. We say that
a Toeplitz representation $\psi$ of $X$ is \emph{Nica
covariant} if
\[
\displaystyle \psi^{(p)}(S)\psi^{(q)}(T) =
\begin{cases}
\psi^{(p\vee q)}\big(\iota^{p\vee q}_p(S)\iota^{p\vee q}_q(T)\big)
& \text{if $p\vee q < \infty$} \\
0 &\text{otherwise}
\end{cases}
\]
for all $S \in \Kk(X_p)$ and $T \in \Kk(X_q)$ (see also
\cite[Definition 5.7]{F99}). Let $\Tc(X)$ be the quotient of
$\Tt_X$ by the ideal generated by the elements
$$
i^{(p)}(S)i^{(q)}(T) - i^{(p \vee q)}(\iota^{p \vee q}_p(S)
\iota^{p \vee q}_q(T))
$$ where $p,q \in P$, $S \in \Kk(X_p)$, $T
\in \Kk(X_q)$, and by convention, $\iota^{p \vee q}_p(S)
\iota^{p \vee q}_q(T) = 0$ if $p \vee q = \infty$. The
composition  of  the quotient map from $\Tt_X$ onto $\Tc(X)$
with $i$ is a Nica covariant Toeplitz representation $i_X
\colon X\to \Tc(X)$ with the following universal property: if
$\psi$ is a Nica covariant Toeplitz representation of $X$ in
$B$ there is a $*$-homomorphism $\psi_* : \Tc(X) \to B$ such
that $\psi_*\circ i_X=\psi$. Thus if $X$ is a compactly aligned
product system of essential Hilbert bimodules, then $\Tc(X)$
coincides with Fowler's algebra (denoted by the same symbol)
from \cite{F99} defined for not necessarily compactly aligned
product systems of essential Hilbert bimodules. By an argument
similar to \cite[Theorem~6.3]{F99} we have
\begin{equation}\label{eq:Tc spanners}
\Tc(X)=\clsp\{\, i_X(x)i_X(y)^* \mid x, y \in X \,\}.
\end{equation}
It follows from \eqref{eq:Tc spanners} that if the image of a Nica covariant
Toeplitz representation $\psi$ of $X$ generates $B$ as a $C^*$-algebra, then
$B= \clsp\{\, \psi(x)\psi(y)^* \mid x,y \in X \,\}$.

\subsection{The algebra $\NO{X}$}\label{subsection_def_NOX}
To define Cuntz-Pimsner covariance of representations, we must
first summarise some definitions from \cite[Section~3]{SY}. We
say that a predicate statement $\mathcal{P}(s)$ (where $s \in
P$) is true for large $s$ if for every $q \in P$ there exists
$r \ge q$ such that $\mathcal{P}(s)$ is true whenever $r \le
s$.

Assume $(G,P)$ is quasi-lattice ordered and $X$ is a compactly
aligned product system over $P$. Define $I_e = A$, and for each
$q \in P \setminus \{e\}$ write $I_q := \bigcap_{e < p \le q}
\ker(\phi_p)$. We then write $\widetilde{X}_q$ for the
right-Hilbert $A$--$A$ bimodule
\[
\widetilde{X}_q := \textstyle{\bigoplus}_{p \le q} X_p \cdot I_{p^{-1}q}.
\]
The homomorphism implementing the left action is denoted
$\tilde\phi_q$. We say that $X$ is
\emph{$\tilde\phi$-injective} if the homomorphisms
$\tilde\phi_q$ are all injective.

For $p \not\le q \in P$ we define $\iota^q_p(T) = 0_{\Ll(X_q)}$
for all $T \in \Ll(X_p)$. Recalling the definitions of the maps
$\iota^{pq}_p$ from Section~\ref{subsec_prod_syst}, we then
have homomorphisms $\tilde\iota^q_p \colon \Ll(X_p) \to
\Ll(\widetilde{X}_q)$ for all $p,q \in P$ with $p \not= e$
defined by $\tilde\iota^q_p(T) = \bigoplus_{r \le q}
\iota^r_p(T)$ for all $p,q\in P$ with $p\neq e$. When $p=e$,
similar to the above there is a homomorphism $\tilde\iota^q_e :
\Kk(X_e) \to \Ll(\widetilde{X}_q)$.

Suppose that $X$ is $\tilde\phi$-injective. We say that a Nica covariant
Toeplitz representation $\psi$ of $X$ in a $C^*$-algebra $B$ is
\emph{Cuntz-Nica-Pimsner covariant} (or \emph{CNP-covariant})
if it has the following property:
\[\parbox{0.9\textwidth}{
$\sum_{p \in F} \psi^{(p)}(T_p) = 0_B$ whenever $F \subset P$
is finite, $T_p \in \Kk(X_p)$ for each $p \in F$, and
$\sum_{p \in F} \tilde\iota^q_p(T_p) = 0$ for large $q$.
}\]

As in \cite[Proposition 3.12]{SY}, if $X$ is
$\tilde\phi$-injective, we write $\NO{X}$ for the universal
$C^*$-algebra generated by a CNP-covariant representation $j_X$
of $X$, and call it the Cuntz-Nica-Pimsner algebra of $X$. By
\cite[Remark~4.2]{SY}, the hypothesis that $X$ is
$\tilde\phi$-injective ensures that $j_X$ is an injective
representation. By \cite[Theorem~4.1]{SY}, $X$ is
$\tilde\phi$-injective (and hence $\NO{X}$ is defined and $j_X$
is an injective representation) whenever each $\phi_p$ is injective, and also
whenever each bounded subset of $P$ has a maximal element.

\section{Analysis of the core}\label{sec:core}

In this section we lay the foundation for the proof of our main
result Theorem \ref{thm:projective property}. To do this we
shall analyse the fixed-point algebra of $\Tc(X)$ under a
canonical coaction $\delta$. As a
corollary we show that
% Added 22/2
under certain conditions
$\NO{X}$ satisfies criterion (B) of \cite[Section~1]{SY}. Throughout the
rest of the article, we write $q_{\CNP} \colon \Tc(X) \to \NO{X}$ for the canonical surjection
arising from the universal property of $\Tc(X)$.

\begin{lemma}\label{lem:psi equal}
Let $(G,P)$ be a quasi-lattice ordered group and let $X$ be a
product system over $P$ of right-Hilbert $A$--$A$ bimodules.
Let $\psi \colon X \to B$ be a Toeplitz representation of $X$.
Then:
\begin{enumerate}
\item If $p \le t \in P$, $T \in
\Kk(X_p)$, and $x \in X_t$, then
$\psi_t(\iota_{p}^{t}(T)(x)) = \psi^{(p)}(T) \psi_t(x)$;
\item If $t < r \le s \in P$, $T \in
\Kk(X_r)$, and $x \cdot a \in X_t \cdot I_{t^{-1}s}$, then
$\psi^{(r)} (T) \psi_t (x \cdot a) = 0$.
\end{enumerate}
\end{lemma}

\begin{proof}
(1) If $p=e$ then~(1) follows from the observations that
$\Kk(X_e) \cong A$ and $\iota_{e}^{t} := \phi_{t}$, so suppose
$p \ne e$. Since $\lsp \{\, xy \mid x \in X_p,\, y \in X_{p^{-1}t} \,\}$
is dense in $X_t$, and since $\lsp\{\, w \otimes z^* \mid w,z \in X_p \,\}$ is
dense in $\Kk(X_p)$, to prove~(1) it suffices to show that for
$x,w,z \in X_p$ and $y \in X_{p^{-1}t}$ we have
\[
\psi_t (\iota_{p}^{t} (w \otimes z^*) (xy))
= \psi^{(p)} (w \otimes z^*) \psi_t (xy).
\]
Using~\eqref{eq:iotapq def}, we calculate:
\begin{align*}
\psi_t (\iota_{p}^{t} (w \otimes z^*) (xy))
&= \psi_p (w \cdot \langle z, x \rangle_{A}^{p}) \psi_{p^{-1}t} (y)
\\&= \psi_p (w) \psi_e (\langle z, x \rangle_{A}^{p}) \psi_{p^{-1}t} (y)
\\&= \psi_p (w) \psi_{p}^* (z) \psi_p (x) \psi_{p^{-1}t} (y)
\\&= \psi^{(p)} (w \otimes z^*) \psi_{t} (xy)
\end{align*}
as required in (1).

(2) If $t=e$ then $x \cdot a \in I_{s}$, so $x \cdot a \in
\ker (\phi_{r})$. By using that $\Kk(X_r)=\clsp\{y\otimes z^*: x,y\in X_r\}$ one easily
checks that $\psi_e(b)\psi^{(r)}(S)=\psi^{(r)}(\phi_r(b)S)$ for $b\in X_r$ and $S\in\Kk(X_r)$.
By taking adjoints and letting $b=(x\cdot a)^*$ and $S=T^*$, it follows that
\[
\psi^{(r)} (T) \psi_{e} (x \cdot a) = \psi^{(r)} (T \phi_r (x \cdot a)) = 0.
\]
Now suppose $t \ne e$. Fix $y \in X_t$, $z \in X_{t^{-1}r}$,
and $v \in X_r$. It suffices to show that $\psi^{(r)} (v
\otimes (yz)^*) \psi_t (x \cdot a) =
0$. Since $a \in I_{t^{-1}s} = \bigcap_{e < q \le t^{-1}s} \ker
(\phi_{q})$, we have $\phi_{t^{-1}r} (a) = 0$, and hence
\begin{align*}
\psi^{(r)} (v \otimes (yz)^*) \psi_t (x \cdot a)
&= \psi_r (v) \psi_{t^{-1}r} (z)^* \psi_t (y)^*  \psi_t (x \cdot a)
\\&= \psi_r (v) \psi_{t^{-1}r} (z)^* \psi_e (\langle y,x \cdot a \rangle_{A}^{t})
\\&= \psi_r (v) \psi_{t^{-1}r}(z)^* \psi_{e} (\langle y ,x \rangle_{A}^{t}a)
\\&= \psi_r (v) (\psi_e (a^* \langle x,y \rangle_{A}^{t}) \psi_{t^{-1}r} (z))^*
\\&= \psi_r (v) \psi_{t^{-1}r} (\phi_{t^{-1}r} (a^* \langle x,y \rangle_{A}^{t})z)^*
\\&= \psi_r (v) \psi_{t^{-1}r} (\phi_{t^{-1}r} (a)^* \phi_{t^{-1}r} (\langle x,y \rangle_{A}^{t})z)^*
\\&=0.\qedhere
\end{align*}
\end{proof}

Lemma~\ref{lem:psi equal}(2) says, roughly, if $r \in
tP\setminus\{t\}$, then $\psi^{(r)}(T) \in B$ annihilates
$\psi_{t}(X_{t} \cdot I_{t^{-1}s})$ whenever $s \in rP$. The
next corollary says that when $X$ is compactly aligned and
$\psi$ is Nica covariant, we can replace the requirement that
$r \in tP\setminus\{t\}$ with the much weaker requirement that
$t \not\in rP$, and that $s$ is a common upper bound for $t$
and $r$.

\begin{cor}\label{lem:psi zero 2}
Let $(G, P)$ be a quasi-lattice ordered group and let $X$ be a
compactly aligned product system over $P$ of right-Hilbert
$A$--$A$ bimodules. Let $\psi \colon X \to B$ be a Nica
covariant representation of $X$. Suppose $p,t \le s
\in P$ and $p \not\le t$. Then for $T \in \Kk(X_p)$ and $x
\cdot a \in X_t \cdot I_{t^{-1}s}$, we have $\psi^{(p)} (T)
\psi_t (x \cdot a) = 0$.
\end{cor}
\begin{proof}
Let $(E_k)_{k \in K}$ be an approximate identity for $\Kk(X_t
\cdot I_{t^{-1}s})$. Since $p,t \le s$, we have $p \vee t <
\infty$. Hence Nica covariance and the fact that each $E_k\in
\Kk(X_t)$ imply that
\begin{align*}
\psi^{(p)}(T)\psi_t (x \cdot a)
 &= \lim_{k \in K} \psi^{(p)}(T)\psi_t (E_k(x \cdot a)) \\
 &= \lim_{k \in K} \psi^{(p)}(T) \psi^{(t)}(E_k) \psi_t(x \cdot a) \\
 &=\lim_{k \in K} \psi^{(p \vee t)}\big(\iota^{p \vee t}_p(T)\iota^{p \vee t}_t(E_k)\big)
 \psi_t(x \cdot a).\\
\end{align*}
Since $p \not\le t$ forces $t <  p \vee t$ and since $p \vee t \le
s$, the result now follows from statement~(2) of
Lemma~\ref{lem:psi equal}.
\end{proof}

\begin{lemma}\label{prp:ev zero}
Let $(G, P)$ be a quasi-lattice ordered group and let $X$ be a
compactly aligned product system over $P$ of right-Hilbert
$A$--$A$ bimodules. Suppose either that the left action on each
fibre is by injective homomorphisms, or that $P$ is directed.
Let $\psi \colon X \to B$ be an injective Nica covariant
representation of $X$. Fix a finite subset $F \subset
P$ and fix operators $T_p \in \Kk(X_p)$ for each $p\in F$
satisfying $\sum_{p \in F} \psi^{(p)}(T_p) = 0$. Then $\sum_{p
\in F} \iotale^s_p(T_p) = 0$ for large $s$.
% in the sense of \cite[Definition~3.8]{SY}.
\end{lemma}

\begin{proof}
Fix $q \in P$. We must show that there exists $r \ge q$ such
that for every $s \ge r$, we have $\sum_{p \in F} \tilde\iota_p^s(T_p) = 0_{\Ll(\widetilde{X}_s)}$.

List the elements of $F$ as $p_1,
\dots, p_{|F|}$. Define $r_0 := q$, and inductively, for $1 \le
i \le |F|$, define
\[
r_i := \begin{cases}
r_{i-1} \vee p_i &\text{if $r_{i-1} \vee p_i <\infty$}\\
r_{i-1} &\text{otherwise.}
\end{cases}
\]
Set $r := r_{|F|}$, and note that this satisfies $r \ge q$.
With no extra assumptions on the quasi-lattice ordered group we
also have $r \ge p$ whenever $p \in F$ satisfies $r \vee p<
\infty$. If $P$ is directed then $r = q \vee \big(\bigvee_{p\in
F}p\big)$, and is an upper bound for $F$.

Let $s \ge r$. To show that $\sum_{p \in F} \tilde
\iota_{p}^{s} (T_p) = \sum_{p \in F} \big(\bigoplus_{t \le s}
\iota_{p}^{t} (T_p)\big)$ is equal to the zero operator on
$\widetilde X_s = \bigoplus_{t \le s} {X_{t} \cdot
I_{t^{-1}s}}$ we shall prove that $\sum_{p \in F, p \le t}
\iota_{p}^{t} (T_p) |_{{X_{t} \cdot I_{t^{-1}s}}} = 0_{\Ll
({X_{t} \cdot I_{t^{-1}s}})}$ for each $t \le s$. Indeed, for
$x \cdot a \in X_{t} \cdot I_{t^{-1}s}$, using
Lemma~\ref{lem:psi equal}(1) we have
\begin{align}
\textstyle{\psi_t \big(\sum_{p \in F, p \le t} \iota_{p}^{t} (T_p) (x \cdot a)\big)}
 &= \textstyle{\sum_{p \in F, p \le t} \psi_t (\iota_{p}^{t} (T_p) (x \cdot a))} \notag\\
 &= \textstyle{\sum_{p \in F, p \le t} \psi^{(p)} (T_p) \psi_t (x \cdot a).} \label{eq:tricky_sum}
 \end{align}

We claim that \eqref{eq:tricky_sum} is equal to $\sum_{p \in F} \psi^{(p)} (T_p) \psi_t (x \cdot a)$.
We will establish this claim under each of the additional hypotheses of the lemma. Note that the claim
comes down to proving
\begin{equation}\label{eq:zero_term}
\psi^{(p)}(T_p)\psi_t(x\cdot a)=0\text{ if }p\in F, p\not\le t.
\end{equation}

Suppose that the $\phi_p$ are injective. Then $I_{t^{-1}s}=0$
for $t<s$, and so $a=0$ in \eqref{eq:zero_term} unless $t=s$.
Thus it suffices in this case to show that
$\psi^{(p)}(T_p)\psi_s(x)=0$ for $T_p\in\Kk(X_p)$ and $x\in
X_s$. By choice of $r$ and $s$ and the assumption $p\not\le s$,
we necessarily have $p\vee s=\infty$. Let $(E_k)_{k \in K}$ be
an approximate identity for $\Kk(X_s)$. By Nica covariance,
$\psi^{(p)}(T_p)\psi_s(x)=\lim_{k \in K} \psi^{(p)}(T_p)
\psi^{(s)}(E_k) \psi_s(x)=0$.

Now suppose that $P$ is directed. Then $p\le r\le s$ and
equation \eqref{eq:zero_term} follows from Corollary~\ref{lem:psi zero
2}.

Thus we have in both cases that $\psi_t(\sum_{p\in F,p\le
t}\iota_p^t(T_p)(x\cdot a))=\sum_{p\in F}\psi^{(p)}(T_p)\psi_t(x\cdot
a)$. Since this last sum is equal to $0$ by hypothesis, and since the representation $\psi$ is injective,
so that in particular every $\psi_t$ is injective, it follows that  $\sum_{p \in F, p \le t}
\iota_{p}^{t} (T_p) (x \cdot a) = 0$, as needed.
\end{proof}

\begin{rmk}
The hypotheses that either the left actions are all injective,
or $P$ is directed are genuinely necessary in Lemma~\ref{prp:ev
zero}; see Example~\ref{eg:degenerate}.
\end{rmk}

%\subsection{The core}\label{subsec:core}

Let $(G, P)$ be a quasi-lattice ordered group and let $X$ be a
compactly aligned product system over $P$ of right-Hilbert
$A$--$A$ bimodules. It follows (see \cite[Proposition
5.10]{F99}) from the Nica-covariance of $i_X$ that
\begin{equation}\label{eq:F-def}
\Ff := \clsp\{\, i_X(x) i_X(y)^* \mid x,y\in X, d(x) = d(y) \,\}
\end{equation}
is closed under multiplication, and thus that it is a $C^*$-subalgebra of $\Tc(X)$.
We call this subalgebra the \emph{core} of $\Tc(X)$.

For any discrete group $G$ there is a homomorphism $\delta_G
\colon C^*(G)\to C^*(G) \otimes C^*(G)$ given by
$\delta_G(g)=i_G(g)\otimes i_G(g)$. Recall that a full coaction
of $G$ on a $C^*$-algebra $A$ is an injective homomorphism
$\delta \colon A\to A\otimes C^*(G)$ which is nondegenerate (in
the sense that  $\clsp\delta(A)(A\otimes C^*(G))=A\otimes
C^*(G)$) and satisfies the coaction identity $(\delta\otimes
\id_{C^*(G)})\circ \delta= (\id_A \otimes \delta_G)\circ
\delta$ (see, for example, \cite{qui:discrete coactions}. All
coactions in this paper are full. The generalised fixed-point
algebra of $A$ with respect to $\delta$ is $A_e^\delta:=\{a\in
A\mid \delta(a)=a\otimes i_G(e)\}$.

We will now show that there is a coaction of $G$ on $\Tc(X)$
whose generalised fixed-point algebra is equal to the core
$\Ff$. For Fowler's $\Tc(X)$ associated to a not-necessarily
compactly aligned product system over $P$ of essential $A$--$A$
bimodules (where $(G, P)$ is quasi-lattice ordered),
Proposition 4.7 in \cite{F99} and the discussion
 preceding \cite[Theorem 6.3]{F99} imply the existence of a coaction with similar properties as
 $\delta$ in the next result. We present a different and more direct proof here.

\begin{prop}\label{prop:existence_coaction}
  Let $(G, P)$ be a quasi-lattice ordered group and let $X$ be a
  compactly aligned product system over $P$ of right-Hilbert
  $A$--$A$ bimodules. Then there is a coaction $\delta$ of $G$ on
  $\Tc(X)$ such that $\delta(i_X(x))=i_X(x)\otimes i_G(d(x))$ for all $x\in X$.
\end{prop}

\begin{proof}
  Let $\psi: X\to \Tc(X)\otimes C^*(G)$ be the map $x\mapsto i_X(x)\otimes i_G(d(x))$.
  It is straightforward to check that $\psi$ is a Nica covariant  representation
  of $X$. It follows from the universal property of $\Tc(X)$ that there is a
  $*$-homomorphism $\delta:\Tc(X)\to \Tc(X)\otimes C^*(G)$ such that
  $\delta(i_X(x))=\psi(x)=i_X(x)\otimes i_G(d(x))$ for all $x\in X$. We will show that
  $\delta$ is a coaction.

  We first show that $\delta$ is nondegenerate.
  Let $(\theta_\lambda)_{\lambda\in\Lambda}$ be an approximate identity for $\Ff$. We claim
  that $(\theta_\lambda)_{\lambda\in\Lambda}$ is also an approximate identity for $\Tc(X)$.
  Since $\Tc(X)$ is the closure of the span of elements of the form
  $i_X(x)i_X(y)^*$, it suffices to show that $\theta_\lambda i_X(x)i_X(y)^*\to i_X(x)i_X(y)^*$
  for all $x,y\in X$. Fix $x,y\in X$ and let $p=d(x)\in P$.
  By \cite[Proposition 2.31]{TFB} we may write
  $x=z\cdot\langle z,z\rangle_A^p=(z\otimes z^*)(z)$ for some $z\in X$, and then
  $i_X(x)i_X(y)^*=i_X^{(p)}(z\otimes z^*)i_X(z)i_X(y)^*$. Since $i_X^{(p)}(z\otimes z^*)\in \Ff$,
  we have $\theta_\lambda i_X^{(p)}(z\otimes z^*)\to i_X^{(p)}(z\otimes z^*)$, and hence
  $\theta_\lambda i_X(x)i_X(y)^*\to i_X(x)i_X(y)^*$ as claimed. Since
  $\delta(\theta_\lambda)=\theta_\lambda\otimes 1$ for each
  $\lambda\in\Lambda$, the approximate identity $(\theta_\lambda)_{\lambda\in\Lambda}$ is
  mapped under $\delta$ to an approximate identity for $\Tc(X)\otimes C^*(G)$, and it
  follows that $\clsp\delta(\Tc(X))(\Tc(X)\otimes C^*(G))=\Tc(X)\otimes C^*(G)$.
  By checking on generators, it is easy to see that $\delta$ satisfies the coaction identity
  $(\delta\otimes \id_{C^*(G)})\circ \delta= (\id_{\Tc(X)} \otimes \delta_G)\circ \delta$,
  and $\delta$ is injective since
  $\id_{\Tc(X)}=(\id_{\Tc(X)}\otimes\epsilon)\circ\delta$ where
  $\epsilon:C^*(G)\to\CC$ is the integrated form of the representation $g\mapsto 1$.
\end{proof}

We call the above coaction $\delta$ of $G$ on  $\Tc(X)$ for the \emph{gauge coaction} on
$\Tc(X)$.
 It follows from equation \eqref{eq:Tc spanners} that the generalised fixed-point algebra $\Tc(X)_e^\delta = \{\, a \in
\Tc(X) \mid \delta(a) = a \otimes i_G(e) \,\}$ is equal to the core $\Ff$.

Since $i^{(p)}_X \colon \Kk(X_p)\to \Tt_X$ satisfies $i^{(p)}_X(x \otimes
y^*)=i_X(x) i_X(y)^*$ and each $\Kk(X_p) = \clsp\{\, x \otimes y^*
\mid x,y \in X_p \,\}$ by definition, we have
\begin{equation}\label{eq:F-spanning}
\Ff = \clsp\{\, i^{(p)}_X(T) \mid p \in P \text{ and } T \in \Kk(X_p) \,\}.
\end{equation}

 We say that a
subset $F$ of $P$ is \emph{$\vee$-closed} if, whenever $p,q \in
F$ satisfy $p \vee q < \infty$, we have $p \vee q \in F$. Let
$\fvcl{P}$ denote the set of finite $\vee$-closed subsets of
$P$; then $\fvcl{P}$ is directed under set inclusion
(see~\cite[p. 367]{F99}). If $F\in \fvcl{P}$ is bounded, then
$\bigvee_{p\in F}p$ is a maximal element in $F$.

For $p \in P$, we write $B_p$ for
the $C^*$-subalgebra $i_X^{(p)}(\Kk(X_p)) \subset \Tc(X)$. For
each finite $\vee$-closed subset $F$ of $P$, we denote by $B_F$
the linear subspace
\begin{equation}\label{eq:Bp-def}
B_F := \textstyle{\sum}_{p \in F} B_p = \textstyle{\Big\{\, \sum_{p \in F} i_X^{(p)}(T_p) \mid \mbox{$T_p \in \Kk(X_p)$ for each $p \in F$} \,\Big\}} \subset \Tc(X).
\end{equation}
Equation~\eqref{eq:F-spanning} implies that
\begin{equation}\label{eq:BPs span F}
\Ff = \textstyle{\overline{\bigcup_{F \in \fvcl{P}} B_{F}}}.
\end{equation}

\begin{lemma}\label{lem:core structure}
Let $(G, P)$ be a quasi-lattice ordered group and let $X$ be a
compactly aligned product system over $P$ of right-Hilbert
$A$--$A$ bimodules.  For each finite $\vee$-closed subset $F$
of $P$, the space $B_F$ is a $C^*$-subalgebra of $\Ff$.
\end{lemma}
\begin{proof}
Fix a finite $\vee$-closed subset $F$ of $P$. Then $B_F$ is a
subspace of $\Ff$ by definition. One can check on spanning
elements that it is closed under adjoints and multiplication
(for the latter, one uses the Nica covariance of the universal
representation $i_X$ of $X$ in $\Tc(X)$). It therefore suffices
to show that $B_F$ is norm-closed.

We proceed by induction on $|F|$. If $|F| = 1$, then $F = \{p\}$ for
some $p \in P$, and then $B_F = B_p = i_X^{(p)}(\Kk(X_p))$ is the
range of a $C^*$-homomorphism and hence closed.

Now suppose that $B_F$ is closed whenever $|F| \le k$. Suppose
that $F \subset P$ is $\vee$-closed with $|F| = k+1$. Since $F$
is finite, we may fix an element $m$ of $F$ which is minimal in
the sense that for $p \in F\setminus\{m\}$, we have $p \not\le
m$. The sets $\{m\}$ and $F \setminus \{m\}$ are both finite
and $\vee$-closed, and it follows from our induction hypothesis that $B_m$
and $B_{F\setminus\{m\}}$ are $C^*$-subalgebras of $\mathcal{F}$.
For $p \in F \setminus \{m\}$, we have $p
\not\le m$ by choice of $m$ and it follows that if $p \vee m <
\infty$, then $p \vee m \in F \setminus \{m\}$. Hence for $S
\in \Kk(X_p)$ and $T \in \Kk(X_m)$, we have
\[
i_X^{(p)}(S) i_X^{(m)}(T)
 = i_X^{(p \vee m)}(\iota^{p \vee m}_p(S) \iota^{p \vee m}_m(T))
 \in i_X^{(p \vee m)}(\Kk(X_{p \vee m})) \subset B_{F \setminus \{m\}}.
\]
Similarly, $i_X^{(m)}(T) i_X^{(p)}(S) \in B_{F \setminus
\{m\}}$, so by linearity, $ab, ba \in B_{F \setminus \{m\}}$
for all $a \in B_{F \setminus\{m\}}$ and $b \in B_m$.
Corollary~1.8.4 of \cite{Dix:C*-algebras} now shows that $B_F =
B_m + B_{F \setminus \{m\}}$ is norm closed.
\end{proof}

The following proposition is the key technical result which we will use in
the proof of our main theorem in the next section.

\begin{prop} \label{lemma:inj}
Let $(G, P)$  be a quasi-lattice ordered group and let $X$ be a
compactly aligned product system over $P$ of right-Hilbert
$A$--$A$ bimodules. Suppose either that the left action on each
fibre is by injective homomorphisms, or that $P$ is directed.
Let $\psi \colon X \to B$ be an injective Nica covariant
representation of $X$ and let $\psi_* : \Tc(X) \to B$
be the homomorphism characterised by $\psi = \psi_* \circ i_X$.
Then $\ker(\psi_*)\cap\Ff\subset \ker (q_{\CNP})$.
\end{prop}

\begin{proof}
By
\cite[Lemma~1.3]{ALNR}, equation \eqref{eq:BPs span F} and Lemma
\ref{lem:core structure}, it suffices to show that $\ker(\psi_*)
\cap B_F \subset \ker(q_{\CNP})$ for each $F \in \fvcl{P}$. For
this, we fix $F \in \fvcl{P}$ and generalised compact operators
$T_p \in \Kk(X_p)$ for $p \in F$, so that $c := \sum_{p \in F}
i_X^{(p)}(T_p)$ is a typical element of $B_F$. Suppose that $c
\in \ker(\psi_*)$; we must show that $c \in \ker(q_{\CNP})$ as
well. Since the representation $\psi$ is injective, Lemma~\ref{prp:ev zero}
implies that $\sum_{p \in F} \iotale^s_p(T_p) = 0$ for large
$s$ in the sense of \cite[Definition~3.8]{SY}. Since $j_X$ is
CNP-covariant, it follows that
\[
\textstyle{q_{\CNP}(c) = \sum_{p \in F} j_X^{(p)}(T_p) = 0}
\]
as well, so $c \in \ker(q_{\CNP})$ as required.
\end{proof}

We now have enough machinery to confirm that $\NO{X}$ indeed
satifies criterion (B) of \cite[Section 1]{SY} when the left
actions on the fibres of $X$ are all injective, or $P$ is
directed. One could use the following theorem to prove directly
a gauge-invariant uniqueness theorem for $\NO{X}$ when $G$ is
amenable, but since this will be an easy corollary of our more
general main result, we will not pursue this line of attack.
Recall from \cite{SY} that $\NO{X}$ has the following universal
property: for each CNP-covariant representation $\psi$ of $X$
there is a homomorphism $\intfrm{\psi}$ such that
$\intfrm{\psi} \circ j_X = \psi$.

\begin{theorem}\label{thm:inj on core}
Let $(G,P)$ be a quasi-lattice ordered group and let $X$ be a
compactly aligned product system over $P$ of right-Hilbert
$A$--$A$ bimodules. Assume either that the left actions on the
fibres of $X$ are all injective, or that $P$ is directed and
$X$ is $\tilde\phi$-injective. Let $\psi \colon X \to B$ be a
CNP-covariant representation of $X$ in a $C^*$-algebra $B$.
Then the induced homomorphism  $\intfrm{\psi} : \NO{X} \to B$
is injective on $q_{\CNP}(\Ff)$ if and only if $\psi$ is
injective as a Toeplitz representation.
\end{theorem}

\begin{proof}
Suppose that $\intfrm{\psi}$ is injective on $q_{\CNP}(\Ff)$.
By \cite[Theorem4.1]{SY}, $j_X$ is injective on
$A$. Hence $\psi_e = \intfrm{\psi}\circ(j_X)_e$ is also
injective, and thus $\psi$ is an injective Toeplitz representation.

Now suppose that $\psi$ is injective as a Toeplitz representation; we must show that
$\intfrm{\psi}$ is injective on $q_{\CNP}(\Ff)$. By definition
of $\NO{X}$ and of $\Pi\psi$, we have $\Pi\psi \circ q_{\CNP} =
\psi_*$. Proposition~\ref{lemma:inj} therefore implies that
$\ker(\Pi\psi \circ q_{\CNP}) \cap \Ff \subset \ker(q_{\CNP})$.
Hence $\Pi\psi$ is injective on $q_{\CNP}(\Ff)$ as claimed.
\end{proof}

\begin{example}\label{eg:degenerate}
We present an example of a product system $X$ in which the left
actions are not injective, and $P$ is not directed, and the
conclusion of Lemma~\ref{prp:ev zero} fails. It is easy to see
that the conclusions of Proposition~\ref{lemma:inj} and
Theorem~\ref{thm:inj on core} both fail in this example (see
also Remark~\ref{rmk:troublesomeex}).

Let the quasi-lattice ordered group be $(G,P)=(\FF_2,\FF_2^+)$,
and denote by $a$ and $b$ the generators of $\FF_2^+$. Define a
product system over $\FF_2^+$ by $X_{a^n}=\CC$ for $n\in \NN$
and $X_p=0$ for all other elements of $\FF_2^+$. This is
compactly aligned since $\mathcal{L}(X_p)=\mathcal{K}(X_p)$ for each $p$,
but the left actions are not all
injective (and $a\vee b=\infty$). Define $\psi \colon X\to \CC$
by $\psi_p(x)=x$ for $x\in X_p$ and $p\in \FF_2^+$. Then $\psi$
is an injective Nica covariant Toeplitz representation of $X$.
Let $1_p$ be the identity in $\Ll(X_p)$ for each $p$, and note
that $1_e\in \Kk(X_e)$ and $1_a\in \Kk(X_a)$. We have that
$\psi^{(e)}(1_e)=\psi^{(a)}(1_a) =1$, so
$\psi^{(e)}(1_e)-\psi^{(a)}(1_a)=0$. However, we claim that
$\iotale^e_s(1_e)-\iotale^a_s(1_a)$ is not equal to $0$ for
large $s$. Indeed, note that
$$
I_q=\begin{cases}0&\text{ if }q=a^n\text{ for some }n\in \NN\\
\CC & \text{ otherwise}
\end{cases}
$$
for $q\in P\setminus\{e\}$. It follows that if $q\geq b$, then
$X_e\cdot I_q=X_e$, and so
$$
(\iota_e^e(1_e)-\iota_a^e(1_a))\vert_{X_e\cdot I_q}=\iota_e^e(1_e)-\iota_a^e(1_a)=
\iota_e^e(1_e)\not= 0,
$$
which shows that $\iotale^e_s(1_e)-\iotale^a_s(1_a)\not= 0$ for
all $s \ge b$.
\end{example}

\section{The co-universal $C^*$-algebra and the uniqueness
theorems}\label{sec:couniversal}

We begin this section with our main theorem. Before stating it,
we introduce some terminology: given
a quasi-lattice ordered group $(G, P)$ and a product system $X$
over $P$ of right-Hilbert $A$--$A$ bimodules, a Toeplitz representation
$\psi \colon X \to B$ is \emph{gauge-compatible} if there is a coaction
$\beta$ of $G$ on $B$ such that
\begin{equation}\label{E:gauge-compatible}
\mbox{$\beta(\psi(x)) = \psi(x) \otimes i_G(d(x))$\quad for all $x \in X$.}
\end{equation}

Suppose that $\psi_1:X\to B_1$ and $\psi_2:X\to B_2$ are two
gauge-compatible Toeplitz representations of $X$, that
$\beta_i$ is a coaction of $G$ on $B_i$ satisfying
$\beta_i(\psi_i(x))=\psi_i(x)\otimes i_G(d(x))$ for all $x\in
X$ and $i=1,2$, and that $\phi :B_1\to B_2$ is a
$*$-homomorphism satisfying $\phi\circ \psi_1=\psi_2$. Then
$\phi$ is equivariant for $\beta_1$ and $\beta_2$, meaning that
$(\phi\otimes \id_{C^*(G)})\circ \beta_1=\beta_2\circ \phi$.

Since our main result depends on the technical hypothesis that
$X$ is $\tilde\phi$-injective, we emphasise that the results of
\cite{SY} imply that this is automatic whenever either the left
actions on the fibres of $X$ are all injective or every bounded
subset of $P$ has a maximal element.

\begin{theorem}\label{thm:projective property}
Let $(G,P)$ be a quasi-lattice ordered group and  $X$ a
compactly aligned product system over $P$ of right-Hilbert
$A$--$A$ bimodules. Suppose either that the left action on each
fibre is injective, or that $P$ is directed and $X$
is $\tilde\phi$-injective. Then there exists a triple
$(\NO{X}^\reduced, j_X^\reduced, \redgaug)$ which is
co-universal for gauge-compatible injective Nica covariant
representations of $X$ in the following sense:
\begin{enumerate}
\item $\NO{X}^\reduced$ is a $C^*$-algebra, $j^\reduced_X$ is an injective
    Nica covariant representation of $X$ whose image
    generates $\NO{X}^\reduced$, and $\redgaug$ is a coaction of $G$ on $\NO{X}^\reduced$
    such that $\redgaug(j^\reduced_X(x))=j^\reduced_X(x)\otimes i_G(d(x))$ for all $x\in X$.
\item\label{it:c-u property} If $\psi \colon X \to B$ is an
    injective Nica covariant gauge-compatible representation whose image
    generates $B$ then there is a
    surjective $*$-homomorphism $\phi \colon B \to
    \NO{X}^\reduced$ such that $\phi(\psi(x)) =
    j^\reduced_X(x)$ for all $x \in X$.
\end{enumerate}
Moreover, the representation $j^\reduced_X$ is CNP-covariant, the coaction $\redgaug$ is normal, and
$(\NO{X}^\reduced, j^\reduced_X, \redgaug)$ is the unique
triple satisfying \textnormal{(1)} and \textnormal{(2)}: if $(C, \rho, \gamma)$ satisfies the same
two conditions, then there is an isomorphism $\phi \colon C\to
\NO{X}^\reduced$ such that $j^\reduced_X = \phi\circ \rho$ and $\phi$ is equivariant for $\gamma$ and $\redgaug$.
\end{theorem}

\begin{rmk}\label{rmk:troublesomeex}
Although Example~\ref{eg:degenerate} does not satisfy the
assumptions of Theorem~\ref{thm:inj on core}, it nevertheless
does admit a co-universal algebra as described in
Theorem~\ref{thm:projective property}; but this co-universal
algebra is a proper quotient of the algebra $\NO{X}^\reduced$
that we shall construct later. Specifically,
%the co-universal
%algebra for the system $X$ described in
%Example~\ref{eg:degenerate} is the co-universal $C^*$-algebra
%generated by a partial isometry $i_X(1_a)$ and carrying a
%circle action satisfying $\gamma_z(i_X(1_a)) = zi_X(1_a)$, namely $C(\TT)$. However,
it is not difficult to see that every Toeplitz representation of the system $X$
described in Example~\ref{eg:degenerate} is automatically Nica covariant, and that
there is a bijective correspondence between Toeplitz representations
of $X$ and Toeplitz representations of $X_a$ which takes injective
representations to injective representations and gauge-compatible
representations to gauge-compatible representations. It thus follows
that both the Toeplitz algebra and the covariant Toeplitz algebra of
$X$ are equal to the classical Toeplitz algebra $\mathcal{T}$ (generated by a
single isometry), and
that $C(\TT)$ has the co-universal property described in Theorem~\ref{thm:projective property}
with respect to the system $X$. Moreover, one can check that the Toeplitz representation of $X$ into
$\mathcal{T}$ is CNP-covariant, and it thus follows that $\NO{X}$, and
hence also the $\NO{X}^\reduced$ which we will construct later, are
both isomorphic to $\mathcal{T}$ and not to $C(\TT)$.
\end{rmk}

To prove Theorem~\ref{thm:projective property}, we first need
to recall a few facts about Fell bundles and their
$C^*$-algebras, and about coactions. The main reference to Fell
bundles and properties of the full cross-sectional algebra of a
bundle is \cite[Section VIII.17.2]{Fell}. For the relationship
between topologically graded $C^*$-algebras and $C^*$-algebras
associated to Fell bundles, in particular the {reduced}
$C^*$-algebra of a bundle, we refer to \cite{Exel}. The
connection between discrete coactions and Fell bundles was
explored in \cite{qui:discrete coactions}. In \cite[Definition
3.5]{qui:discrete coactions}, Quigg introduced a reduced
$C^*$-algebra of a Fell bundle together with a coaction. The
subtle point that the reduced constructions from \cite{Exel}
and \cite{qui:discrete coactions} are compatible (although far
from obviously so) was clarified in \cite[page 749]{EQ}.

\begin{notation}\label{ntn:fell bundle}
Suppose that $\delta$ is a coaction of a discrete group $G$ on
a  $C^*$-algebra $A$. For every $g\in G$, let
$A^\delta_g:=\{\, a\in A \mid \delta(a)=a\otimes i_G(g) \,\}$ be the
spectral subspace of $A$ at $g$. By \cite{qui:discrete
coactions}, the disjoint union of the spectral subspaces
$A^\delta_g\times \{g\}$ for $g\in G$ forms a Fell bundle over $G$,
which we call \emph{the Fell bundle associated to} $\delta$
(see \cite[page 748]{EQ}).

Conversely, if $(\mathcal{A},G)$ is a Fell bundle then it follows
from \cite[Proposition 3.3]{qui:discrete coactions} that there is
a canonical coaction $\delta_{\mathcal{A}}$ on the full
cross-sectional algebra $C^*(\mathcal{A})$ such that
$\delta_{\mathcal{A}}(a_g)=a_g\otimes i_G(g)$ for all $a_g$ in
the fiber of $\mathcal{A}$ over $g$ and
 all $g$ in $G$.

If $(\mathcal{A},G)$ is a Fell bundle over $G$ and $A$ a
cross-sectional algebra of $(\mathcal{A},G)$ (in the sense that
$A$ is a $C^*$-completion of the algebra of finitely supported
sections on $\mathcal{A}$), then we say that $A$ is
topologically graded if there exists a contractive conditional
expectation from $A$ to $A_e$ which vanishes on $A_g$ for each
$g \in G \setminus \{e\}$ (see \cite[Definition 3.4]{Exel}).
The \emph{reduced cross-sectional algebra} $C^*_r(\mathcal{A})$ defined in
\cite{Exel} was shown to be minimal among topologically graded
cross-sectional algebras $A$. To be more precise, if $A$ is any
topologically graded cross-sectional algebra of $(\mathcal{A},G)$,
\cite[Theorem 3.3]{Exel} shows that there exists a surjective
homomorphism $\lambda_\mathcal{A}:A\to C^*_r(\mathcal{A})$ such that
$\lambda_{\mathcal{A}}\circ\eta_\mathcal{A}=\kappa_\mathcal{A}$
where $\eta_\mathcal{A}$ and $\kappa_\mathcal{A}$ are the
embeddings of the algebra of finitely supported
sections on $\mathcal{A}$ into $A$ and $C^*_r(\mathcal{A})$, respectively.
On the other hand, the universal property of $C^*(\mathcal{A})$
(see \cite[VIII.16.11]{Fell}) gives a surjective homomorphism
\begin{equation}\label{def:map_from_full_algebra}
    \phi_{\mathcal A} \colon C^*(\mathcal{A})\to A
\end{equation}
such that $\phi_\mathcal{A}\circ \gamma_\mathcal{A}=\eta_\mathcal{A}$
where $\gamma_\mathcal{A}$ is the embedding of the algebra of finitely supported
sections on $\mathcal{A}$ into $C^*(\mathcal{A})$.

If $\delta$ is a coaction of $G$ on $A$ and $(\mathcal{A}, G)$
is the Fell bundle associated to $\delta$, it follows from
\cite[Lemma 1.3]{qui:discrete coactions} (see also \cite[page 749]{EQ})
that $A$ is a topologically graded cross-sectional algebra of
$\mathcal{A}$. We shall adopt the {notation $A^\reduced$} for
the reduced cross-sectional algebra of the bundle
$(\mathcal{A}, G)$ arising from the coaction $\delta$ on $A$.
(We choose not to cram $\delta$ into the notation $A^\reduced$
for the sake of readability: the coaction $\delta$ will always
be clear from context.) By the considerations of the previous
paragraph applied to $A$ and $C^*({\mathcal A})$, there are
surjective homomorphisms
\begin{equation}\label{def:lambda_homom}
\lambda_{\mathcal A} \colon A\to A^\reduced \text{ and }\Lambda_{\mathcal A} \colon C^*({\mathcal A})\to A^\reduced
\end{equation}
such that
$\lambda_{\mathcal{A}}\circ\eta_\mathcal{A}=\kappa_\mathcal{A}$
and
$\Lambda_{\mathcal{A}}\circ\gamma_\mathcal{A}=\kappa_\mathcal{A}$;
and hence $\lambda_{\mathcal{A}}\circ
\phi_{\mathcal{A}}=\Lambda_{\mathcal{A}}$ (see, for example,
\cite{EQ}).
% Let $\iota_A$ denote the map defined on $\bigcup_{g\in G}A^\delta_g$
% satisfying that its restriction to each $A^\delta_g$ is the
% composition of the identification of $A^\delta_g$ with
% $A^\delta_g\times\{g\}$ and the canonical embedding of
% $A^\delta_g\times\{g\}$ into the algebra of finitely supported
% sections on $\mathcal{A}$. Then
% $\phi_\mathcal{A}\circ\gamma_\mathcal{A}\circ\iota_A=\eta_\mathcal{A}\circ\iota_A$
% is the inclusion of $\bigcup_{g\in G}A^\delta_g$ into $A$, and
% $\kappa_\mathcal{A}\circ\iota_A$ is the restriction of
% $\lambda_\mathcal{A}$ to $\bigcup_{g\in G}A^\delta_g$.

As explained in  \cite[page 749]{EQ}, $A^\reduced$ (defined by
its minimality, or co-universal property) is the same as the
reduced algebra from \cite{qui:discrete coactions} associated
to $(\mathcal{A}, G)$. By \cite[Definition 3.5]{qui:discrete
coactions}, there exists a coaction $\delta^\normal$ (the
$\normal$ stands for ``normal"; see
Remark~\ref{rmk:reduced/normal}) on $A^\reduced$ with the
property that
\begin{equation}\label{def_delta_reduced}
    \delta^\normal\bigl(\lambda_\mathcal{A}(a_g)\bigr)=
    \lambda_\mathcal{A}(a_g)\otimes i_G(g)\text{ for all }a_g\in A^\delta_g.
\end{equation}
Recall that a coaction $\eta$ of $G$ on a $C^*$-algebra $C$ is
called \emph{normal} if $(\id \otimes \lambda_G)\otimes \eta$
is injective. Every coaction $\eta$ of $G$ on $C$ has a
\emph{normalisation}: the quotient $C^{\normal}$ of $C$ by
$\ker ((\id \otimes \lambda_G)\otimes \eta)$ carries a coaction
$\tilde{\eta}$ which is automatically normal. In our set-up,
$A^{\reduced}$ is isomorphic to $A^{\normal}$, and this
isomorphism identifies the coaction $\delta^{\normal}$ defined
by~\eqref{def_delta_reduced} with the normalisation
$\tilde{\delta}$ of $\delta$ (see \cite[Lemma 2.1]{EQ}).
Moreover, $\delta^{\normal}$ may also be identified with the
normalisation of $\delta_{\mathcal{A}}$ by construction (see
\cite{qui:discrete coactions}). In particular, as the notation
suggests, $\delta^\normal$ is a normal coaction on
$A^\reduced$.

\begin{rmk}\label{rmk:reduced/normal}
Our choice of notation $(A^{\reduced}, \delta^{\normal})$ for
the system obtained above from $(A,\delta)$ may seem a little
perverse when either $(A^{\reduced}, \delta^{\reduced})$ or
$(A^{\normal}, \delta^{\normal})$ would at least be internally
consistent. We have our reasons. The notation $A^{\reduced}$
is, for us, much more appealing than $A^{\normal}$ for two
reasons: firstly, it coincides with our key
reference~\cite{Exel}; and secondly, there is strong evidence
that the object obtained in this way from the algebra $\NO{X}$
of \cite{SY} should be regarded as a reduced crossed product
(see Section~\ref{sec:examples}). However, the notation
$\delta^{\reduced}$ would be a most unfortunate choice because
it suggests a reduced coaction (that is, one taking values in
$A \otimes C^*_\reduced(G)$) whereas we have been careful to use only
full coactions throughout this paper for the sake of
consistency and self-containment; in particular,
$\delta^{\normal}$ is a full normal coaction. See the notation
after \cite[Definition 3.5]{qui:discrete coactions} for a
similar point of view.
\end{rmk}

%It follows from \cite{qui:discrete coactions} that $\{A^\delta_g\}_{g\in
%  G}$ is a topological $G$-grading of $A$ (cf. \cite{Exel}). On the
%other hand, if $\{A_g\}_{g\in G}$ is a topological $G$-grading, then
%it follows from \cite[Proposition 3.3]{qui:discrete coactions} that
%there exists a unique coaction $\delta$ of $G$ on $A$ such that
%$A^\delta_g=A_g$ for every $g\in G$.

%If $G$ is a discrete group $\{A_g\}_{g\in G}$ is a topological
%$G$-grading on a $C^*$-algebra $A$, then $\{A_g\times\{g\}\}_{g\in G}$
%is a \emph{Fell bundle} (also called a $C^*$-algebraic bundle, cf. \cite{Fell})
%over $G$ (cf. \cite[page 748]{EQ}).
%On the other hand, if $\mathcal{B}$ is a Fell bundle and
%$C^*(\mathcal{B})$ is its full cross-sectional algebra (cf. \cite[17.2]{Fell}, \cite[page 45]{Exel}
%and \cite[page 749]{EQ}), then there exists a $G$-topological grading $(C^*(\mathcal{B})_g)_{g\in G}$
%of $C^*(\mathcal{B})$ such that the two Fell bundles $\mathcal{B}$ and $(C^*(\mathcal{B})_g\times\{g\})_{g\in G}$
%are isometrically isomorphic, cf. \cite[16.11]{Fell}.

%It follows that if $\delta$ is a coaction of a discrete group
%$G$ on a $C^*$-algebra $A$, then $\{A_g\times\{g\}\}_{g\in G}$
%is a Fell bundle over $G$, which we call \emph{the Fell bundle
%associated to} $\delta$.
\end{notation}

\begin{rmk}\label{rmk:NO(X)=x-sect. alg}
Let $(G,P)$ be a quasi-lattice ordered group, and let $X$ be a
compactly aligned product system of Hilbert bimodules over $P$.
By Proposition~\ref{prop:existence_coaction}, $\Tc(X)$ admits a
coaction $\delta$, and hence gives rise to a Fell bundle
$\mathcal{B}=(\Tc(X)^\delta_g\times\{g\})_{g\in G}$ over $G$.
The generalised fixed-point algebra $\Tc(X)^\delta_e$ is
precisely the algebra $\Ff$ of~\eqref{eq:F-def}.

Let $\iota$ denote the map from $\bigcup_{g\in
G}\Tc(X)^\delta_g$ to the algebra of finitely supported
sections on $\mathcal{B}$ such that the restriction of $\iota$
to each $\Tc(X)^\delta_g$, identified with
$\Tc(X)^\delta_g\times\{g\}$, is the canonical embedding of
$\Tc(X)^\delta_g\times\{g\}$. We then have that
$\eta_\mathcal{B}\circ\iota$ is the inclusion of $\bigcup_{g\in
G}\Tc(X)^\delta_g$ into $\Tc(X)$. We claim that the map
$\gamma_\mathcal{B}\circ\iota\circ i_X:X\to C^*(\mathcal{B})$
is a Nica covariant Toeplitz representation. Indeed, to check
this we use that $i_X$ is a Nica covariant representation and
that $\gamma_\mathcal{B}\circ\iota$ is compatible with the
multiplication and involution and restricts to a linear map on
each fiber $\Tc(X)^\delta_g$ and to a $*$-homomorphism on the
fiber $\Tc(X)^\delta_e$. Then the universal property of
$\Tc(X)$ supplies a $*$-homomorphism $\zeta:\Tc(X)\to
C^*(\mathcal{B})$ such that $\zeta\circ i_X=
\gamma_\mathcal{B}\circ\iota\circ i_X$. By
\eqref{def:map_from_full_algebra}, there is a surjective
homomorphism $\phi_\mathcal{B} \colon C^*(\mathcal{B})\to
\Tc(X)$ such that $\phi_\mathcal{B}
\circ\gamma_\mathcal{B}=\eta_\mathcal{B}$. We then have
\begin{equation*}
  \phi_\mathcal{B}(\zeta(i_X(x))) = \phi_\mathcal{B}(\gamma_\mathcal{B}(\iota(i_X(x))))
  = \eta_\mathcal{B}(\iota(i_X(x))) = i_X(x)
\end{equation*}
and
\begin{equation*}
  \zeta(\phi_\mathcal{B}(\gamma_\mathcal{B}(\iota(i_X(x))))
  = \zeta(\eta_\mathcal{B}(\iota(i_X(x))))
  = \zeta(i_X(x)) = \gamma_\mathcal{B}(\iota(i_X(x)))
\end{equation*}
for each $x\in X$, from which it follows that $\zeta$ is the inverse of $\phi_\mathcal{B}$.
Hence $\phi_\mathcal{B}$ is an isomorphism from $C^*(\mathcal{B})$ to $\Tc(X)$ which is equivariant for $\delta_\mathcal{B}$ and $\delta$.
% The
% universal property of $\Tc(X)$ supplies an inverse to
% $\phi_\mathcal{B}$. Thus $\psi := \phi_{\mathcal{B}}^{-1} \circ
% i_X$ is a Nica covariant representation of $X$ in
% $C^*(\mathcal{B})$.

Suppose that $X$ is $\tilde{\phi}$-injective. Let $\NO{X}$ be
the Cuntz-Nica-Pimsner algebra of $X$ and $j_X$ the universal
CNP-covariant representation. By the proof of \cite[Proposition 3.12]{SY}
and \eqref{eq:BPs span F}, the kernel of the canonical
homomorphism $q_{\CNP} \colon \Tc(X)\to \NO{X}$ is generated by its
intersection with $\Ff$. Therefore Proposition~\ref{prp:deltaI
existence} applied to the coaction $\delta$ on $\Tc(X)$ yields
a \emph{gauge coaction} $\CNPgaug$ on $\NO{X}$. The spectral
subspaces
\[(\NO{X})^{\CNPgaug}_g:= \{\, c\in\NO{X} \mid \CNPgaug(c)=c\otimes
i_G(g) \,\}
\]
give rise to a Fell bundle $\mathcal{N}$, and it follows as above
from the universal property of $\NO{X}$ (see \cite[Proposition 3.12]{SY})
that $\phi_{\mathcal{N}} \colon C^*(\mathcal{N})\to \NO{X}$
is an isomorphism which is equivariant for $\delta_{\mathcal N}$
and $\CNPgaug$.
% . The universal
% property of $\NO{X}$ (see \cite[Proposition 3.12]{SY}) implies
% that there is an inverse for $\phi_{\mathcal{N}}$, so we may
% identify the coaction $\CNPgaug$ on $\NO{X}$ with the coaction
% $\delta_{\mathcal N}$ on $C^*(\mathcal{N})$.
\end{rmk}

\begin{proof}[Proof of Theorem~\ref{thm:projective property}]
Let $\NO{X}^\reduced$ be the reduced cross-sectional algebra of
the Fell bundle $\mathcal{N}$, and $\lambda_{\mathcal{N}}
\colon \NO{X} \to\NO{X}^\reduced$ the homomorphism from
\eqref{def:lambda_homom}. Put $j^\reduced_X :=\lambda_{\mathcal
N} \circ j_X$, and  define $\redgaug$ to be the normal coaction
on $\NO{X}^\reduced$ described in \eqref{def_delta_reduced}.
% ; as
% in Notation~\ref{ntn:fell bundle}, we identify $\redgaug$ with
% the normalisation of $\CNPgaug$.

To prove property~(1), note that $\NO{X}$ is generated by the
injective CNP-covariant representation $j_X$, see
\cite[Proposition 3.12 and Theorem~4.1]{SY}. So $j_X^\reduced$
is CNP-covariant. Since $\lambda_{\mathcal N}$ is surjective,
$\NO{X}^\reduced$ is generated by $j_X^\reduced$. Further,
$\lambda_{\mathcal N}$ restricts to a bijection from
$(\NO{X})^{\CNPgaug}_e$ to $(\NO{X}^\reduced)^{\redgaug}_e$,
and since $j_X(A)\subset (\NO{X})^{\CNPgaug}_e$ the
representation $j_X^\reduced$ is injective. Finally it follows
from~\eqref{def_delta_reduced} that
$\redgaug(j^\reduced_X(x))=j^\reduced_X(x)\otimes i_G(d(x))$
for all $x\in X$.
% Since $\redgaug$ is isomorphic to the normalisation
% of $\CNPgaug$, we have $(\lambda_{\mathcal
% N}\otimes\id_{C^*(G)})\circ\CNPgaug=\redgaug\circ\lambda_{\mathcal
% N}$. Thus
% \[
% \redgaug(j^\reduced_X(x)) = j^\reduced_X(x)\otimes i_G(d(x))\quad\text{for all $x \in X$,}
% \]
% so $j^\reduced_X$ is gauge-compatible.

We next show that $(\NO{X}^\reduced, j^\reduced_X)$ has
property~(2). Suppose $\psi \colon X \to B$ is as in~(2), and
let $\beta$ be a coaction on $B$ such that
\eqref{E:gauge-compatible} holds. For $g\in G$, let
$B^\beta_g=\{\, b\in B \mid \beta(b)=b\otimes i_G(g) \,\}$.
Then \cite[Lemma 1.3 and 1.5]{qui:discrete coactions} and
\cite[Theorem 3.3]{Exel} imply that $\{B^\beta_g\}_{g\in G}$ is
a topological grading of $B$. The universal property of
$\Tc(X)$ gives a $*$-homomorphism $\psi_* \colon \Tc(X) \to B$
such that $\psi = \psi_* \circ i_X$. Since the image of $\psi$
generates $B$, $\psi_*(\Ff)=B^\beta_e$, and so
$I_0:=\psi_*(\ker(q_{\CNP}))\cap
B^\beta_e=\psi_*(\ker(q_{\CNP})\cap \Ff)$ is an ideal of
$B_e^\beta$. Let $I$ be the ideal of $B$ generated by $I_0$.

By construction, $I$ is an induced ideal in the sense of
\cite[Definition 3.10]{Exel}. Let $\pi \colon B\to B/I$ be the
quotient map. By \cite[Proposition 3.11]{Exel},
$\{\pi(B^\beta_g)\}_{g\in G}$ is a topological grading of
$B/I$.

Since the image of $\psi$ generates $B$, we have
$\pi(\psi_*(\Tc(X)^\delta_g))=\pi(B^\beta_g)$ for all $g\in G$.
We aim to show that for every $g\in G$ we have
\begin{equation}\label{same_kernel}
\ker
(\pi\circ\psi_*)\cap \Tc(X)^\delta_g=\ker (q_{\CNP})\cap \Tc(X)^\delta_g.
\end{equation}

It will then follow that the two Fell bundles $\tilde{\mathcal{B}}:=(\pi(B^\beta_g)\times\{g\})_{g\in G}$ and
 $((\NO{X})^{\CNPgaug}_g \times\{g\})_{g\in G}$ are isometrically
isomorphic. Indeed, \eqref{same_kernel} implies that for every
$g\in G$ there is an isomorphism $\varphi_g$ from
$(\NO{X})^{\CNPgaug}_g$ onto $\pi(B^\beta_g)$  given by
$\varphi_{d(x)}(q_{CNP}(i_X(x)))=\pi\circ\psi_*(i_X(x))$ for
all $x\in X$. These isomorphisms are compatible with the Fell
bundle structure, in the sense that:
$\varphi_{g_1}(c_1)\varphi_{g_2}(c_2)=\varphi_{g_1g_2}(c_1c_2)$
for $g_1,g_2\in G$, $c_1\in (\NO{X})^\nu_{g_1}$ and $c_2\in
(\NO{X})^\nu_{g_2}$; and $\varphi_g(c)^*=\varphi_{g^{-1}}(c^*)$
for $g\in G$ and $c\in (\NO{X})^\nu_g$. Hence the isomorphisms
$\varphi_g$ induce the claimed isomorphism between
$\tilde{\mathcal{B}}$ and $\mathcal{N}$. Thus every cross
sectional algebra of $\tilde{\mathcal{B}}$ is also a cross
sectional algebra of $\mathcal{N}$ and vice versa, and
% By \cite{Exel} there are
% surjective homomorphisms $\Lambda_{\mathcal{N}} \colon
% \mathcal{N}\to \NO{X}^\reduced$ and $\Lambda_{\mathcal{B}}
% \colon \mathcal{B}\to C^*_r(\mathcal{B})$ which restrict to the
% identity homomorphism on the fibre over $g$ for all $g\in G$
% (note that this notation is not in conflict with the similar
% notation in \eqref{def:lambda_homom}; each of
% $\Lambda_{\mathcal{N}}$, $\Lambda_{\mathcal{B}}$ is a regular
% representation, of the Fell bundle and of its full
% cross-sectional $C^*$-algebra).
the co-universal properties
of $\NO{X}^\reduced$ and $C_r^*(\tilde{\mathcal{B}})$ then imply that
there is an isomorphism
$\tilde{\varphi} \colon \NO{X}^\reduced\to C^*_r(\tilde{\mathcal{B}})$
such that $\lambda_{\tilde{\mathcal{B}}}\circ \varphi_g$ and
$\tilde\varphi \circ \lambda_{\mathcal{N}}$ agree on
$(\NO{X})^{\CNPgaug}_g$ for all $g$.

Let $\phi:=(\tilde\varphi)^{-1}\circ
\lambda_{\tilde{\mathcal{B}}}\circ \pi$. Then $\phi$ is a
homomorphism from $B$ onto $\NO{X}^\reduced$, and we have
\begin{align*}
\phi(\psi(x))
&=(\tilde\varphi)^{-1}(\lambda_{\tilde{\mathcal{B}}}(\pi(\psi_*(i_X(x)))))=(\tilde\varphi)^{-1}
(\lambda_{\tilde{\mathcal{B}}}(\varphi_{d(x)}(q_{\CNP}(i_X(x)))))\\
&=\lambda_{\mathcal{N}}(q_{\CNP}(i_X(x)))=\lambda_{\mathcal{N}}(j_X(x))=j_X^r(x),
 \end{align*}
for all $x \in X$, as claimed.

%It will then follow that the two gradings
%$(\pi(B^\beta_g))_{g\in G}$ and $((\NO{X})^{\CNPgaug}_g)_{g\in
%G}$ are isometrically isomorphic, and thus $\lambda_{\mathcal
%N}$ induces a $*$-homomorphism $\tilde\phi \colon B/I\to
%\NO{X}^\reduced$.
%%$\tilde\phi(\pi(\psi_*(x)))=\lambda_{\NO{X}}(q_{\CNP}(x))$ for each $x\in
%\Tc(X)^\delta_g$ and each $g\in G$.
%The composition
%$\phi:=\tilde\phi\circ\pi$ is then a $*$-homomorphism $\phi \colon B
%\to \NO{X}^\reduced$ such that $\phi(\psi(x)) =
%j^\reduced_X(x)$ for all $x \in X$.

We first prove \eqref{same_kernel} when $g=e$. So we claim that $\ker
(\pi\circ\psi_*)\cap \Ff=\ker (q_{\CNP})\cap \Ff$.
If $c\in \ker (q_{\CNP})\cap \Ff$, then
$\psi_*(c)\in I_0\subset I$, proving the right to left
inclusion. To prove the other inclusion, note that since $\pi\circ \psi_*=(\pi\circ \psi)_*$, it
suffices by Proposition \ref{lemma:inj} to show that the Toeplitz representation $\pi\circ\psi$ is injective.

Fix  $a\in A$ with $\pi(\psi(a))=0$.  Then $\psi(a)\in I\cap
B_e^\beta$. Since $\psi_*(\ker(q_{\CNP}))$ is an
ideal, it contains $I$,  and therefore $ I\cap B_e^\beta\subset
\psi_*(\ker(q_{\CNP}))\cap B_e^\beta=I_0$. It follows that
there exists a $y\in\ker(q_{\CNP})\cap \Ff$ such that
$\psi_*(y)=\psi(a)$.
%We have that
%\begin{equation*}
%(\psi_*\otimes1_{C^*_\reduced(G)})(\delta(y))=\beta(\psi_*(y))
%=\beta(\psi(a))=\psi(a)\otimes i_G(e)=\psi_*(y)\otimes i_G(e),
%\end{equation*}
%from which it follows that
%$y\in \Ff$.
Hence $y-i_X(a)\in\ker(\psi_*)\cap\Ff$.
Since $\ker(\psi_*)\cap\Ff\subset \ker (q_{\CNP})$ by Proposition \ref{lemma:inj} applied to $\psi$,
it follows that $i_X(a)\in\ker (q_{\CNP})\cap i_X(A)$.
However, $\ker (q_{\CNP})\cap i_X(A)=\{0\}$, hence
$i_X(a)=0$, and therefore  $a=0$.
This proves that the representation $\pi\circ\psi$ is injective, and thus that $\ker
(\pi\circ\psi_*)\cap \Tc(X)^\delta_e=\ker (q_{\CNP})\cap \Tc(X)^\delta_e$.

Now let $g$ be an arbitrary element of $G$.
Then
\begin{align*}
c\in\ker (\pi\circ\psi_*)\cap \Tc(X)^\delta_g&\Longleftrightarrow c^*c\in \ker
(\pi\circ\psi_*)\cap \Tc(X)^\delta_e=\ker (q_{\CNP})\cap \Tc(X)^\delta_e,\\
&\Longleftrightarrow
c\in\ker (q_{\CNP})\cap \Tc(X)^\delta_g.
\end{align*}
Hence (2) is established.

Finally, for the uniqueness assertion, suppose that $(C, \rho,\gamma)$
also satisfies (1) and (2). Then property~(2) for $\NO{X}^\reduced$
gives a homomorphism $\phi$ from $C$ to $\NO{X}^\reduced$, the
corresponding property for $C$ gives a homomorphism from
$\NO{X}^\reduced$ to $C$, and these two homomorphisms are mutually
inverse.
\end{proof}

We saw in Remark~\ref{rmk:NO(X)=x-sect. alg} that $\NO{X}$ is
isomorphic to the full cross sectional algebra
$C^*(\mathcal{N})$ of its associated Fell bundle arising from
the gauge coaction $\nu$. This fact has interesting
implications. To explain this point, we need to recall some
terminology from \cite{EKQ}. First, a coaction $\eta$ of a
group $G$ on $C$ is maximal if the canonical map from the
iterated coaction crossed product $C\times_\eta
G\times_{\hat\eta} G$ to $C\otimes \mathcal{K}(l^2)$ is an
isomorphism. Second, a maximal coaction system $(D,G,
\epsilon)$ is a maximalisation of $(C, G, \eta)$ if there is an
equivariant surjective homomorphism from $D$ to $C$ which
induces an  isomorphism of the coaction crossed products
$D\times_\epsilon G$ and $C\times_\eta G$. Then
\cite[Proposition 4.2]{EKQ} implies that the canonical coaction
$\nu_{\mathcal{N}}$ on $C^*(\mathcal{N})$ is a maximalisation
of $\nu$, in fact the unique one with the same underlying Fell
bundle; moreover $\NO{X}\cong C^*(\mathcal{N})$ means precisely
that $\nu$ itself is maximal. At the other extreme, a system
$(D,G, \epsilon)$ is called a normalisation of $(C, G, \eta)$
if the coaction $\epsilon$ is normal and there is an
equivariant surjective homomorphism from $C$ to $D$ which
induces an  isomorphism of the coaction crossed products
$C\times_\eta G$ and $D\times_\epsilon G$.

\begin{cor}\label{cor:maximalisation}
Assume the hypotheses of Theorem~\ref{thm:projective property}.
Let $\psi:X\to B$ be an injective CNP-covariant representation
of $X$ which is gauge-compatible for a coaction $\beta$ on $B$,
and such that $\psi$ generates $B$. Then the following hold.
\begin{enumerate}\renewcommand{\theenumi}{\alph{enumi}}
\item\label{it:maximalisation} The coaction system
    $(\NO{X},G, \nu)$ is a maximalisation of $(B,G,\beta)$.
\item\label{it:normalisation} The coaction system
    $(\NO{X}^\reduced,G,\redgaug)$ is a normalisation of
    $(B,G,\beta)$.
\end{enumerate}
\end{cor}
\begin{proof}
The universal property of $\NO{X}$ gives a surjective
homomorphism $\intfrm{\psi}:\NO{X}\to B$ which is equivariant
for $\nu$ and $\beta$. Theorem~\ref{thm:projective property}
gives a homomorphism $\phi:B\to \NO{X}^\reduced$ which is
equivariant for $\beta$ and $\redgaug$. We have $\phi\circ
\intfrm{\psi}=\lambda_{\mathcal{N}}$. Then it follows as in the
proof of \cite[Lemma 2.1]{EQ} that the induced map
$\lambda_{\mathcal{N}}\times G$ is an isomorphism from
$\NO{X}\times_\nu G$ onto $\NO{X}^\reduced \times_{\redgaug} G$
and satisfies $\lambda_{\mathcal{N}}\times G=(\phi \times
G)\circ (\intfrm{\psi} \times G)$. Therefore $\intfrm{\psi}
\times G$ is also an isomorphism, which shows that $\nu$ is a
maximalisation of $\beta$, as claimed
in~(\ref{it:maximalisation}). Likewise, $\phi \times G$ becomes
an isomorphism from $B\times_\beta G$ onto
$\NO{X}^\reduced\times_{\redgaug} G$,
so~(\ref{it:normalisation}) follows.
\end{proof}

The main reason for proceeding to a gauge-invariant uniqueness
theorem for $\NO{X}$ via Theorem~\ref{thm:projective property}
rather than proving it directly from Theorem~\ref{thm:inj on
core} is that we feel that the co-universal property as a
defining property of $\NO{X}^\reduced$ is just as important as
--- and in some ways more natural than --- the defining universal property
of $\NO{X}$.

In particular when $\NO{X}$ and $\NO{X}^\reduced$ coincide, the
definition as a co-universal $C^*$-algebra has the advantage
over the definition as a universal $C^*$-algebra that it
involves only the natural defining relations for $\Tc(X)$ which
are present in the Fock representation, and does not involve
the complicated (and difficult to check) Cuntz-Pimsner
covariance condition. This has clear advantages when trying to
establish Cuntz-Nica-Pimsner algebras as models for other
classes of examples (see Remark~\ref{rmk:BQ analysis}).
Moreover, when $\NO{X}$ and $\NO{X}^\reduced$ do not coincide,
it is unclear what makes $\NO{X}$ worthy of singling out for
study beyond the bare fact that it is defined by a universal
property with respect to a relation which holds in the
co-universal $C^*$-algebra.

\begin{example}
If the advantage of the definition of $\NO{X}^\reduced$ as a
co-universal algebra is that it bypasses the troublesome
Cuntz-Pimsner covariance condition, then a natural next
question is whether or not $\NO{X}^\reduced$ is actually
co-universal for injective (not necessarily Nica covariant) gauge-compatible
Toeplitz representations of $X$, thus allowing
us to bypass the Nica covariance condition as well.

The answer in general is ``No:" there exist product systems
which admit no such co-universal $C^*$-algebra. To see this,
let $(G,P)=(\FF_2,\FF_2^+)$ and let $X_p=\CC$ for every $p\in
P$. Then $\mathcal{L}(X_p)=\mathcal{K}(X_p)$ for all $p$, so $X$ is compactly aligned,
and all left actions are injective. Denote the generators
of $\FF_2$ by $a$ and $b$. Suppose that $(C,\rho)$ is a
co-universal pair for injective gauge-compatible Toeplitz representations of $X$.
Note that there is an
injective Toeplitz representation $\psi$ of $X$ on $C^*_\reduced(\FF_2)$
determined by $\psi(1_p) := \lambda_{\FF_2}(p)$. By \cite[Example 1.15]{qui:full reduced} or
  \cite[Proposition 2.4]{qui:full duality}, there is a (full) coaction
  $\delta_{\FF_2}$ of $\FF_2$ on $C^*_\reduced(\FF_2)$ such that
  $\delta_{\FF_2}(\lambda_{\FF_2}(g)) = \lambda_{\FF_2}(g) \otimes i_{\FF_2}(g)$
  for all $g \in \FF_2$. The integrated form of $\psi$ is therefore
gauge-compatible. By the co-universal property of
$(C,\rho)$, there is a surjective homomorphism $\phi \colon
C^*_\reduced (\FF_2)\to C$ satisfying $\phi(\psi(1_p)) =
\rho(1_p)$ for all $p \in P$. In particular, since
$\lambda_{\FF_2}(a)$ and $\lambda_{\FF_2}(b)$ are unitaries,
surjectivity of $\phi$ implies that
\begin{equation}\label{eq:prod=1}
\rho(1_a)\rho(1_a)^* \rho(1_b)\rho(1_b)^*
    = \phi(\lambda_{\FF_2}(a)\lambda_{\FF_2}(a)^*\lambda_{\FF_2}(b)\lambda_{\FF_2}(b)^*)
    = \phi(1_{C^*_\reduced(\FF_2)})
    = 1_C.
\end{equation}
Since $j_X^\reduced$ is also an injective gauge-compatible
Toeplitz representation of $X$,  the co-universal property of
$(C,\rho)$ gives a surjective homomorphism $\eta \colon
\NO{X}^\reduced \to C$ such that $\eta(j_X^\reduced(1_p)) =
\rho(1_p)$ for all $p \in \FF_2$. Since $j_X^\reduced$ is Nica
covariant, and since $a \vee b = \infty$ in $\FF_2$, we have
\begin{equation}\label{eq:prod=0}
\rho(1_a)\rho(1_a)^* \rho(1_b)\rho(1_b)^*
    = \eta(j_X^\reduced(1_a)j_X^\reduced(1_a)^*j_X^\reduced(1_b)j_X^\reduced(1_b)^*)
    = \eta(0)
    = 0.
\end{equation}
Combining \eqref{eq:prod=1}~and~\eqref{eq:prod=0}, we obtain
$1_C = 0$, so $C = \{0\}$, which contradicts the assumption
that $\rho$ is an injective representation of $X$.
\end{example}

The preceding example shows that a co-universal $C^*$-algebra
for gauge-compatible injective Toeplitz representations need not exist.
We will now show that if such a co-universal $C^*$-algebra does
exist, then it must be isomorphic to $\NO{X}^\reduced$, and
prove that it satisfies a kind of rudimentary gauge-invariant
uniqueness theorem.

\begin{cor}
Let $(G,P)$ be a quasi-lattice ordered group and let $X$ be a
compactly aligned product system over $P$ of right-Hilbert
$A$--$A$ bimodules. Suppose either that the left action on each
fibre is injective, or that $P$ is directed and $X$ is
$\tilde\phi$-injective.
\begin{enumerate}
\item If $\phi \colon \NO{X}^\reduced\to B$ is a surjective
    $*$-homomorphism, then $\phi$ is injective if and only
    if $\phi{|j_X^\reduced(A)}$ is injective and there is a
    coaction
    $\beta$ of $G$ on $B$ such that
    $\beta\circ\phi=(\phi\otimes
    \id_{C^*(G)})\circ\redgaug$.
\item Suppose that $(C,\rho)$ is a co-universal pair for
    injective gauge-compatible (not necessarily Nica covariant)
    Toeplitz representations of $X$. Then
    there is a $*$-isomorphism $\phi \colon \NO{X}^\reduced\to C$
    such that $\phi(j_X^\reduced(x))=\rho(x)$ for all $x\in
    X$.
\end{enumerate}
\end{cor}
\begin{proof}
(1) The ``only if" assertion is trivial. For the ``if"
assertion note that $\phi\circ j_X^\reduced$ is an injective
Nica covariant representation of $X$ in $B$ which is
gauge-compatible, and whose image generates $B$. An application
of Theorem~\ref{thm:projective property} yields a surjection
$B\to \NO{X}^\reduced$ which is an inverse for $\phi$.

(2) Since $\NO{X}^\reduced$ is generated by an injective
Toeplitz representation of $X$ which is gauge-compatible, the
co-universal property of $(C,\rho)$ implies that there is a
surjective homomorphism $\phi \colon \NO{X}^\reduced\to C$ such that
$\phi(j_X^\reduced(x))=\rho(x)$ for all $x\in X$. Part~(1) then
implies that $\phi$ is injective and hence an isomorphism.
\end{proof}

Part~(1) of the preceding corollary is used to prove
statement~(2), but is somewhat unsatisfactory as a stand alone
result because there is no universal property to induce
homomorphisms $\phi \colon \NO{X}^\reduced \to B$ of the desired form
(compare with Definition~\ref{dfn:giup}
below). %
%\footnote{Question: is it true that $\NO{X}^\reduced$ is
%universal for CNP covariant representations of $X$ into
%$C^*$-algebras $B$ which carry normal coactions? Is this
%interesting? The point is that to use the gauge-invariant
%uniqueness theorem, one has to establish the existence of a
%coaction on $B$ anyway...}
The following result is much more
satisfactory in that it provides an injectivity criterion for
the homomorphism $\phi \colon B \to \NO{X}^\reduced$ induced by
Theorem~\ref{thm:projective property}(\ref{it:c-u property}).

\begin{cor}\label{cor:phi_inj}
Let $(G,P)$ be a quasi-lattice ordered group and let $X$ be a
compactly aligned product system over $P$ of right-Hilbert
$A$--$A$ bimodules. Suppose either that the left action on each
fibre is injective, or that $P$ is directed and $X$ is
$\tilde\phi$-injective. Let $\psi \colon X \to B$ be an
injective Nica covariant gauge-compatible
representation whose image generates $B$. Then the surjective
$*$-homomorphism $\phi \colon B \to \NO{X}^\reduced$ of
Theorem~\ref{thm:projective property}(\ref{it:c-u property}) is
injective if and only if $\psi$ is Cuntz-Pimsner covariant and
$\beta$ is normal.
\end{cor}
\begin{proof}
If $\phi$ is injective, then $\psi$ is Cuntz-Nica-Pimsner covariant
and $\beta$ is normal because $j_X^\reduced$ and $\redgaug$
have these properties.

Now suppose that $\psi$ is Cuntz-Nica-Pimsner covariant and $\beta$
is normal. The universal property of $\NO{X}$ gives a
homomorphism $\Pi\psi \colon \NO{X} \to B$. By definition,
$\lambda_{\mathcal{N}}
\colon \NO{X} \to\NO{X}^\reduced$ satisfies
$\lambda_{\mathcal{N}}
= \phi \circ \Pi\psi$. The map
$\lambda_{\mathcal{N}}$ restricts to an isomorphism of
$(\NO{X})^{\CNPgaug}_e$  to $(\NO{X}^\reduced)^{\redgaug}_e$,
hence $\phi$ restricts to an isomorphism $B^\beta_e \to
(\NO{X}^\reduced)^{\redgaug}_e$. Since $\beta$ is normal, it
determines a faithful conditional expectation $\Phi^\beta \colon B
\to B^\beta_e$. But $\phi$ intertwines $\Phi^\beta$ and the
conditional expectation from $\NO{X}^\reduced$ to
$(\NO{X}^\reduced)^{\redgaug}_e$, and so the standard argument
implies that $\phi$ is injective.
\end {proof}

For the next corollary, we first define some additional
terminology.

\begin{dfn}\label{dfn:giup}
Fix a quasi-lattice ordered group $(G,P)$ and a
$\tilde\phi$-injective compactly aligned product system $X$
over $P$. We say that $\NO{X}$ has the \emph{gauge-invariant
uniqueness property} provided that the following is satisfied.

A surjective $*$-homomorphism $\phi:\NO{X}\to B$ is injective if
and only if:
\begin{enumerate}\renewcommand{\theenumi}{GI\arabic{enumi}}
\item\label{it:carries coaction} there is a coaction
    $\beta$ of $G$ on $B$ such that $\beta\circ\phi=(\phi\otimes\id_{C^*(G)})\circ \nu$, and
\item\label{it:inj on A} the homomorphism $\phi\vert_{j_X(A)}$ is injective.
\end{enumerate}
%Suppose that for each CNP-covariant representation
%$\psi:X\to B$ whose image generates $B$, the induced homomorphism $\intfrm{\psi}$ is
%injective if and only if both of the following conditions hold:
%\begin{enumerate}\renewcommand{\theenumi}{GI\arabic{enumi}}
%\item\label{it:carries coaction} There is a coaction
%    $\beta$ of $G$ on $B$ satisfying \eqref{E:gauge-compatible};
%\item\label{it:inj on A} The homomorphism $\psi_e \colon A = X_e
%    \to B$ is injective.
%\end{enumerate}
%Then we say that $\NO{X}$ has the \emph{gauge-invariant
%uniqueness property.}
\end{dfn}

\begin{cor}[The gauge-invariant uniqueness theorem]\label{cor:guit}
Let $(G,P)$ be a quasi-lattice ordered group and let $X$ be a
compactly aligned product system over $P$ of right-Hilbert
$A$--$A$ bimodules.  Suppose either that the left action on
each fibre is injective, or that $P$ is directed and $X$ is
$\tilde\phi$-injective. The following are equivalent.
\begin{enumerate}
\item The coaction $\CNPgaug$ is normal. \label{item:6}
\item The coaction $\redgaug$ is maximal.
    \label{item:redgaug max}
\item The Fell bundle
  $\bigl((\NO{X})^{\CNPgaug}_g\times\{g\}\bigr)_{g\in G}$ is
  amenable. \label{item:7}
\item The $*$-homomorphism
  % Changed 22/2
  % $\Lambda_{\mathcal N}
  $\lambda_{\mathcal{N}}
  \colon \NO{X}\to\NO{X}^\reduced$
  is an isomorphism. \label{item:8}
\item $\NO{X}$ has the gauge-invariant uniqueness property. \label{item:9}
\item If $\psi_1:X\to B_1$ and $\psi_2:X\to B_2$ are two injective gauge-compatible
CNP-covariant representations of $X$ whose images generate $B_1$ and $B_2$
respectively, then there exists a $*$-isomorphism
$\phi:B_1\to B_2$ such that $\phi\circ\psi_1=\psi_2$. \label{item:10}
\end{enumerate}
\end{cor}
\begin{proof}
That (\ref{item:6})~and~(\ref{item:8}) are equivalent follows
from the fact that $\redgaug$ is the normalisation of
$\CNPgaug$ (see the last paragraph of Notation~\ref{ntn:fell
bundle}). The equivalence of (\ref{item:redgaug
max})~and~(\ref{item:8}) follows from
\cite[Proposition~4.2]{EKQ} and the fact, established in
Remark~\ref{rmk:NO(X)=x-sect. alg}, that $(\NO{X},G,\CNPgaug)$ is isomorphic
to $(C^*(\mathcal{N}),G,\CNPgaug_\mathcal{N})$.
The equivalence of
(\ref{item:7})~and~(\ref{item:8}) follows from the definition
of amenability for Fell bundles in \cite{Exel} and the fact
%established in Remark~\ref{rmk:NO(X)=x-sect. alg},
that $\NO{X}$ is isomorphic to $C^*(\mathcal{N})$. We will show that
(\ref{item:8})~implies~(\ref{item:9}), that
(\ref{item:9})~implies~(\ref{item:10}), and that
(\ref{item:10})~implies~(\ref{item:8}).

Suppose first that (\ref{item:8}) holds. Let $\phi:\NO{X}\to B$
be a surjective $*$-homomorphism.  We must show that  $\phi$ is
injective if and only if conditions (\ref{it:carries
coaction})~and~(\ref{it:inj on A}) hold. If $\phi$ is
injective, then we may define $\beta$ by $\beta := (\phi
\otimes \id_{C^*(G)}) \circ \CNPgaug \circ \phi^{-1}$, so
condition~(\ref{it:carries coaction}) is satisfied. Moreover,
\cite[Theorem~4.1]{SY} implies that $(\phi\circ j_X)|_A$ is
injective, so condition~(\ref{it:inj on A}) is also satisfied.
Now suppose that (\ref{it:carries coaction})~and~(\ref{it:inj
on A}) hold. Then $\phi\circ j_X$ is a  gauge-compatible Nica
covariant representation whose image generates $B$, so by
Theorem~\ref{thm:projective property} there is a
$*$-homomorphism $\tilde\phi:B\to \NO{X}^\reduced$ satisfying
$\tilde\phi\circ(\phi\circ j_X)=j_X^\reduced
=\lambda_{\mathcal{N}}\circ j_X$.  Since $\NO{X}$ is generated
by the elements $\{\, j_X(x) \mid x \in X \,\}$ and
$\lambda_{\mathcal{N}}$ is an isomorphism, it follows that
$(\lambda_{\mathcal{N}})^{-1}\circ \tilde\phi$ is an inverse
for $\phi$, and hence $\phi$ is injective as required.

Suppose that~(\ref{item:9}) holds. If $\psi:X\to B$ is an
injective gauge-compatible CNP-covariant representation of $X$
whose images generate $B$, then $\Pi\psi:\NO{X}\to B$ is a
surjective $*$-homomorphism such that $\Pi\psi\circ j_X=\psi$
and (\ref{it:carries coaction})~and~(\ref{it:inj on A}) are
satisfied, hence is an isomorphism. Statement~(\ref{item:10})
follows.

Finally suppose~(\ref{item:10}) holds. It follows from
\cite[Proposition 3.12 and Theorem~4.1]{SY} and
Theorem~\ref{thm:projective property} that $j_X:X\to\NO{X}$ and
$j_X^\reduced:X\to\NO{X}^\reduced$ are two injective
gauge-compatible CNP-covariant representations whose images
generate $\NO{X}$ and $\NO{X}^\reduced$ respectively. Thus
there exists a  $*$-isomorphism $\phi:\NO{X}^\reduced\to\NO{X}$
such that $\phi\circ j_X^\reduced=j_X$. We then have that
$\lambda_\mathcal{N}\circ\phi=\id_{\NO{X}^\reduced}$, from
which~(\ref{item:8}) holds.
\end{proof}

It is of course interesting to know under which conditions
$\NO{X}$ has the gauge-invariant uniqueness property. Using
Exel's work \cite{Exel}, we obtain the following conditions
under which $\NO{X}$ has the gauge-invariant uniqueness
property.

\begin{cor} \label{cor:amenable}
Let $(G,P)$ be a quasi-lattice ordered group and $X$  a
compactly aligned product system over $P$ of right-Hilbert
$A$--$A$ bimodules.  Suppose either that the left action on
each fibre is injective, or that $P$ is directed and $X$ is
$\tilde\phi$-injective. Then $\NO{X}$ has the gauge-invariant
uniqueness property in the following cases:
\begin{enumerate}
\item The group $G$ is exact and the coaction $\delta$ of $G$ on
  $\Tc(X)$ is normal. \label{item:1}
%\item \label{item:4} The group $G$ is exact and there is a quasi-lattice ordered
 % group $(\mathcal{G}, \mathcal{P})$ with $\mathcal{G}$ amenable and a
  %homomorphism $\pi \colon G \to \mathcal{G}$ such that whenever
 % $g, h \in G$ satisfy $g \vee h < \infty$, we have
  %\[
  %\pi(g) \vee \pi(h) =
  %\pi(g \vee h)\qquad\text{ and }\qquad \pi(g) = \pi(h) \implies g = h.
  %\]
\item The Fell bundle $\mathcal{B}=(\Tc(X)^\delta_g\times\{g\})_{g\in G}$ has the
approximation property. \label{item:2}
\item The Fell bundle
  $\mathcal{N}=\bigl((\NO{X})^{\CNPgaug}_g\times\{g\}\bigr)_{g\in G}$ has the
  approximation property. \label{item:5}
\item The group $G$ is amenable. \label{item:3}
\end{enumerate}
\end{cor}

\begin{proof}
Statement~(\ref{item:1}) follows from Corollary~\ref{cor:guit}
because normality of $\delta$ implies normality of $\nu$ by
Proposition \ref{prop:exact}.
%\eqref{item:4} follows from statement \eqref{item:1} and \cite[Theorem~8.1]{F99}, and in
%turn implies \eqref{item:3}.
For statement~(\ref{item:2}), let $J=\ker (q_{\CNP})$, let
$\Phi^\delta \colon \Tc(X)\to\Tc(X)^\delta_e$ be the
conditional expectation induced by the coaction $\delta$, and
let $\Phi^{\CNPgaug} \colon \NO{X}\to (\NO{X}^{\CNPgaug})_e$ be
the conditional expectation induced by the coaction $\CNPgaug$.
Then it follows from \cite[Proposition 3.6]{Exel} that
$\ker(\lambda_\mathcal{N})=\{\, b\in\NO{X} \mid
\Phi^{\CNPgaug}(b^*b)=0 \,\}$ and thus that
$\ker((j_X^\reduced)_*)=q_{\CNP}^{-1}(\ker(\lambda_\mathcal{N}))=
\{\, c\in\Tc(X) \mid \Phi^\delta(c^*c)\in J \,\}$. Hence, if
the Fell bundle $\mathcal{B}=(\Tc(X)^\delta_g\times\{g\})_{g\in
G}$ has the approximation property, then \cite[Proposition
4.10]{Exel} implies that
$\ker((j_X^\reduced)_*)=\ker(q_{\CNP})$. Thus
$\lambda_{\mathcal{N}}$ is an isomorphism, and
statement~(\ref{item:2}) then follows from
Corollary~\ref{cor:guit}. If the Fell bundle $\mathcal{N}$ has
the approximation property, then it is amenable by
\cite[Theorem 4.6]{Exel}, hence~(\ref{item:5}) follows from
Corollary~\ref{cor:guit}. Finally, if $G$ is amenable, then
\cite[Theorem 4.7]{Exel} shows that the bundle $\mathcal{N}$
has the approximation property, so~(\ref{item:3}) follows
from~(\ref{item:5}).
\end{proof}

\begin{rmk}
Observe that the universal property of $\NO{X}$ together with
the co-universal property of $\NO{X}^\reduced$ stated in
Theorem~\ref{thm:projective property} imply that the canonical
homomorphism from $\NO{X}$ to $\NO{X}^\reduced$ factors through
the image $B$ of any injective CNP-covariant representation
$\psi$ of $X$ which generates $B$ and respects $\delta$, hence
$\CNPgaug$. By Corollary~\ref{cor:maximalisation}, we see that
when the gauge-invariant uniqueness theorem  applies, it
implies that the universal and co-universal algebras for
gauge-compatible injective CNP-covariant representations of $X$
agree, the gauge coaction is both normal and maximal, and all
of the $C^*$-algebras $B$ coincide.
\end{rmk}

Our motivating example was $(G,P) = (\ZZ^k, \NN^k)$, in which
case condition~(\ref{it:carries coaction}) can be stated in the
familiar terms of an action of the dual group $\TT^k$.

\begin{cor}\label{cor:Nk giut}
Let $X$ be a compactly aligned product system over $\NN^k$.
%and let $\psi \colon X \to B$ be a surjective CNP-covariant
%representation of $X$. The induced map $\intfrm{\psi} \colon
%\NO{X} \to B$ is injective if and only if:
%\begin{enumerate}\renewcommand{\theenumi}{\roman{enumi}}
%\item\label{it:gauge action} there is a strongly continuous
 %   action $\alpha$ of $\TT^k$ on $B$ such that
  %  $\alpha_z(\psi(x)) = z^{d(x)}\psi(x)$ for all $x \in
   % X$ and $z\in\TT^k$;
%\item $\psi_e \colon A \to B$ is injective.
%\end{enumerate}
A surjective $*$-homomorphism $\phi:\NO{X}\to B$ is injective if and only if:
\begin{enumerate}\renewcommand{\theenumi}{\arabic{enumi}}
\item\label{it:gauge action} there is a strongly continuous
    action $\alpha$ of $\TT^k$ on $B$ such that
    $\alpha_z(\phi(j_X(x))) = z^{d(x)}\phi(j_X(x))$ for all $x \in
    X$ and $z\in\TT^k$, and
\item $\phi\vert_{j_X(A)} \colon A \to B$ is injective.
\end{enumerate}

\end{cor}
\begin{proof} Certainly $\NN^k$ is directed, and Lemma~4.3 of \cite{SY} implies that $X$ is
$\tilde\phi$-injective. Since $\ZZ^k$ is amenable, $\NO{X}$
has the gauge-invariant uniqueness property by
Corollary~\ref{cor:amenable}.

Every coaction $\beta$ of $\ZZ^k$ determines and is determined
by a strongly continuous action $\alpha$ of $\TT^k =
\widehat{\ZZ^k}$: specifically, $\beta(a) = a \otimes
i_{\ZZ^k}(m)$, if and only if $\alpha_z(a) = z^m a$ for all $z
\in \TT^k$. Hence condition~(\ref{it:carries coaction}) is
equivalent, in this setting, to condition~(\ref{it:gauge
action}).
\end{proof}

\section{Applications and examples}\label{sec:examples}

In this section we investigate a number of examples which
illustrate both the class of $C^*$-algebras $\NO{X}^\reduced$
and the utility of its co-universal property as set out in
Theorem~\ref{thm:projective property}.

\subsection*{Group crossed products}

Let $(G,P)$ be a quasi-lattice ordered group, and let $\alpha \colon
G \to \Aut(A)$ be an action of $G$ on a $C^*$-algebra $A$.
Suppose that $\omega$ is a $\TT$-valued cocycle on $G$; that is
$\omega \colon G \times G \to \TT$ satisfies $\omega(e,e) = 1$ and
\[
\omega(gh,k) \omega(g,h) = \omega(g,hk)\omega(h,k)\text{ for all $g,h,k \in G$.}
\]
Recall from \cite[Lemma~3.2]{F99} that there is a product
system $X:=X(\alpha,\omega)$ over the opposite semigroup $P^{\op}$ defined as follows:
for $p \in P$, let $X_p$ be the right-Hilbert
$A$--$A$ bimodule which is equal to $A$ as a normed vector
space with operations
\[
\langle x,y \rangle_A := x^*y \qquad a \cdot x = \alpha_p(a)x
\quad\text{ and }\quad x \cdot a = xa
\]
for all $x,y \in X_p$ and $a \in A$,
and define isomorphisms $X_p \otimes_A X_q \to X_{qp}$ by $x
\otimes_A y \mapsto \overline{\omega(q,p)}\alpha_q(x)y$.
Then $X(\alpha,\omega)=\coprod_{p\in P^{op}}X_p$ is the claimed product system (note that the left and right actions
are compatible with the product in $X$ because $\omega(p,e)=\omega(e,p)=1$ for all $p\in P$).
Moreover, the left action $\phi_p$ satisfies $\phi_p(A)\subset \mathcal{K}(X_p)$ for all $p\in P$,
and since each $X_p$ is essential in Fowler's sense,
\cite[Proposition~5.8]{F99} implies that
$X(\alpha, \omega)$ is compactly aligned.

The twisted crossed product $A \times_{\alpha,\omega} G$ is the
universal $C^*$-algebra generated by a covariant representation
of $(A,G,\alpha,\omega)$: that is, a homomorphism
$i^{\alpha,\omega}_A \colon A \to A \times_{\alpha,\omega} G$
and multiplier unitaries $\{\, i^{\alpha,\omega}_G(g) \mid g \in
G \,\}$ such that for $g,h \in G$ and $a \in A$,
\[
i_G^{\alpha,\omega}(g) i_G^{\alpha,\omega}(h) = \omega(g,h) i_G^{\alpha, \omega}(gh)
\quad\text{ and }\quad
i_G^{\alpha,\omega}(g) i_A^{\alpha,\omega}(a) i_G^{\alpha, \omega}(g)^* =
i_A^{\alpha,\omega}(\alpha_g(a));
\]
we have used $(i^{\alpha,\omega}_A, i^{\alpha,\omega}_G)$ in
this example, rather than the traditional $(i_A, i_G)$ to
distinguish the inclusion of $G$ as unitaries in $A
\times_{\alpha,\omega} G$ from its inclusion as unitaries in
$C^*(G)$. There is a coaction $\widehat{\alpha}$ of $G$ on $A
\times_{\alpha,\omega} G$ determined by
$\widehat{\alpha}(i^{\alpha,\omega}_A(a)i^{\alpha,\omega}_G(g))
= i^{\alpha,\omega}_A(a)i^{\alpha,\omega}_G(g) \otimes i_G(g)$
for all $a \in A$ and $g \in G$, and the universal properties
of the two algebras involved show that the crossed product $A
\times_{\alpha,\omega} G$ is isomorphic to the full
cross-sectional algebra of the resulting Fell bundle. The
reduced crossed product $A \times^r_{\alpha,\omega} G$ is the
reduced cross-sectional algebra $(A \times_{\alpha,\omega} G)^r$
of the same bundle, and is a quotient of $A \times_{\alpha,\omega} G$.
We write $(\lambda_A^{\alpha, \omega}, \lambda_G^{\alpha,\omega})$ for the generating covariant
representation of $(A, G, \alpha, \omega)$ in $A \times^r_{\alpha,\omega} G$.
Recall that $\widehat{\alpha}^n$ denotes the normalisation of $\widehat{\alpha}$
and is a normal coaction on $A \times^r_{\alpha,\omega} G$.

\begin{lemma}\label{lem:reduced CPs}
Let $(G,P)$ be a quasi-lattice ordered group such that $P$ is directed and $G$ is generated as a
group by $P$. Let $\alpha \colon G \to \Aut(A)$, $\omega \colon G \times
G \to \TT$ and $X$ be as above. Then there
is an isomorphism $\phi \colon A \times^r_{\alpha,\omega} G \to
\NO{X}^\reduced$ which takes
$(\lambda_G^{\alpha,\omega}(p))^* \lambda_A^{\alpha, \omega}(x)$ to
$j^\reduced_{X }(x)$ for all $x \in X_p = A$ and satisfies
$\redgaug\circ\phi=(\phi\otimes \id_{C^*(G)}) \circ
\widehat{\alpha}^n$.
\end{lemma}

\begin{proof}
We will first show that $A \times^r_{\alpha,\omega} G$ is
generated by a Nica covariant representation of $X$. We will
then apply Theorem~\ref{thm:projective property} to obtain a
surjective homomorphism $\phi$ from $A \times^r_{\alpha,\omega} G$
to $\NO{X}^\reduced$. Finally, we will use
the canonical faithful conditional expectation from $A
\times^r_{\alpha,\omega} G$ to $A$ to see that $\phi$ is
injective.

We begin by constructing a representation $\psi$ of $X$ in $A
\times^r_{\alpha,\omega} G$. For $p \in P$, define
$\psi_p \colon X_p \to A \times^r_{\alpha,\omega} G$ by
\[
\psi_p(x) := (\lambda_G^{\alpha,\omega}(p))^* \lambda_A^{\alpha, \omega}(x)\quad\text{ for all $x \in X_p = A$.}
\]
In the following calculations, we use $\diamond$ for the
multiplication in $P^{\op}$; so $p \diamond q = qp$ for all
$p,q \in P$. Fix $p,q \in P$ and elements $x \in X_p$ and $y
\in X_q$, and calculate:
\begin{align*}
\psi(x)\psi(y)
    &= (\lambda_G^{\alpha,\omega}(p))^* \lambda_A^{\alpha, \omega}(x) (\lambda_G^{\alpha,\omega}(q))^* \lambda_A^{\alpha, \omega}(y) \\
    &= (\lambda_G^{\alpha, \omega}(p))^* (\lambda_G^{\alpha,\omega}(q))^* (\lambda_G^{\alpha,\omega}(q))
    \lambda_A^{\alpha, \omega}(x) (\lambda_G^{\alpha,\omega}(q))^* \lambda_A^{\alpha, \omega}(y) \\
    &= \overline{\omega(q,p)} (\lambda_G^{\alpha,\omega}(qp))^* \lambda_A^{\alpha, \omega}(\alpha_q(x))\lambda_A^{\alpha, \omega}(y) \\
    &= \psi_{p\diamond q}(xy).
\end{align*}
Moreover, for $p \in P$ and $x,y \in X_p$, we have
\[\begin{split}
\psi(x)^* \psi(y)
    &= \big((\lambda_G^{\alpha,\omega}(p))^* \lambda_A^{\alpha, \omega}(x)\big)^*\big((\lambda_G^{\alpha,\omega}(p))^* \lambda_A^{\alpha, \omega}(y)\big) \\
    &= \lambda_A^{\alpha, \omega}(x)^* (\lambda_G^{\alpha,\omega}(p)) (\lambda_G^{\alpha,\omega}(p))^* \lambda_A^{\alpha, \omega}(y)
    = \lambda_A^{\alpha, \omega}(x^*y)
    = \psi_e(\langle x,y \rangle_A).
\end{split}\]
Each $X_p$ is essential, the left action of $A$ on each $X_p$
is injective and by compacts, and $P$ is directed.
Hence~\cite[Corollary~5.2]{SY} implies that $\psi$ is
CNP-covariant provided the condition $\psi^{(p)} \circ \phi_p =
\lambda_A^{\alpha, \omega}$ holds for all $p\in P$. To verify
this condition, fix an approximate identity $(e_k)_{k \in K}$
for $A$, and note that then $\phi_p(a)=\lim_{k \in K}
\alpha_p(a) \otimes e_k^*$ for $p\in P$ and $a\in P$. Hence
\[
\psi^{(p)}(\phi_p(a))
    = \lim_{k \in K} \psi_p(\alpha_p(a)) \psi_p(e_k)^*
    = \lim_{k \in K} \lambda_G^{\alpha, \omega}(p)^* \lambda_A^{\alpha, \omega}(\alpha_p(a) e_k^*)
    \lambda_G^{\alpha, \omega}(p)
    = \lambda_A^{\alpha, \omega}(a),
\]
showing that $\psi$ is CNP-covariant.

The image of $\psi$ generates $A \times^r_{\alpha,\omega} G$
because the latter is spanned by elements of the form
$\lambda_G^{\alpha,\omega}(g) \lambda_A^{\alpha, \omega}(a)$, and each $\lambda_G^{\alpha,\omega}(g) \in
C^*(\{\, \lambda_G^{\alpha,\omega}(p) \mid p \in P \,\})$ because $G$ is generated as a
group by $P$. Moreover, $\psi$ is injective as a representation because $\lambda_A^{\alpha, \omega}$
is automatically injective.

Since the left action on each fibre of $X$ is injective,
\cite[Lemma~4.3]{SY} implies that $X$ is
$\tilde\phi$-injective. Since the normalisation $\widehat{\alpha}^n$
of $\widehat{\alpha}$ of $G$ on $A \times^r_{\alpha,\omega} G$
satisfies $\widehat{\alpha}^n (\psi(x)) = \psi(x) \otimes
i_G(p)$ for all $p \in P$ and $x \in X_p$,
Theorem~\ref{thm:projective property} gives a surjective
homomorphism
$\phi \colon A \times^r_{\alpha,\omega} G \to \NO{X}^\reduced$
such that $\phi\circ\psi=j_{X}^r$.
Finally, the coaction $\widehat{\alpha}^n$ is normal by
definition, and so Corollary~\ref{cor:phi_inj} applies to give injectivity of $\phi$.
\end{proof}

\begin{cor}\label{cor:full CPs}
Resume the hypotheses of Lemma~\ref{lem:reduced CPs}. Then
there is an isomorphism $A \times_{\alpha,\omega} G \to \NO{X}$ which
takes $i_G^{\alpha,\omega}(p)^* i_A^{\alpha,\omega}(x)$ to
$j_{X }(x)$ for all $x \in X_p = A$.
\end{cor}
\begin{proof}
Since the isomorphism $\phi \colon A \times_{\alpha,\omega}^r G \to
\NO{X}^\reduced$ of Lemma~\ref{lem:reduced CPs}
intertwines the coactions on the two algebras, the
corresponding Fell bundles are isometrically isomorphic. Since
$A \times_{\alpha,\omega} G$ and $\NO{X}$ are the full
cross-sectional algebras of these Fell bundles (see
Remark~\ref{rmk:NO(X)=x-sect. alg}), the result follows.
\end{proof}

\begin{cor}\label{cor:NO(X)=C*(G)}
Let $(G,P)$ be a quasi-lattice ordered group such that $G$ is
generated as a group by $P$. Suppose that $P$ is directed. Let
$X$ be the product system over $P$ such that $X_p = \CC$ for
all $p \in P$ with all operations given by the usual operations
in $\CC$. Then $\NO{X} \cong C^*(G)$, and $\NO{X}^\reduced
\cong C^*_\reduced(G)$. Specifically, if $1_p$ denotes the element $1
\in \CC$ when regarded as an element of $X_p$, then $j_X(1_p)
\mapsto i_G(p)$ determines an isomorphism $\NO{X} \cong
C^*(G)$, and $j_X^\reduced(1_p) \mapsto \lambda_G(p)$
determines an isomorphism $\NO{X}^\reduced \cong C^*_\reduced(G)$.
\end{cor}
\begin{proof}
We apply Lemma~\ref{lem:reduced CPs} and
Corollary~\ref{cor:full CPs} to the trivial action $\alpha$ of
$P$ on $\CC$ and the trivial cocycle $\omega$ on $G$.
\end{proof}

\begin{rmk}
Let $G$ be a nonabelian finite-type Artin group and $P$ its
standard positive cone. By \cite[Proposition~29]{CL2002}, $P$
is directed and $G$ is not amenable. Therefore
Corollary~\ref{cor:NO(X)=C*(G)} implies that for the product
system considered there, we have
\[
\NO{X} \cong C^*(G) \not\cong C^*_\reduced(G) \cong \NO{X}^\reduced
\]
(cf. \cite[Theorem~30]{CL2002}). In particular $\NO{X}$ does
not have the gauge-invariant uniqueness property; so the
gauge-invariant uniqueness property is not automatic even for
systems where the left action is compact and injective and $P$
is directed. We thank Marcelo Laca for bringing this example to
our attention.
\end{rmk}

%\subsection*{Product systems over right-angled Artin semigroups}

%We pause here to observe that every product system over a
%finitely generated right-angled Artin group has the
%gauge-invariant uniqueness property.

%\begin{lemma}
%Let $G$ be a finitely generated right-angled Artin group, and
%let $P$ be its positive cone. For any product system over $P$,
%the $C^*$-algebra $\NO{X}$ has the gauge-invariant uniqueness
%property.
%\end{lemma}
%\begin{proof}
%We aim to apply Corollary~\ref{cor:amenable}(\ref{item:4}). To
%begin, we note that \cite[Theorem~13]{CN2005} implies that
%every finitely generated right-angled Artin group is exact (see
%\cite[p.~292]{CN2005} for the statement that their Theorem~13 applies
%to finitely generated right-angled Artin groups).

%Let $\{a_i \mid i \in I\}$ be the generators of $G$ in its
%presentation as a right-angled Artin group. Let $\{e_i \mid i \in
%I\}$ be the canonical generators for the free abelian group
%$\bigoplus_{i \in I} \ZZ$. The canonical homomorphism $\pi \colon G
%\to \ZZ^I$ satisfying $\pi(a_i) = e_i$ has all the properties
%set forth in Corollary~\ref{cor:amenable}(\ref{item:4}).
%Lemma~4.3 of \cite{SY} implies that $X$ is
%$\tilde\phi$-injective. Hence
%Corollary~\ref{cor:amenable}(\ref{item:4}) applies to prove the
%result.
%\end{proof}

\subsection*{Boundary quotient algebras}

The results in this section refer to the boundary quotient
algebras studied in~\cite{CL}. Throughout, given a
quasi-lattice ordered group $(G,P)$, we write $\Omega$ for the
Nica spectrum of $(G,P)$ and $\alpha$ for the partial action of
$G$ on $\Omega$ considered in \cite{CL} (see also \cite{ELQ}), and we let $\partial\Omega$
be the boundary of $\Omega$ defined in \cite{CL} and \cite{L}.

If $\alpha$ is a partial action of a discrete group $G$ on a $C^*$-algebra $A$,  \cite[Proposition 3.2]{QuiggRa}
shows that there is a canonical (dual) coaction $\widehat\alpha$ on the full partial crossed product
$A\rtimes_\alpha G$.
% such that $\widehat{\alpha}(am_s)=am_s \otimes i_G(s)$ for
%$a\in A$, $s\in G$ and $m_s$ the partial isometries in $(A\times_\alpha G)^{**}$ implementing the partial
%automorphisms $\alpha_s$.
Moreover, the discussion following \cite[Remark 3.7]{QuiggRa}
shows that the normalisation of $\widehat{\alpha}$ is naturally a coaction on the
reduced partial crossed product $A\rtimes_r G$.

\begin{lemma}\label{lem:boundary quotients}
Let $(G,P)$ be a quasi-lattice ordered group. Let $C_0(\partial
\Omega)\rtimes_r G$ be the reduced partial crossed product
corresponding to the partial crossed product $C_0(\partial
\Omega)\rtimes G$ considered in \cite{CL}, and let
$\beta \colon C_0(\partial \Omega)\rtimes_r G \to (C_0(\partial
\Omega)\rtimes_r G)\otimes C^*(G)$ be the canonical coaction of
$G$ on $C_0(\partial \Omega)\rtimes_r G$. Let $X$ be the
product system over $P$ such that $X_p=\CC$ for all $p$. Then
there is an isomorphism $\phi \colon C_0(\partial \Omega)\rtimes_r
G\to \NO{X}^\reduced$ such that $\redgaug\circ\phi=(\phi\otimes
\id_{C^*(G)}) \circ \beta$. In particular, $\phi|_{C_0(\partial
\Omega)}$ is an isomorphism from $C_0(\partial \Omega)$ to
$(\NO{X}^\reduced)_e^{\redgaug}$.
\end{lemma}
\begin{proof}
Let $(\pi,u)$ be the universal covariant representation of
$(C(\Omega),G,\alpha)$. By \cite[Proposition 6.1 and Theorem
6.4]{ELQ} the collection $\{\, u(p) \mid p\in P \,\}$ is a family of
isometries satisfying Nica's covariance relation which
generates $C(\Omega)\rtimes_\alpha G$. Since
$C_0(\partial\Omega) \rtimes_r G$ is a quotient of
$C(\Omega)\rtimes_\alpha G$, it follows that
$C_0(\partial\Omega) \rtimes_r G$ is generated by nonzero
isometries $\{\, W_p \mid p \in P \,\}$ (that these isometries are nonzero follows
from the fact that $\partial\Omega$ is nonempty, cf. \cite[Lemma 3.5]{CL} and
\cite[Theorem 3.7]{L}) satisfying Nica's covariance
relation such that the canonical coaction $\beta$ satisfies
$\beta(W_p) = W_p \otimes \lambda_G(p)$ for all $p$. For $p\in
P$ let $1_p$ denote the complex number $1$ regarded as an
element of $X_p$. The assignment $\psi\mapsto \{\, \psi(1_p) \mid
p\in P \,\}$ is a bijective correspondence between injective Nica
covariant Toeplitz representations of $X$, and families of
nonzero isometries satisfying Nica's covariance relation. Thus
$C_0(\partial\Omega) \rtimes_r G$ is generated by an injective
Nica covariant Toeplitz representation $\psi$ of $X$ satisfying
$\beta(\psi(x)) = \psi(x) \otimes i_G(d(x))$ for all $x \in X$.
Since the left action on each fibre of $X$ is implemented by an
injective homomorphism, Theorem \ref{thm:projective
property} gives a surjective homomorphism $\phi \colon
C_0(\partial \Omega)\rtimes_r G\to \NO{X}^\reduced$ such that
$\redgaug\circ\phi=(\phi\otimes \id_{C^*(G)}) \circ \beta$, and
we need only show that $\phi$ is injective.

By \cite[Lemma 3.5]{CL},  $\partial \Omega$ is the
unique minimal closed invariant subset of the Nica spectrum
$\Omega$, and since $\phi$ is nonzero, it follows that $\phi$
is injective on $C_0(\partial \Omega)$. Since $\beta$ is
normal, the expectation $\Phi^\beta \colon C_0(\partial
\Omega)\rtimes_r G \to C_0(\partial \Omega)$ is faithful. Since
$\redgaug\circ\phi=(\phi\otimes \id_{C^*(G)}) \circ \beta$,
it follows that $\phi$ intertwines the expectation $\Phi^\beta$
and the expectation from $\NO{X}^\reduced$ to
$(\NO{X}^\reduced)_e^{\redgaug    }$. The standard argument
now shows that $\phi$ is injective.
\end{proof}

\begin{cor}\label{cor:full bq}
Resume the hypotheses of Lemma~\ref{lem:boundary quotients}.
Let $\overline{\beta}$ be the canonical coaction on the
universal partial crossed product $C^*$-algebra $C_0(\partial
\Omega)\rtimes G$. There is an isomorphism $\phi \colon C_0(\partial
\Omega)\rtimes G\to \NO{X}$ such that
$\CNPgaug\circ\phi=(\phi\otimes \id_{C^*(G)}) \circ
\overline{\beta}$.
\end{cor}
\begin{proof}
We use the same trick as in the proof of
Corollary~\ref{cor:full CPs}.
\end{proof}

\begin{rmk}\label{rmk:BQ analysis}
The proofs of Lemma~\ref{lem:boundary quotients} and
Corollary~\ref{cor:full bq} are excellent examples of the
utility of the co-universal property of $\NO{X}^\reduced$. To
prove the same results otherwise one would first have to show
that $\NO{X}$ is isomorphic to the universal partial crossed
product of $C_0(\partial\Omega)$ by $G$ and then argue that the
normalisations of the two coactions on these universal
$C^*$-algebras yield  $\NO{X}^\reduced$ and $C_0(\partial
\Omega)\rtimes_{r} G$. Moreover, proving equality of universal
$C^*$-algebras would require verifying condition~(CP)
of~\cite{SY} in one direction, and the elementary relations
associated with the essential spectrum from~\cite{CL} in the
other direction --- for non-amenable $G$, there would be no
gauge-invariant uniqueness theorem to apply in either
direction.
\end{rmk}

\begin{rmk}
Combining Lemma~\ref{lem:boundary quotients} with
Lemma~\ref{lem:reduced CPs}, we see that if $(G,P)$ is a
quasi-lattice ordered group such that $G$ is generated as a
group by $P$ and each pair of elements in $P$ has a common
upper bound, then the boundary quotient algebra $C_0(\partial
\Omega) \rtimes_{r} G$ is isomorphic to the reduced
group $C^*$-algebra $C^*_\reduced(G)$; and under the same hypotheses,
Corollary~\ref{cor:full bq} combined with
Corollary~\ref{cor:full CPs} shows that the universal partial
crossed product $C^*$-algebra $C_0(\partial \Omega)
\rtimes G$ is isomorphic to the full group $C^*$-algebra
$C^*(G)$.
\end{rmk}

\subsection*{Topological higher-rank graphs}

In this section we show that each compactly aligned topological
higher-rank graph $\Lambda$ in the sense of Yeend can be used
to construct a compactly aligned product system $X$
of Hilbert bimodules over $C_0(\Lambda^0)$ with resulting $\Tc(X)$ isomorphic
to the $C^*$-algebra of Yeend's path groupoid and with corresponding
$\NO{X}$ isomorphic to the $C^*$-algebra of Yeend's boundary-path groupoid.

Recall \cite{y} that, for $k \in \NN$, a \emph{topological
$k$-graph} is a pair $( \Lambda, d )$ consisting of: (1)~ a
small category $\Lambda$ endowed with a second countable
locally compact Hausdorff topology under which the composition map
is continuous and open, the range map $r$ is continuous and the source
map $s$ is a local homeomorphism; and (2)~a continuous functor $d \colon \Lambda
\to \NN^{k}$, called the \emph{degree map}, satisfying the
factorisation property: if $d(\lambda) = m+n$ then there exist
unique $\mu,\nu$ with $d(\mu) = m$, $d(\nu) = n$ and $\lambda =
\mu\nu$.

Elements of $\Lambda$ are called paths, and paths of degree $0$
are called vertices. For $m \in \NN^{k}$ we define $\Lambda^{m}
:= d^{-1} (m)$. If $0 \le m \le n \le p$ in $\NN^{k}$ and
$\lambda \in \Lambda^{p}$ then we write $\lambda (m,n)$ for the
unique path in $\Lambda^{n-m}$ such that $\lambda = \mu \lambda
(m,n) \nu$, where $\mu \in \Lambda^{m}$ and $\nu \in
\Lambda^{p-n}$. For $0 \le m \le p$ in $\NN^{k}$ and $\lambda
\in \Lambda^{p}$ we write $\lambda (m)$ for $s (\lambda (0,m))
= \lambda (m,m)$. If $m,p \in \NN^{k}$ with $m \le p$ then the
map $\sigma^{m} \colon \Lambda^{p} \to \Lambda^{p-m}$ such that
$\sigma^{m} (\lambda) = \lambda (m,p)$ is continuous. For $U, V
\subset \Lambda$, we write
\[
UV := \{\, \lambda\mu \mid \lambda\in U,\ \mu\in V,\ s(\lambda)=r(\mu)\,\}.
\]
For $U \subset \Lambda^{p}$ and $V \subset \Lambda^{q}$,
\[
U \vee V := U \Lambda^{(p \vee q)-p} \cap V \Lambda^{(p \vee q)-q}
\]
is the set of \emph{minimal common extensions} of paths from
$U$ and $V$. A topological $k$-graph $( \Lambda, d )$ is
\emph{compactly aligned} if $U \vee V$ is compact whenever $U$
and $V$ are compact.

Fix $k \in \NN$, and let $( \Lambda, d)$ be a topological
$k$-graph. Define $A := C_{0} ( \Lambda^{0} )$. For each $n \in
\NN^{k}$ let $X_{n}$ be the right-Hilbert $A$-$A$ bimodule
associated to the topological graph $( \Lambda^{0},
\Lambda^{n}, s |_{\Lambda^{n}}, r |_{\Lambda^{n}} )$ (see
\cite{ka1}). So $X_n$ is the completion of the pre-Hilbert
$A$--$A$ bimodule $C_c(\Lambda^n)$ with operations
\[
\langle f, g \rangle_{A}^{n} ( v ) = \textstyle{\sum_{\eta \in \Lambda^{n}v}} \overline{f ( \eta )} g (\eta)
\quad\text{ and }\quad
( a \cdot f \cdot b) ( \lambda ) = a (r(\lambda))
f(\lambda) b(s(\lambda)).
\]
Katsura shows that $X_n$ is a subspace of $C_0(\Lambda^n)$
\cite[Section~1]{ka1}.

\begin{prop}\label{prp:prod sys}
Let $\Lambda$ be a topological $k$-graph and let $(X_n)_{n\in\NN^k}$ be as above.
For $f \in X_{m}$ and $g \in X_{n}$, define $fg \colon
\Lambda^{m+n} \to \CC$ by $(fg)(\lambda) := f(\lambda (0,m))
g(\lambda (m,m+n))$. Under this multiplication the family $X :=
\bigsqcup_{n \in \NN^{k}} X_{n}$ of right-Hilbert
$C_0(\Lambda^0)$-$C_0(\Lambda^0)$ bimodules
is a product system over $\NN^{k}$.
\end{prop}

\begin{proof}
We first show that for $f_1, f_2 \in X_{m}$ and $g_1, g_2 \in
X_{n}$ we have $\langle f_1g_1, f_2g_2 \rangle_{A}^{m+n} =
\langle g_1, \langle f_1,f_2 \rangle_{A}^{m} \cdot g_2
\rangle_{A}^{n}$. The functions $f_1g_1, f_2g_2$ are continuous
on $\Lambda^{m+n}$ and for $v \in \Lambda^{0}$,
\begin{align*}
\langle f_1g_1, f_2g_2 \rangle_{A}^{m+n} ( v )
&= \textstyle{\sum_{\lambda \in \Lambda^{m+n} v}}  \overline{f_1(\lambda(0,m))} f_2(\lambda(0,m))
            \overline{g_1(\lambda(m,m+n))} g_2(\lambda(m,m+n))
\\&= \textstyle{\sum_{\nu \in \Lambda^{n} v}} \big( \textstyle{\sum_{\mu \in \Lambda^{m} r( \nu )}}
            \overline{f_1(\mu)} f_2(\mu) \big) \overline{g_1(\nu)} g_2(\nu)
\\&= \textstyle{\sum_{\nu \in \Lambda^{n} v}} \langle f_1, f_2 \rangle_{A}^{m} ( r( \nu ) ) \overline{g_1(\nu)} g_2(\nu)
\\&= \langle g_1, \langle f_1,f_2 \rangle_{A}^{m} \cdot g_2 \rangle_{A}^{n} ( v ).
\end{align*}
Taking $f_2 = f_1 = f$ and $g_2 = g_1 = g$, we deduce that $fg
\in X_{m+n}$ by definition of $X_{m+n}$ (see \cite{ka1}).
Further, since $\langle g_1, \langle f_1,f_2 \rangle_{A}^{m}
\cdot g_2 \rangle_{A}^{n} = \langle f_1 \otimes_{A} g_1, f_2
\otimes_{A} g_2 \rangle_{A}$, the map $f \otimes_{A} g \mapsto
fg$ extends to an isometric adjointable operator from $X_{m}
\otimes_{A} X_{n}$ to $X_{m+n}$. To show that it is surjective
it is enough to show that it has dense range. The subset $\Aa
:= \lsp \{\, fg \mid f \in C_{c} (\Lambda^{m}), g \in C_{c}
(\Lambda^{n}) \,\}$ of $C_{c} (\Lambda^{m+n})$ is a subalgebra
of $C_{0} (\Lambda^{m+n})$. For every open subset $U$ of
$\Lambda^{m+n}$, an application of the Stone-Weierstrass
Theorem shows that $\Aa \cap C_{c} (U)$ is uniformly dense in
$C_{c} (U)$. So $\Aa$ is dense in $X_{m+n}$ with respect to $\|
\cdot \|_{A}$ by \cite[Lemma 1.26]{ka1}, and the multiplication
operator has dense range.
\end{proof}

To apply our results, we need to show that the product system
$X$ is compactly aligned if $\Lambda$ is compactly aligned.
We first need a couple of technical lemmas.

\begin{notation}\label{ntn:Fms}
For $m \in \NN^{k}$ we denote by $F_{m}$ the set of functions
$f \in C_{c} (\Lambda^{m})$ such that the source map restricts
to a homeomorphism of $\supp(f)$. By definition of $F_m$, for
each $f \in F_m$ and $v \in \Lambda^0$ such that $\Lambda^{m} v$
is non-empty we may fix an element, henceforth denoted $\lambda_{f,v}$,
of $\Lambda^m v$ such that $f(\mu) = 0$ for all $\mu \in \Lambda^m v \setminus
\{\lambda_{f,v}\}$.
\end{notation}

\begin{lemma}\label{dense}
For $m \in \NN^{k}$, $\lsp F_{m}$ is dense in $X_{m}$ with
respect to the norm $\| \cdot \|_{A}$.
\end{lemma}
\begin{proof}
A partition of unity argument using that the source map in
$\Lambda$ is a local homeomorphism shows that each element of
$C_{c} (\Lambda^{m})$ is a finite sum of elements of $F_m$,
and the result follows.
\end{proof}

\begin{lemma}\label{tensor}
Fix $m,n \in \NN^{k}$. Then, for $f \in X_{m}$, $g \in F_{m}$, $c \in X_{m
\vee n}$, and $\xi \in \Lambda^{m \vee n}$ we have
\[
\big(\iota_{m}^{m \vee n} (f \otimes g^{*}) (c)\big)(\xi)
=
f(\xi(0,m)) \overline{g(\lambda_{g, \xi(m)})} c(\lambda_{g, \xi(m)}\xi (m,m \vee n)).
\]
\end{lemma}
\begin{proof}
For $h \in X_{m}$ and $l \in X_{(m \vee n)-m}$ we have
\begin{align*}
\big(\iota_{m}^{m \vee n} (f \otimes g^{*}) (hl)\big)(\xi)
    &= f(\xi (0,m)) \big(\textstyle{\sum_{\mu \in \Lambda^{m} \xi (m)}} \overline{g(\mu)} h(\mu)\big) l(\xi (m,m \vee n)) \\
    &= f(\xi (0,m)) \overline{g(\lambda_{g, \xi(m)})} h(\lambda_{g, \xi(m)}) l(\xi(m,m \vee n)) \\
    &= f(\xi (0,m)) \overline{g(\lambda_{g, \xi(m)})} hl(\lambda_{g, \xi(m)}\xi(m,m \vee n)).
\end{align*}
Since elements of the form $hl$ are dense in $X_{m \vee n}$,
the result follows.
\end{proof}

\begin{lemma}[{cf. \cite[Lemma 5.1]{RS}}] \label{lemma:compact}
  Fix $m\in\NN^k$ and $T\in\Kk(X_m)$. Let $B_1F_m:=\{f\in F_m\mid ||f||_\infty\le 1\}$.
  Then the function $\chi_T:\Lambda^m\to\RR$ defined by $\chi_T(\lambda)
  =\sup_{f\in B_1F_m}\abs{T(f)(\lambda)}$ vanishes at infinity on $X_m$.
\end{lemma}

\begin{proof}
  Since
  $\chi_{\alpha S+\beta T}(\lambda)\le
  \abs{\alpha}\chi_S(\lambda)+\abs{\beta}\chi_T(\lambda)$,
  and $\abs{T(f)(\lambda)}\le \norm{T}\ \norm{f}_{X_m}$ for
  $\alpha,\beta\in \CC$, $S,T\in \Kk(X_m)$,
  $f\in F_m$ and $\lambda\in\Lambda^m$, it is enough to prove that $\chi_T$ vanishes at
  infinity when $T=g\otimes h^*$ for $g,h\in X_m$.

  If $g,h\in X_m$, $f\in B_1F_m$ and $\lambda\in\Lambda^m$, then we have
  \begin{equation*}
    \abs{(g\otimes h^*)(f)(\lambda)}=\bigl\lvert g(\lambda)
      \textstyle{\sum_{\eta\in \Lambda^ms(\lambda)}}
      \overline{h(\eta)}f(\eta)\bigr\rvert
    =\abs{g(\lambda)\overline{h(\lambda_{f,s(\lambda)})} f(\lambda_{f,s(\lambda)})}
    \le \abs{g(\lambda)}\norm{h}_\infty
  \end{equation*}
  and since $g$ vanishes at infinity on $\Lambda^m$, so will $\chi_{g\otimes h^*}$.
\end{proof}

Most of the work in proving that $X$ is compactly
aligned when $\Lambda$ when is compactly aligned,
is involved in proving the following technical lemma.
We state this lemma separately because we will use it again to
prove Proposition~\ref{prop:THRG NicaCov}.

\begin{lemma}\label{compact}
Assume that $\Lambda$ is compactly aligned.
Fix $m,n \in \NN^{k}$. Let $f_{m} \in F_{m}$ and $f_{n} \in
F_{n}$. Then $C := \supp (f_{m}) \vee \supp (f_{n}) \subset
\Lambda^{m \vee n}$ is compact. For each of $p=m,n$: let $\{\,
V_{i}^{p} \mid 1 \le i \le r_{p} \,\}$ be a finite open cover of
$\sigma^{p} (C)$ such that each $\overline{V_{i}^{p}}$ is
compact and $s$ restricts to a homeomorphism on each
${V_{i}^{p}}$; let $\phi_{i}^{p} \colon \sigma^{p} (C) \to
[0,1]$, $1 \le i \le r_{p}$, be a partition of unity on $C(C)$
subordinate to $\{\, V_{i}^{p} \cap \sigma^{p} (C) \mid 1 \le i
\le r_{p} \,\}$; and fix functions $\rho_{i}^{p} \colon
\Lambda^{(m \vee n)-p} \to [0,1]$ such that each
$\rho_i^p|_{\sigma^{p} (C)} = \sqrt{\phi_i^p}$, and each
$\rho_i^p$ vanishes off $V_i^p$. Fix $g_{m} \in F_{m}$ and
$g_{n} \in F_{n}$. For $1 \le i \le r_{m}$ and $1 \le j \le
r_{n}$, define $a_{ij}, b_{j} \in C_{c} (\Lambda^{m \vee n})$
by $a_{ij} := g_{m}(\rho_{i}^{m} \cdot \langle f_{m}
\rho_{i}^{m},f_{n} \rho_{j}^{n} \rangle_{A}^{m \vee n})$ and
$b_{j} := g_{n} \rho_{j}^{n}$. Then
\begin{equation}\label{compactprod}
\iota_{m}^{m \vee n} (g_{m} \otimes f_{m}^{*}) \iota_{n}^{m \vee n} (f_{n} \otimes g_{n}^{*})
= \textstyle{\sum_{i=1}^{r_{m}} \sum_{j=1}^{r_{n}}} a_{ij} \otimes b_{j}^{*}.
\end{equation}
\end{lemma}
\begin{proof}
Since $\Lambda$ is compactly aligned, $C$ is compact. Since the
source map is injective on each $\supp(\rho_i^m)$, for $\mu \in
\Lambda^{(m \vee n)-m}$ such that $\Lambda^{m} r(\mu)$ is non-empty and $1 \le j \le r_{n}$ we have
\begin{align*}
(\textstyle{\sum_{i=1}^{r_{m}}}& \rho_{i}^{m} \cdot \langle f_{m} \rho_{i}^{m},f_{n} \rho_{j}^{n} \rangle_{A}^{m \vee n}) (\mu) \\
    &= \textstyle{\sum_{i=1}^{r_{m}}} \rho_{i}^{m} (\mu) \sum_{\alpha \in \Lambda^{m \vee n} s(\mu)}
        \overline{f_{m} (\alpha (0,m))} \rho_{i}^{m} (\alpha (m,m \vee n)) (f_{n} \rho_{j}^{n})(\alpha) \\
    &= \textstyle \sum_{i=1}^{r_{m}} (\rho_{i}^{m} (\mu))^{2} \overline{f_{m}(\lambda_{f_m, r(\mu)})} (f_{n} \rho_{j}^{n})(\lambda_{f_m, r(\mu)}\mu) \\
    &= \overline{f_{m}(\lambda_{f_m, r(\mu)})} (f_{n} \rho_{j}^{n})(\lambda_{f_m, r(\mu)} \mu).
\end{align*}
Fix $c \in X_{m \vee n}$ and $\xi \in \Lambda^{m \vee n}$, and
write $\lambda_m := \lambda_{f_m, \xi (m)}$, $\lambda' := \xi(m, m \vee n)$, and $\beta := (\lambda_m \lambda')(n,m \vee n)$.
If $\beta \in \sigma^{n} (C)$ then, since the source map is injective on each $\supp
(\rho_{j}^{n})$, we have
\begin{align*}
\textstyle{\sum_{j=1}^{r_{n}}} \rho_{j}^{n} (\beta) \langle b_{j},c \rangle_{A}^{m \vee n} (s(\xi))
    &= \textstyle{\sum_{j=1}^{r_{n}}} \rho_{j}^{n} (\beta)
        \sum_{\lambda \in \Lambda^{m \vee n} s(\xi)} \overline{g_{n} (\lambda (0,n))} \rho_{j}^{n} (\lambda (n,m \vee n)) c (\lambda) \\
    &= \textstyle{\sum_{j=1}^{r_{n}}} (\rho_{j}^{n} (\beta))^{2} \overline{g_{n}(\lambda_{g_n, r(\beta)})} c(\lambda_{g_n, r(\beta)} \beta) \\
    &= \overline{g_{n}(\lambda_{g_n, r(\beta)})} c(\lambda_{g_n, r(\beta)} \beta).
\end{align*}
Hence we have
\begin{align*}
\big(\textstyle \sum_{i=1}^{r_{m}}&\textstyle \sum_{j=1}^{r_{n}} a_{ij} \otimes b_{j}^{*}\big) (c)(\xi) \\
    &\;= \textstyle{\sum_{i=1}^{r_{m}} \sum_{j=1}^{r_{n}} g_{m}(\xi (0,m))
        (\rho_{i}^{m} \cdot \langle f_{m} \rho_{i}^{m},f_{n} \rho_{j}^{n} \rangle_{A}^{m \vee n}) (\lambda')
        \langle b_{j},c \rangle_{A}^{m \vee n} (s(\xi))} \\
    &\;= g_{m}(\xi (0,m)) \textstyle{\sum_{j=1}^{r_{n}}
        (\sum_{i=1}^{r_{m}} \rho_{i}^{m} \cdot \langle f_{m} \rho_{i}^{m},f_{n} \rho_{j}^{n} \rangle_{A}^{m \vee n}) (\lambda')
        \langle b_{j},c \rangle_{A}^{m \vee n} (s(\xi))} \\
    &\;= g_{m}(\xi (0,m)) \overline{f_{m}(\lambda_m)} f_{n}((\lambda_m\lambda') (0,n))
        \overline{g_{n}(\lambda_{g_n, r(\beta)})} c(\lambda_{g_n, r(\beta)} \beta)
\end{align*}
which equals $\big(\iota_{m}^{m \vee n} (g_{m} \otimes
f_{m}^{*}) \iota_{n}^{m \vee n} (f_{n} \otimes g_{n}^{*})\big)
(c)(\xi)$ by two applications of Lemma~\ref{tensor}.
\end{proof}

The next result generalises \cite[Theorem 5.4]{RS}.

\begin{prop}\label{CA}
Let $\Lambda$ be a topological $k$-graph.
The product system $X$ defined by Proposition \ref{prp:prod sys}
is compactly aligned if and only if $\Lambda$ is
compactly aligned.
\end{prop}

\begin{proof}
Assume that $\Lambda$ is compactly aligned.
Fix $m,n \in \NN^{k}$. By Lemma~\ref{dense} it suffices to show
that, for $f,g \in F_{m}$ and $h,l \in F_{n}$, $\iota_{m}^{m
\vee n} (f \otimes g^{*}) \iota_{n}^{m \vee n} (h \otimes l^{*}) \in
\Kk (X_{m \vee n})$, and this follows from Lemma~\ref{compact}.

Assume next that $\Lambda$ is not compactly aligned. Then there exist $m,n\in\NN^k$
and $U\subset\Lambda^m$, $V\subset\Lambda^n$ such that $U$ and $V$ are compact, but $U\vee V$ is not compact. For each of $C=U,V$: let $\{\,
V_{i}^C \mid 1 \le i \le r_C \,\}$ be a finite open cover of
$C$ such that each $\overline{V_{i}^C}$ is
compact and $s$ restricts to a homeomorphism on each
${V_{i}^C}$; let $\phi_{i}^C \colon C \to
[0,1]$, $1 \le i \le r_C$, be a partition of unity on $C(C)$
subordinate to $\{\, V_{i}^C \cap C \mid 1 \le i
\le r_C \,\}$; and fix functions $\rho_{i}^C\in C_C(\Lambda^p)$, where $p=m$ if $C=U$
and $p=n$ if $C=V$, such that each $\rho_i^C|_{C} = \sqrt{\phi_i^C}$, and each
$\rho_i^C$ vanishes off $V_i^C$. Let $T_C=\sum_{i=1}^{r_C}\rho_i^C\otimes(\rho_i^C)^*$.
Then $T_U\in\Kk(X_m)$ and $T_V\in\Kk(X_n)$, but we claim that
$T:=\iota_m^{m\vee n}(T_U) \iota_n^{m\vee n}(T_V)$ is not compact which will show that
$X$ is not compactly aligned. Notice first that if $f\in X_p$ and $\lambda\in C$, then
we have
\begin{equation} \label{eq:1}
  \begin{split}
    T_C(f)(\lambda) &= \sum_{i=1}^{r_C}\bigl(\rho_i^C\otimes (\rho_i^C)^*\bigr)(f)(\lambda) \\
    &= \sum_{i=1}^{r_C}\rho_i^C(\lambda)
    \sum_{\eta\in\Lambda^ps(\lambda)}\overline{\rho_i^C(\eta)}f(\eta)
    =\sum_{i=1}^{r_C}\rho_i^C(\lambda) \overline{\rho_i^C(\lambda)}f(\lambda)
    =f(\lambda).
  \end{split}
\end{equation}
For each $\lambda\in U\vee V$ choose $f_\lambda\in B_1F_{m\vee n}$ such that $f_\lambda(\lambda)=1$.
Equation \eqref{eq:1} implies that
\begin{equation*}
  T(f_\lambda)(\lambda)=\iota_m^{m\vee n}(T_U) \iota_n^{m\vee n}(T_V)(f_\lambda)(\lambda)
  =\iota_n^{m\vee n}(T_V)(f_\lambda)(\lambda)
  =f_\lambda(\lambda)=1,
\end{equation*}
so it follows from Lemma \ref{lemma:compact} that $T$ is not compact.
\end{proof}

% We may now apply our results to see that to the compactly aligned
% topological higher-rank graph $\Lambda$ there is associated a
% $C^*$-algebra $\NO{X}$ which is generated by isometric
% copies of the $C_c(\Lambda^n)$ and is co-universal amongst all
% such $C^*$-algebras carrying a gauge-action of $\TT^k$. In
% particular, Corollary~\ref{cor:Nk giut} applies, and suggests
% that this $\NO{X}$ should be regarded as \emph{the}
% topological higher-rank graph $C^*$-algebra of $\Lambda$.
In~\cite{y}, Yeend associated two groupoids $G_\Lambda$ and
$\mathcal{G}_\Lambda$ and hence two $C^*$-algebras
$C^*(G_\Lambda)$ and $C^*(\mathcal{G}_\Lambda)$ to each compactly
aligned topological higher-rank graph $\Lambda$, and proposed $C^*(G_\Lambda)$
as a model for the Toeplitz algebra of $\Lambda$, and
$C^*(\mathcal{G}_\Lambda)$ as the Cuntz-Krieger algebra.
%It is therefore natural to ask how $\NO{X}$ relates to the two.
We will  show that $\Tc(X)$ is isomorphic to $C^*(G_\Lambda)$ and
that $\NO{X}$ is isomorphic to $C^*(\mathcal{G}_\Lambda)$.

In the following, we use the notation for paths in topological
$k$-graphs established in \cite[Lemma 3.3]{y}. We denote by
$G_{\Lambda}$ the \emph{path groupoid} associated to $\Lambda$
\cite[Definition 3.4]{y}. So $G_\Lambda$ consists of triples
$(x,m,y)$ where $x$ and $y$ are (possibly infinite) paths in
$\Lambda$, $m \in \ZZ^{k}$, and there exist $p,q \in \NN^{k}$
such that $p \le d(x)$, $q \le d(y)$, $p-q=m$ and $\sigma^p(x) = \sigma^q(y)$.
By \cite[Theorem 3.16]{y}, $G_{\Lambda}$ is a locally compact
$r$-discrete topological groupoid admitting a Haar system
consisting of counting measures. A basis for the topology on $G_{\Lambda}$
is as follows \cite[Proposition 3.6]{y}. Define $\Lambda \ast_{s} \Lambda :=
\{\, (\lambda, \mu) \in \Lambda \times \Lambda \mid s(\lambda)=s(\mu) \,\}$,
and for $U,V \subset \Lambda$ define $U \ast_{s} V := (U \times V) \cap (\Lambda \ast_{s} \Lambda)$.
For $F \subset \Lambda \ast_{s} \Lambda$ and $m \in \ZZ^{k}$, define $Z(F,m) :=
\{\, (\lambda x,d(\lambda)-d(\mu),\mu x) \in G_{\Lambda} \mid (\lambda,\mu) \in F, d(\lambda)-d(\mu)=m\,\}$.
Then the family of sets of the form $Z(U \ast_{s} V,m) \cap Z(F,m)^{c}$, where $m \in \ZZ^{k}$,
$U,V \subset \Lambda$ are open and $F \subset \Lambda \ast_{s} \Lambda$ is compact,
is a basis for the topology on $G_{\Lambda}$.

Let $V\subset \Lambda^0$. A set $E\subset V\Lambda$ is called \emph{exhaustive}
(cf. \cite[Definition 4.1]{y}) for $V$ if for all $\lambda\in V\Lambda$ there exists $\mu\in E$ such that $\{\lambda\}\vee\{\mu\}\ne\emptyset$. For $v\in\Lambda^0$, let $v\mathcal{CE}(\Lambda)$ denote the set of all compact sets $E\subset \Lambda$ such that $r(E)$ is a neighbourhood of $v$ and $E$ is exhaustive for $r(E)$. A (possibly  infinite) path $x$ is called a \emph{boundary path} (cf. \cite[Definition 4.2]{y}) if for all $m\in\NN^k$ with $m\le d(x)$, and for all $E\in x(m)\mathcal{CE}(\Lambda)$, there exists $\lambda\in E$ such that $x(m,m+d(\lambda))=\lambda$. We write $\partial \Lambda$ for the set of all boundary paths. It is shown in \cite{y} that $\partial\Lambda$ is a closed and invariant subset of $G_\Lambda^{(0)}$ (we are here identifying a path $x$ with the element $(x,0,x)$ in $G_\Lambda^{(0)}$) and that $v\partial\Lambda\ne \emptyset$ for all $v\in\Lambda^0$. The \emph{boundary-path groupoid} $\mathcal{G}_\Lambda$ is then defined in \cite[Definition 4.8]{y} to be the reduction of $G_\Lambda$ to $\partial\Lambda$. We will now show that $\partial\Lambda$ is in fact the smallest closed and invariant subset $Y$ of $G_\Lambda^{(0)}$ such that $vY$ is nonempty for all $v\in\Lambda^0$. Let $X_\Lambda$ denote the
  collection of finite and infinite paths in $\Lambda$.  For $V\subset \Lambda$ write
  $Z(V):=\{x\in X_\Lambda\mid\text{there exists } \lambda\in V\text{ such that }x(0,d(\lambda))=\lambda\}$.

\begin{prop} \label{prop:small}
  Let $\Lambda$ be a compactly aligned topological $k$-graph. Then $\partial\Lambda$
  is the smallest closed and invariant subset $Y$ of $G_\Lambda^{(0)}$ such that $vY$ is
  nonempty for all $v\in\Lambda^0$.
\end{prop}

\begin{proof}
  It follows from \cite[Propositions 4.3, 4.4 and 4.7]{y} that $\partial\Lambda$ is a closed
  and invariant subset of $G_\Lambda^{(0)}$ and that $v\partial\Lambda\ne \emptyset$ for
  all $v\in\Lambda^0$. We will show that $\partial\Lambda$ is contained in any other
  closed and invariant subset $Y$ of $G_\Lambda^{(0)}$ satisfying that $vY$ is
  nonempty for all $v\in\Lambda^0$. So let $Y$ be such a subset and assume for
  contradiction that there is a boundary path $x$ such that $x\notin Y$. Since $Y$ is
  closed, it follows from \cite[Lemma 3.8]{y} that there is a relatively compact and open
  subset $U$ of $\Lambda$ and a compact subset $F$ of $\Lambda$ satisfying that
  $\mu\in F$ implies that $\mu=\lambda\mu'$ for some $\lambda\in\overline{U}$ such
  that $x\in Z(U)\setminus Z(F)\subset X_\Lambda\setminus Y$. We will
  show that there is a $\lambda\in \Lambda$ such that
  $\lambda X_\Lambda \subset Z(U)\setminus Z(F)\subset X_\Lambda\setminus Y$.
  It will then follow from the invariance of $Y$ that $s(\lambda)Y=\emptyset$,
  and we have our contradiction.

  Choose $m\in\NN^k$ such that $x(0,m)\in U$ and a compact neighbourhood
  $V\subset U\cap\Lambda^m$ of $x(0,m)$ such that $s$ restricted to $V$ is injective. We let
  $$C:=\sigma^m((V\vee F)\cap F)
  =\{\mu\in\Lambda\mid\text{there exists }\lambda\in V\text{ such that }\lambda\mu\in F\}.$$
  Since $\Lambda$  is compactly aligned and  $\sigma^m$ is continuous, the set $C$ is compact.
  It follows from the assumption $x\in\partial\Lambda$ that if
  $C\in x(m)\mathcal{CE}(\Lambda)$, then there  exists $\mu\in C$ such that
  $x(m,m+d(\mu))\in C$, and then $x(0,m+d(\mu))=x(0,m)\mu \in F$,
   a contradiction. Thus either
  $r(C)$ is not a neighbourhood of $x(m)$ or $C$ is not exhaustive for $r(C)$.

  If $r(C)$ is not a neighbourhood of $x(m)$, then since $s(V)$ is a neighbourhood of
  $x(m)$ there exists $\lambda\in V$ such that $s(\lambda)\ne r(\mu)$ for all
  $\mu\in C$, and then $\lambda X_\Lambda\subset Z(V)\setminus Z(F)\subset
  Z(U)\setminus Z(F)$.

  If $C$ is not exhaustive for $r(C)$, then there exists $\mu\in r(C)\Lambda$ such that
  $\mu\Lambda\cap C\Lambda=\emptyset$. Let $\lambda:=\eta\mu$ where $\eta$
  is the unique element of $V$ such that $s(\eta)=r(\mu)$. Then
  $\lambda X_\Lambda\subset Z(V)\setminus Z(F)\subset
  Z(U)\setminus Z(F)$.
\end{proof}

\begin{prop} \label{prop:thrg-rep}
Let $\Lambda$ be a compactly aligned topological $k$-graph.
Let $X$ be the product system constructed in Proposition \ref{prp:prod sys}.
There is a unique Toeplitz representation $\psi \colon X \to C^{*} ( G_{\Lambda} )$
such that, for $n \in \NN^{k}$ and $f \in C_{c} (\Lambda^{n})$, we have $\psi (f) \in C_{c} (G_{\Lambda})$
and
\begin{equation}\label{TR}
\psi (f) ((x,p,y))
=
\begin{cases}
f(x(0,n)) &\text{if $p = n$ and $\sigma^{n} (x) = y$}\\
0 &\text{otherwise.}
\end{cases}
\end{equation}
Moreover, the $\psi (f)$ for $n \in \NN^{k}$ and $f \in C_{c} (\Lambda^{n})$
generate $C^{*} ( G_{\Lambda} )$.
\end{prop}
\begin{proof}
We first show that, for $f \in C_{c} ( \Lambda^{n} )$, $\psi
(f)$ defined as in \eqref{TR} is in $C_{c} (G_{\Lambda})$. To see that $\psi (f)$ is
continuous, fix $(x,p,y) \in G_{\Lambda}$ and $\epsilon > 0$.
If $p \ne n$ then $Z(\Lambda \ast_{s} \Lambda,p)$ is an open
neigbourhood of $(x,p,y)$ on which $\psi (f)$ is zero. If $p = n$
and $\sigma^{n} (x) \ne y$ then, since $U \ast_{s} s(U)$ is
compact where $U := \supp(f)$, the subset $Z(\Lambda \ast_{s}
\Lambda,n) \cap Z(U \ast_{s} s(U),n)^{c}$ is an open
neighbourhood of $(x,p,y)$ on which $\psi (f)$ is zero. Suppose
now that $p = n$ and $\sigma^{n} (x) = y$. Since
$f$ is continuous there exists an open neighbourhood $U \subset
\Lambda^{n}$ of $x(0,n)$ such that $\lambda \in U$ implies $|
f(\lambda) - f(x(0,n)) | < \epsilon$. Note that, since
$\Lambda^{n}$ is open in $\Lambda$, $U$ is open in $\Lambda$.
Choose an open neighbourhood $V$ of $x(0,n)$ in $\Lambda$ such
that $V \subset U$, $s(V)$ is open in $\Lambda$, and $s|_{V}$
is a homeomorphism onto $s(V)$. Then $Z(V \ast_{s} s(V),n)$ is
an open neighbourhood of $(x,p,y)$ such that $(w,n,z) \in Z(V \ast_{s}
s(V),n)$ implies $\sigma^{n} (w) = z$ and $w(0,n) \in V \subset
U$, hence $| \psi (f) ((w,n,z)) - \psi (f) ((x,p,y)) | = |
f(w(0,n)) - f(x(0,n)) | < \epsilon$. It follows that $\psi (f)
\in C(G_{\Lambda})$.

To see that $\psi (f)$ has compact support, let $U := \supp
(f)$, then $Z(U \ast_{s} s(U),n)$ is compact by \cite[Proposition 3.15]{y}
and $\supp (\psi (f)) \subset Z(U \ast_{s} s(U),n)$.

We now show that $\| \psi (f) \|_{C^{*} (G_{\Lambda})} \le \| f
\|_{A}$ for $f \in C_{c} (\Lambda^{n})$. By \cite[II.1.5]{r} it
suffices to show that
\begin{equation}\label{int1}
\mbox{$\displaystyle{\sup_{(y,0,y) \in G_{\Lambda}^{(0)}}} \textstyle{\int_{G_\Lambda} | \psi (f)^{*} \ast \psi (f) ((x,p,z)) | \ d \lambda^{(y,0,y)} ((x,p,z))}
\le \| f \|_{A}^{2}$\quad and}
\end{equation}
\begin{equation}\label{int2}
\sup_{(y,0,y) \in G_{\Lambda}^{(0)}} \textstyle{\int_{G_\Lambda} | \psi (f)^{*} \ast \psi (f) ((x,p,z)^{-1}) | \ d \lambda^{(y,0,y)} ((x,p,z))}
\le \| f \|_{A}^{2}.
\end{equation}
Fix $(y,p,z) \in G_\Lambda$. Since the Haar system on
$G_{\Lambda}$ consists of counting measures,
\begin{align*}
\psi (f)^{*} \ast \psi (f) ((y,p,z))
&= \textstyle{\sum_{(y,r,w) \in G_{\Lambda}} \overline{\psi (f) ((w,-r,y))} \psi (f) ((w,p-r,z))}
\\&= \textstyle{\sum_{\alpha \in \Lambda^{n} r(y)} \overline{f(\alpha)} \psi (f) ((\alpha y,p+n,z))}
\\&=
\begin{cases}
\langle f,f \rangle_{A}^{n} (r(y)) &\text{if $p=0$ and $z=y$}\\
0 &\text{otherwise.}
\end{cases}
\end{align*}
So, for $(y,0,y) \in G_{\Lambda}^{(0)}$, we have
\begin{align*}
 \textstyle \int_{G_\Lambda} | \psi (f)^{*} \ast \psi (f) ((x,p,z)) | &\ d \lambda^{(y,0,y)} ((x,p,z)) \\
    &= | \psi (f)^{*} \ast \psi (f) ((y,0,y)) |
    = \langle f,f \rangle_{A}^{n} (r(y))
    \le \| f \|_{A}^{2},
\end{align*}
and this establishes~(\ref{int1}). A similar argument
establishes~(\ref{int2}).

Straightforward calculations show that $\psi \colon C_{c}
(\Lambda^{n}) \to C^{*} (G_{\Lambda})$ is linear and is
multiplicative when $n=0$. Since $\| \psi (f) \|_{C^{*}
(G_{\Lambda})} \le \| f \|_{A}$ for $f \in C_{c}
(\Lambda^{n})$, $\psi$ extends to a linear map $\psi$ from
$X_{n}$ to $C^{*} (G_{\Lambda})$, and $\psi \colon A \to C^{*}
(G_{\Lambda})$ is a homomorphism.

We show that $\psi$ is multiplicative. It suffices to show that
$\psi (fg) = \psi (f) \ast \psi (g)$ for $f \in C_{c}
(\Lambda^{m})$ and $g \in C_{c} (\Lambda^{n})$. Indeed, for $(x,p,y)
\in G_{\Lambda}$ we have
\begin{align*}
\psi (f) \ast \psi (g) ((x,p,y))
&= \textstyle{\sum_{(x,r,w) \in G_{\Lambda}} \psi (f) ((x,r,w)) \psi (g) ((w,p-r,y))}
\\&=
\begin{cases}
f(x(0,m)) g(x(m,m+n)) &\text{if $p=m+n$ and $\sigma^{m+n} (x) = y$}\\
0 &\text{otherwise}
\end{cases}
\\&= \psi (fg) ((x,p,y)).
\end{align*}

We now show that $\psi (\langle f,g \rangle_{A}^{n}) = \psi
(f)^{*} \ast \psi (g)$ for $f,g \in X_{n}$, and for this it is
enough to show that $\psi (\langle f,g \rangle_{A}^{n}) = \psi
(f)^{*} \ast \psi (g)$ for $f,g \in C_{c} (\Lambda^{n})$.
Noting that $\langle f,g \rangle_{A}^{n} \in C_{c}
(\Lambda^{0})$ by \cite[Lemma 1.5]{ka1}, for $(x,p,y) \in
G_{\Lambda}$ we have
\begin{align*}
\psi (f)^{*}& \ast \psi (g) ((x,p,y))
\\&=
\begin{cases}
\textstyle{\sum_{\alpha \in \Lambda^{n} r(x)}} \overline{\psi (f) ((\alpha x,n,x))} \psi (g) ((\alpha x,n,x)) &\text{if $p=0$ and $x=y$}\\
0 &\text{otherwise}
\end{cases}
\\&=
\begin{cases}
\textstyle{\sum_{\alpha \in \Lambda^{n} r(x)}} \overline{f(\alpha)} g(\alpha) &\text{if $p=0$ and $x=y$}\\
0 &\text{otherwise}
\end{cases}
\\&= \psi (\langle f,g \rangle_{A}^{n}) ((x,p,y)),
\end{align*}
and the result follows.

Uniqueness of the Toeplitz representation $\psi$ follows from the facts that each of the linear maps $\psi |_{X_{n}}$
is continuous and $C_{c} (\Lambda^{n})$ is dense in $X_{n}$.

Finally, to see that $\{\, \psi (f) \mid n \in \NN^{k}, f \in C_{c} (X_{n}) \,\}$
generates $C^{*} (G_{\Lambda})$, note that the norm on $C^{*} (G_{\Lambda})$ is
dominated by the norm $\| \cdot \|_{I}$ from \cite[Page 50]{r} which in turn is dominated by $\| \cdot \|_{\infty}$. Our strategy is to use the Stone-Weierstrass Theorem to show that the subalgebra $\Aa$ generated by $\{\, \psi (f) \mid n \in \NN^{k}, f \in C_{c} (X_{n}) \,\}$ is dense in $C_{0} (G_{\Lambda})$. Since $\Aa$ is by definition a subset of $C_c(G_\Lambda)$, it will then follow that $\Aa$ is dense in $C_c(G_\Lambda)$, and hence in $C^*(G_\Lambda)$. So it is enough to show that, for distinct $(x,p,y), (x',p',y') \in G_{\Lambda}$, there
exist $m,n \in \NN^{k}$, $f \in C_{c} (\Lambda^{m})$ and $g \in C_{c} (\Lambda^{n})$ such that
$\psi (f) \ast \psi (g)^{*} ((x,p,y)) \ne \psi (f) \ast \psi (g)^{*} ((x',p',y'))$.
So denote $(x,p,y) = (\lambda z,d(\lambda)-d(\mu),\mu z)$. If $p \ne p'$ then choose
$f \in C_{c} (\Lambda^{d(\lambda)})$ with $f(\lambda)=1$ and $g \in C_{c} (\Lambda^{d(\mu)})$
with $g(\mu)=1$. Then
\[
\psi (f) \ast \psi (g)^{*} ((x,p,y))
= \psi (f) ((x,d(\lambda),z)) \overline{\psi (g) ((y,d(\mu),z))}
= f(\lambda) \overline{g(\mu)}
=1
\]
and $\psi (f) \ast \psi (g)^{*} ((x',p',y')) = 0$.

Suppose $p=p'$ and $x \ne x'$. Assume without loss of generality that $d(x) \not< d(x')$.
If $d(x)=d(x')$ then there exists $m \in \NN^{k}$ such that $d(\lambda) \le m \le d(x)$
and $x(0,m) \ne x'(0,m)$. Choose $f \in C_{c} (\Lambda^{m})$ such that $f(x(0,m)) =1$ and
$f(x'(0,m)) = 0$, and choose $g \in C_{c} (\Lambda^{m-p})$ such that $g(y(0,m-p))=1$. Then,
since $\sigma^{m-p} (y) = \sigma^{m} (x)$,
\begin{align}\label{sep1}
\psi (f) \ast \psi (g)^{*} ((x,p,y))
&= \psi (f) ((x,m,\sigma^{m} (x))) \overline{\psi (g) ((y,m-p,\sigma^{m-p} (y)))}\notag
\\&= f(x(0,m)) \overline{g(y(0,m-p))}
=1\quad \text{and}
\\\psi (f) \ast \psi (g)^{*} ((x',p',y')) &= 0.\notag
\end{align}
If $d(x) \ne d(x')$ then there exists $m \in \NN^{k}$ such that $d(\lambda) \le m \le d(x)$
and $m \not\le d(x')$. So, choosing $f \in C_{c} (\Lambda^{m})$ such that $f(x(0,m)) =1$ and
$g \in C_{c} (\Lambda^{m-p})$ such that $g(y(0,m-p))=1$ gives (\ref{sep1}).

The case where $p=p'$, $x = x'$ and $y \ne y'$ follows from an argument similar to that of
the preceding paragraph.
\end{proof}

To show that $\psi$ is Nica covariant we will use the following
lemma.

\begin{lemma}\label{rank1}
If $m \in \NN^{k}$, $f,g \in C_{c} (\Lambda^{m})$, and $(x,p,y)
\in G_{\Lambda}$ then
\begin{align*}
\psi^{(m)} (f \otimes g^{*}) &((x,p,y))
\\&=
\begin{cases}
f(x(0,m)) \overline{g(y(0,m))} &\text{if $p=0$, $m \le d(x)$ and $\sigma^{m} (x) = \sigma^{m} (y)$}\\
0 &\text{otherwise.}
\end{cases}
\end{align*}
\end{lemma}
\begin{proof}
Recalling that the Haar system on $G_{\Lambda}$ consists of
counting measures, we calculate
\begin{equation}\label{pims}
\psi^{(m)} (f \otimes g^{*}) ((x,p,y))
= \textstyle{\sum_{(x,r,z) \in G_{\Lambda}}} \psi (f) ((x,r,z)) \psi (g)^{*} ((z,p-r,y)).
\end{equation}
The right-hand side of~\eqref{pims} is zero unless $p=0$, $m
\le d(x)$ and $\sigma^{m} (x) = \sigma^{m} (y)$ (noting that
$p=0$ implies $d(x)=d(y)$), and if indeed $p=0$, $m \le d(x)$
and $\sigma^{m} (x) = \sigma^{m} (y)$, then~\eqref{pims}
simplifies to
\[
\psi (f) ((x,m,\sigma^{m} (x))) \psi (g)^{*} ((\sigma^{m} (y),-m,y))
= f(x(0,m)) \overline{g(y(0,m))}.\qedhere
\]
\end{proof}

\begin{prop}\label{prop:THRG NicaCov}
Let $\Lambda$ be a compactly aligned topological $k$-graph.
The Toeplitz representation $\psi \colon X \to C^{*} (
G_{\Lambda} )$ from Proposition \ref{prop:thrg-rep}
is gauge-compatible and Nica covariant.
\end{prop}
\begin{proof}
The canonical continuous cocycle $c:(x,m,y) \mapsto m$ on
$G_\Lambda$ induces a coaction $\beta$ of $\ZZ^k$ on
$C^*(G_\Lambda)$ satisfying $\beta(f) = f \otimes
i_{\ZZ^k}(m)$ whenever $\supp(f) \subset c^{-1}(m)$. In
particular, $\beta(\psi(f)) = \psi(f) \otimes i_{\ZZ^k}(m)$ for
$f \in X_m$, so $\psi$ is gauge-compatible.

For Nica covariance, fix $m,n \in \NN^{k}$, $f_{m},g_{m} \in F_{m}$ and $f_{n},g_{n}
\in F_{n}$. Resume the notation of Lemma~\ref{compact} so that
\[
\iota_{m}^{m \vee n} (g_{m} \otimes f_{m}^{*}) \iota_{n}^{m \vee n} (f_{n} \otimes g_{n}^{*})
= \textstyle{\sum_{i=1}^{r_{m}} \sum_{j=1}^{r_{n}}} a_{ij} \otimes b_{j}^{*}.
\]
By Lemma~\ref{dense} it suffices to
show that for every $(x,p,y) \in G_\Lambda$,
\[
\psi^{(m)} (g_{m} \otimes f_{m}^{*}) \psi^{(n)} (f_{n} \otimes g_{n}^{*}) ((x,p,y))
    = \psi^{(m \vee n)} (\iota_{m}^{m \vee n} (g_{m} \otimes f_{m}^{*}) \iota_{n}^{m \vee n} (f_{n} \otimes g_{n}^{*})) ((x,p,y)).
\]
By Lemma~\ref{rank1}, both sides of this equation are equal to
zero unless $p = 0$. By definition of the multiplication in
$C^*(G_\Lambda)$ and another application of Lemma~\ref{rank1},
we are left to show that for all $(x,0,y) \in G_\Lambda$,
\begin{equation}\label{NicaCov2}
\begin{split}
\textstyle{\sum_{(x,0,z) \in G_{\Lambda}}} \psi^{(m)} (g_{m} \otimes f_{m}^{*}) &((x,0,z)) \psi^{(n)} (f_{n} \otimes g_{n}^{*}) ((z,0,y)) \\
 &= \psi^{(m \vee n)} (\iota_{m}^{m \vee n} (g_{m} \otimes f_{m}^{*}) \iota_{n}^{m \vee n} (f_{n} \otimes g_{n}^{*})) ((x,0,y)).
\end{split}\end{equation}
We may further deduce from Lemma~\ref{rank1} that both sides
are equal to zero unless $m \vee n \le d(x)$ and
$\sigma^{m \vee n}(x) = \sigma^{m \vee n}(y)$ (noting that $p=0$ implies $d(x) = d(y)$).

So fix $(x,0,y) \in G_\Lambda$ such that $m\vee n \le d(x)$ and
$\sigma^{m \vee n}(x) = \sigma^{m \vee n}(y)$. Recall from
Notation~\ref{ntn:Fms} that for each $h \in F_r$ and $v \in
\Lambda^0$ such that $\Lambda^{r} v$ is non-empty we have a fixed
path $\lambda_{h,v} \in \Lambda^{r} v$
such that $h(\lambda) = 0$ for all $\lambda \in \Lambda^{r} v
\setminus \{\lambda_{h,v}\}$. Let $\lambda_m := \lambda_{f_m,
x(m)}$. We calculate, using Lemma~\ref{rank1} yet again, that
\begin{align}
\textstyle\sum_{(x,0,z) \in G_\Lambda} &
            \psi^{(m)} (g_{m} \otimes f_{m}^{*}) ((x,0,z)) \psi^{(n)} (f_{n} \otimes g_{n}^{*}) ((z,0,y)) \nonumber\\
            &= \textstyle{\sum_{\{\,\zeta \in \Lambda^m x(m) \mid \sigma^{n} (\zeta \sigma^{m} (x)) = \sigma^{n} (y)\,\}}} g_{m}(x(0,m))
                \overline{f_{m}(\zeta)} f_{n}([\zeta \sigma^m(x)](0,n)) \overline{g_{n}(y(0,n))} \nonumber\\
            &= \begin{cases}\label{calc6}
            		g_{m}(x(0,m)) \overline{f_{m}(\lambda_m)} f_{n}([\lambda_m \sigma^m(x)](0,n)) \overline{g_{n}(y(0,n))}\\
               \hspace{4cm}\text{if $y(n,m \vee n) = [\lambda_m \sigma^{m} (x)](n,m \vee n)$}\\
                0\hspace{3.8cm}\text{otherwise.}
	             \end{cases}
\end{align}

We will show that the right-hand side of \eqref{NicaCov2} is
equal to~\eqref{calc6}. For $1 \le i \le r_{m}$ and $1 \le j \le r_{n}$, we have
\begin{align*}%\label{calc1}
\psi^{(m \vee n)}&(a_{ij} \otimes b_{j}^{*}) ((x,0,y)) \notag\\
    &= a_{ij} (x(0,m \vee n)) \overline{b_{j} (y(0,m \vee n))}\notag \\
    &= g_{m}(x(0,m)) [\rho_{i}^{m} \cdot \langle f_{m} \rho_{i}^{m},f_{n} \rho_{j}^{n} \rangle_{A}^{m \vee n}]
        (x(m,m \vee n)) \overline{b_{j} (y(0,m \vee n))} \\
    &= \begin{cases}
    g_m(x(0,m)) [\rho_{i}^{m} (x(m,m \vee n))]^{2} \overline{f_{m}(\lambda_m)} f_{n}([\lambda_m \sigma^{m} (x)](0,n))\\
     \hspace{5cm}\overline{g_{n}(y(0,n))} [\rho_{j}^{n} ([\lambda_m \sigma^{m} (x)](n,m \vee n))]^{2}\\
     \hspace{6cm}\text{if $y(n,m \vee n) = [\lambda_m \sigma^{m} (x)](n,m \vee n)$}\\
     0\hspace{5.8cm}\text{otherwise.}
     \end{cases}
\end{align*}
Hence, if $y(n, m \vee n) \not= [\lambda_m \sigma^{m} (x)](n,m \vee n)$
then each of the right-hand side of \eqref{NicaCov2}~and~\eqref{calc6} is
equal to zero; and if $y(n, m \vee n) = [\lambda_m \sigma^{m}
(x)](n,m \vee n)$ we calculate:
\begin{align*}
\psi^{(m \vee n)}&(\iota_{m}^{m \vee n}(g_{m} \otimes f_{m}^{*}) \iota_{n}^{m \vee n} (f_{n} \otimes g_{n}^{*})) ((x,0,y))\\
    &= \textstyle \sum_{i=1}^{r_{m}} \sum_{j=1}^{r_{n}}  \psi^{(m \vee n)} (a_{ij} \otimes b_{j}^{*}) ((x,0,y)) \\
    &= \textstyle \sum_{i=1}^{r_{m}} \sum_{j=1}^{r_{n}}  g_{m}(x(0,m)) [\rho_{i}^{m} (x(m,m \vee n))]^{2} \overline{f_{m}(\lambda_m)}\\
     &\phantom{= \textstyle{\sum_{i=1}^{r_{m}} \sum_{j=1}^{r_{n}}}\quad} f_{n}([\lambda_{m} \sigma^{m} (x)](0,n)) \overline{g_{n}(y(0,n))} [\rho_{j}^{n} ([\lambda_{m} \sigma^{m} (x)](n,m \vee n))]^{2}\\
    &= g_{m}(x(0,m)) \overline{f_{m}(\lambda_{m})} f_{n}([\lambda_{m} \sigma^{m} (x)](0,n)) \overline{g_{n}(y(0,n))}
\end{align*}
as required.
\end{proof}

\begin{theorem} \label{thm:thrg}
Let $\Lambda$ be a compactly aligned topological $k$-graph.
Let $G_\Lambda$ and $\Gg_\Lambda$ be Yeend's path groupoid and
boundary-path groupoid for $\Lambda$, let $C^*(G_\Lambda)$
and $C^*(\Gg_\Lambda)$ and $C^*_\reduced(G_\Lambda)$
and $C^*_\reduced(\Gg_\Lambda)$ be the associated full and reduced $C^*$-algebras,
and let $q \colon C^*(G_\Lambda) \to C^*(\Gg_\Lambda)$ be the quotient map.
Let $X$ be the product system defined by Proposition \ref{prp:prod sys}
and let $\psi_* \colon \Tc(X) \to C^*(G_\Lambda)$ be the homomorphism
obtained from the universal property of $\Tc(X)$ and
Proposition~\ref{prop:THRG NicaCov}.

Then $\psi_*$ is an isomorphism, the canonical maps
$C^*(G_\Lambda)\to C^*_\reduced(G_\Lambda)$ and
$C^*(\Gg_\Lambda)\to C^*_\reduced(\Gg_\Lambda)$ are isomorphisms,
and there is a unique $*$-isomorphism
$\phi \colon C^*(\Gg_\Lambda) \to \NO{X}$
which makes the following diagram commute.
\[
\begin{CD}
\Tc(X) @>{\psi_*}>> C^*(G_\Lambda) \\
@V{q_{\CNP}}VV @V{q}VV \\
\NO{X} @<{\phi}<< C^*(\Gg_\Lambda).
\end{CD}
\]
\end{theorem}
\begin{proof}
% It follows from Proposition~\ref{prop:THRG NicaCov}
% that $\psi$ is a gauge-compatible Nica covariant
% representation of $X$ which generates $C^*(G_\Lambda)$.
% Thus it follows from the universal property of $\Tc(X)$ that
% there exists a $*$-homomorphism $\psi_*:\Tc(X)\to C^*(G_\Lambda)$
% such that $\psi_*\circ i_X=\psi$.
Since $\psi(X)$ generates
$C^*(G_\Lambda)$, it follows that $\psi_*$ is surjective.
Let $\pi_{G_\Lambda}$ denote the canonical factor
map from $C^*(G_\Lambda)$ to
$C^*_\reduced(G_\Lambda)$ and let
$\tilde{\psi}_*:=\pi_{G_\Lambda}\circ\psi_*$. We aim to show that
$\tilde{\psi}_*$ it injective. It will then follow that both $\psi_*$
and $\pi_{G_\Lambda}$ are isomorphisms.
The canonical continuous cocycle $c:(x,m,y) \mapsto m$ from
$G_\Lambda$ to $\ZZ^k$ induces a strongly continuous action of $\TT^k$ on
$C^*_\reduced(G_\Lambda)$ (cf. \cite[Proposition II.5.1]{r}) and thereby a coaction
$\beta$ of $\ZZ^k$ on
$C^*_\reduced(G_\Lambda)$ satisfying $\beta(f) = f \otimes
i_{\ZZ^k}(m)$ whenever $\supp(f) \subset c^{-1}(m)$. Thus
$\tilde{\psi}_*$ is equivariant for $\delta$ and $\beta$.
Since $\ZZ^k$ is amenable, the coaction $\delta$ is normal, so
it is enough to prove that
$\ker\tilde{\psi}_*\cap\Ff=\{0\}$ which we will do now. It follows from equation
\eqref{eq:BPs span F}, Lemma \ref{lem:core structure} and \cite[Lemma~1.3]{ALNR}
that it is enough to prove that
$\ker\tilde{\psi}_*\cap B_F=\{0\}$ for each $F \in \fvcl{P}$. For
this, we fix $F \in \fvcl{P}$ and generalised compact operators
$T_p \in \Kk(X_p)$ for $p \in F$ such that
\begin{equation}\label{eq:psistar_inj}
\sum_{p \in F}
\tilde{\psi}_*(i_X^{(p)}(T_p))=0.
\end{equation}
We use induction over the number of elements in $F$ to show
that $T_p=0$ for each $p\in F$. Proposition~\ref{prop:thrg-rep}
implies that the representation $\pi_{G_\Lambda}\circ\psi$ is
injective. Thus it follows from Lemma 2.4 of \cite{ka1} that
$\tilde{\psi}_*\circ
i_X^{(p)}=(\pi_{G_\Lambda}\circ\psi)^{(p)}$ is injective for
all $p\in\NN^k$. So if $F=\{p\}$, then $T_p=0$. Assume then
that $F$ consists of more than one element. Let $p_0$ be a
minimal element of $F$. Then $p\not\le p_0$ for all $p\in
F\setminus\{p_0\}$. It follows from Lemma \ref{rank1} that if
$(x,0,y)\in G_\Lambda$ with $d(x)=d(y)=p_0$, then
$\tilde{\psi}_* (i_X^{(p)}(T_p))((x,0,y))=0$ for all $p\in
F\setminus\{p_0\}$. The assumption \eqref{eq:psistar_inj} then
implies that $\tilde{\psi}_* (i_X^{(p_0)}(T_{p_0}))((x,0,y))=0$
for all $(x,0,y)\in G_\Lambda$ with $d(x)=d(y)=p_0$. It then
follows from Lemma \ref{rank1} that $\tilde{\psi}_*
(i_X^{(p_0)}(T_{p_0}))((x,p,y))=0$ for all $(x,p,y)\in
G_\Lambda$, and thus that $\tilde{\psi}_*
(i_X^{(p_0)}(T_{p_0}))=0$. As before, this implies that
$T_{p_0}=0$, and it then follows from our inductive hypothesis
that $T_p=0$ for every $p\in F$. Thus $\tilde{\psi}_*$ is
injective.

% Let $\pi_{G_\Lambda}$ denote the canonical map from $C^*(G_\Lambda)$ to
% $C^*_\reduced(G_\Lambda)$. Let $\Phi_{C^*_\reduced(G_\Lambda)}$ be the
% unique conditional expectation from $C^*_\reduced(G_\Lambda)$ onto
% $C^*(G_\Lambda^{(0)})$ (cf. \cite[Proposition II.4.8]{r}). Then $(\pi_{G_\Lambda}\circ\psi_*)\circ\Phi^{\delta}=
% \Phi_{C^*_\reduced(G_\Lambda)}\circ (\pi_{G_\Lambda}\circ\psi_*)$.
% Since $\ZZ^k$ is amenable it follows that $\Phi^{\delta}$
% is faithful. The standard intertwining argument then shows that
% $\pi_{G_\Lambda}\circ\psi_*$, and thus $\pi_{G_\Lambda}$, are injective.

Let $\pi_{\Gg_\Lambda}$ denote the canonical map from $C^*(\Gg_\Lambda)$ to
$C^*_\reduced(\Gg_\Lambda)$
% Since $\Gg_\Lambda=G_{\Lambda|\partial\Lambda}$
% it follows that $\ker q$ is the closure of
% $$\{f\in C_C(G_\Lambda)\mid f((x,m,y))=0\text{ if }(x,m,y)\in
% G_\Lambda\text{ and }x,y\in \partial\Lambda\}.$$
% It follows that there is an approximate identity for
% $\ker q$ in $C_0(G_\Lambda^{(0)}\setminus\partial\Lambda)$, and
% $\ker q$ is therefore generated by its intersection with $C_0(G_\Lambda^{(0)})$.
% Hence $\ker (q\circ\psi_*)$ is generated by its intersection with $\Tc(X)^\delta_e$
% and since $\ZZ^k$ amenable and thus exact, it follows from Proposition
% \ref{prop:exact} that $\delta^{\ker(q\circ\psi_*)}$ is an exa
% and thus with its intersection with $C^*(G_\Lambda[c])$ where
% $G_\Lambda[c]$ is the subgroupoid $c^{-1}(\{0\})= \{(x,m,y)\in G_\Lambda\mid m=0\}$.
and let $\rho : X \to C^*(\Gg_\Lambda)$ be the Toeplitz
representation $\pi_{\Gg_\Lambda}\circ q \circ \psi$. Then
Proposition~\ref{prop:THRG NicaCov} implies that $\rho$ is a
Nica covariant representation of $X$ which generates
$C^*_\reduced(\Gg_\Lambda)$. Proposition~4.3 of~\cite{y}
implies that $\rho_e : A \to C^*_\reduced(\Gg_\Lambda)$ is
injective, so $\rho$ is injective as a representation of $X$.
As above the canonical coaction $(x,m,y)\mapsto m$ from
$\Gg_\Lambda$ to $\NN^k$ induces a coaction $\gamma$ of $\ZZ^k$
on $C^*_\reduced(\Gg_\Lambda)$ such that
$\gamma(\rho(x))=\rho(x)\otimes i_{\ZZ^k}(d(x))$ for $x\in X$.
Thus $\rho$ is gauge-compatible. As noted in the proof of
Corollary~\ref{cor:Nk giut}, a product system over $\NN^k$ is
automatically $\tilde\phi$-injective. Thus it follows from
Theorem~\ref{thm:projective property} that there exists a
surjective $*$-homomorphism $\phi':C^*_\reduced(\Gg_\Lambda)\to
\NO{X}$ such that $\phi'\circ \pi_{\Gg_\Lambda} \circ q \circ
\psi_*=q_{\CNP}$. Let $\phi:=\phi'\circ \pi_{\Gg_\Lambda}$. We
will show that $\phi$ is injective. It will then follow that
$\pi_{\Gg_\Lambda}$ is an isomorphism from $C^*(\Gg_\Lambda)$
to $C^*_\reduced(\Gg_\Lambda)$, and that $\phi$ is an
isomorphism from $C^*(\Gg_\Lambda)$ to $\NO{X}$ such that
$\phi\circ q\circ\psi_*=q_{\CNP}$.  The various maps defined so
far are summarised in the following commuting diagram. We have
established already that all three maps in the top row are
isomorphisms:
\[
\begin{tikzpicture}
    \node[inner sep=1pt] (Tc) at (0,3) {\small$\Tc(X)$};%
    \node[inner sep=1pt] (C*G) at (3,3) {\small$C^*(G_\Lambda)$};%
    \node[inner sep=1pt] (C*rG) at (6,3) {\small$C^*_r(G_\Lambda)$};%
    \node[inner sep=1pt] (NO) at (0,0) {\small$\NO{X}$};%
    \node[inner sep=1pt] (C*Gg) at (3,0) {\small$C^*(\Gg_\Lambda)$};%
    \node[inner sep=1pt] (C*rGg) at (6,0) {\small$C^*_r(\Gg_\Lambda)$};%
    \draw[->] (Tc)--(C*G) node[pos=0.5,anchor=north,inner sep=1pt] {\small $\psi_*$};%
    \draw[->] (C*G)--(C*rG) node[pos=0.5,anchor=north,inner sep=1pt] {\small $\pi_{G_\Lambda}$};%
    \draw[->] (C*Gg)--(NO) node[pos=0.5,anchor=south,inner sep=1pt] {\small $\phi$};%
    \draw[->] (C*Gg)--(C*rGg) node[pos=0.5,anchor=south,inner sep=1pt] {\small $\pi_{\Gg_\Lambda}$};%
    \draw[->] (Tc)--(NO) node[pos=0.5,anchor=east,inner sep=1pt] {\small $q_{{\CNP}}$};%
    \draw[->] (C*G)--(C*Gg) node[pos=0.5,anchor=east,inner sep=1pt] {\small $q$};%
    \draw[<-] (C*rG) .. controls (3,3.75) .. (Tc) node[pos=0.5,anchor=south,inner sep=1pt] {\small $\tilde\psi_*$};%
    \draw[->] (C*rGg) .. controls (3,-0.75) .. (NO) node[pos=0.5,anchor=north,inner sep=1pt] {\small $\phi'$};%
\end{tikzpicture}
\]

To show that $\phi$ is injective it suffices to prove that
$\ker(q_{\CNP}\circ \tilde{\psi}_*^{-1})\subset
\ker (q \circ \pi_{G_\Lambda}^{-1})$.
Since $\Gg_\Lambda=G_{\Lambda|\partial\Lambda}$ it follows that
$\ker (q \circ \pi_{G_\Lambda}^{-1})$ is the closure of
$$
\{f\in C_c(G_\Lambda)\mid f((x,m,y))=0\text{ if }(x,m,y)\in
G_\Lambda\text{ and }x,y\in \partial\Lambda\}.
$$
It follows that there is an approximate identity for
$\ker (q \circ \pi_{G_\Lambda}^{-1})$ in
$C_0(G_\Lambda^{(0)}\setminus\partial\Lambda)$, and
$\ker (q \circ \pi_{G_\Lambda}^{-1})$ is therefore generated
by its intersection with $C_0(G_\Lambda^{(0)})$, and thus by its
intersection with $C^*_\reduced(G_\Lambda[c])$ where
$G_\Lambda[c]$ is the subgroupoid $c^{-1}(\{0\})= \{(x,m,y)\in G_\Lambda\mid m=0\}$.
Thus it follows from Proposition \ref{prp:deltaI existence} that $\beta$
induces a coaction $\beta^{\ker (q\circ\pi_{G_\Lambda}^{-1})}$ of $\ZZ^k$
on $C^*(\Gg_\Lambda)$, which is normal since $\ZZ^k$ is amenable.
Since $q_{\CNP}\circ\tilde{\psi}_*^{-1}$ is equivariant for $\beta$ and $\CNPgaug$,
it follows that it suffices to show that
$\ker(q_{\CNP}\circ\tilde{\psi}_*^{-1})\cap C^*_\reduced(G_\Lambda[c])
\subset \ker (q \circ \pi_{G_\Lambda}^{-1})$.

By \cite[Theorem 3.16]{y}, $G_\Lambda$, and thus $G_\Lambda[c]$, are
$r$-discrete, and since $G_\Lambda[c]$ is
also (essentially) principal, \cite[Proposition II.4.6]{r} implies that
$$
Y:=\{x\in G_\Lambda^{(0)}\mid f((x,0,x))=0\text{ for all }f\in \ker(q_{\CNP}
\circ\tilde{\psi}_*^{-1})\cap C^*_\reduced(G_\Lambda[c])\}$$
is a closed $G_\Lambda[c]$-invariant subset of $G_\Lambda^{(0)}$ such that
$\ker(q_{\CNP}\circ\tilde{\psi}_*^{-1})\cap C^*_\reduced(G_\Lambda[c])$
is the closure of
$$
\{f\in C_c(G_\Lambda[c])\mid f((x,0,y))=0\text{ if }(x,0,y) \in G_\Lambda[c]\text{ and }x,y\in Y\}.
$$
We claim that $Y$ is not just $G_\Lambda[c]$-invariant, but also
$G_\Lambda$-invariant; indeed if $(x,m,y)\in G_\Lambda$ and $x\notin Y$, then
there exists
$f\in \ker(q_{\CNP}\circ\tilde{\psi}_*^{-1})\cap C^*_\reduced(G_\Lambda[c])$
such that $f((x,0,x))\ne 0$. Since $G_\Lambda$ is $r$-discrete there is
$g\in C_c(G_\Lambda)$ such that $g$ is supported on a subset on which the source map
is bijective, and such that $g((x,m,y))=1$. We then have that
\begin{align*}
    &(g^*\ast f\ast g)((y,0,y))\\
    &\qquad =\textstyle{\sum_{(y,m_1,z_1),(z_2,m_2,y)\in
G_\Lambda}}\overline{g((z_1,-m_1,y))}f((z_1,-m_1-m_2,z_2))g((z_2,m_2,y))\\
    &\qquad =\overline{g((x,m,y))}f((x,0,x))g((x,m,y))=f((x,0,x))\ne 0.
\end{align*}
Let  $\Phi^\beta$ be the conditional expectation of $C^*_\reduced(G_\Lambda)$ onto
$C^*_\reduced(G_\Lambda[c])$ induced by the coaction $\beta$. Then $\Phi^{\beta}(g^*\ast f\ast g)((y,0,y))=(g^*\ast f\ast g)((y,0,y))\ne 0$,
and since $\ker(q_{\CNP}\circ\tilde{\psi}_*^{-1})$
is generated by its intersection with the subset
$\tilde{\psi}_*(\mathcal{F})$ of $C^*_\reduced(G_\Lambda[c])$ it follows from \cite[Theorem 3.9]{Exel} that
$\Phi^{\beta}(g^*\ast f\ast g)\in
\ker(q_{\CNP}\circ\tilde{\psi}_*^{-1})\cap
C^*_\reduced(G_\Lambda[c])$. Thus $y\notin Y$, showing that $Y$ is $G_\Lambda$-invariant.
Since $q_{\CNP}\circ\tilde{\psi}_*^{-1}$ is injective on $\pi_{G_\Lambda}(\psi(C_0(\Lambda^{(0)})))$ we must have
that $vY\ne\emptyset$ for all $v\in\Lambda^0$. Thus $Y$ is a closed and invariant
subset of $G_\Lambda^{(0)}$ which satisfies that $vY\ne\emptyset$ for all
$v\in\Lambda^0$.
It therefore follows from Proposition \ref{prop:small} that $\partial\Lambda\subset Y$.
Thus $\ker(q_{\CNP}\circ\tilde{\psi}_*^{-1})\cap C^*_\reduced(G_\Lambda[c])
\subset \ker (q \circ \pi_{G_\Lambda}^{-1})$, as claimed.
\end{proof}

By combining Theorem \ref{thm:thrg} with Corollary \ref{cor:Nk giut} we get the
following gauge-invariant uniqueness result for $C^*(\Gg_\Lambda)$.

\begin{cor} \label{cor:thrg-giut}
  Let $\Lambda$ be a compactly aligned topological $k$-graph.
  Let $\Gg_\Lambda$ be Yeend's
  boundary-path groupoid for $\Lambda$, and let $C^*(\Gg_\Lambda)$
  be the associated full $C^*$-algebra.
  Let $\psi \colon X \to C^{*} (G_{\Lambda} )$ be the map from Proposition
  \ref{prop:thrg-rep}.

  A surjective $*$-homomorphism $\phi:C^*(\Gg_\Lambda)\to B$ is injective if and only if
  \begin{enumerate}\renewcommand{\theenumi}{\arabic{enumi}}
  \item there is a strongly continuous
    action $\alpha$ of $\TT^k$ on $B$ such that
    $\alpha_z(\phi(\psi(x)) = z^{d(x)}\phi(\psi(x))$ for all $x \in
    X$ and $z\in\TT^k$, and
\item $\phi\vert_{\psi(C_0(\Lambda^0))} \colon A \to B$ is injective.
\end{enumerate}
\end{cor}

We note that Yamashita in \cite{Yam09} has studied $\NO{X}$
under the assumption that $\Lambda$ is row-finite and without
sources. Among other things he shows a Cuntz-Krieger type
uniqueness theorem for $\NO{X}$ and gives sufficient conditions
for when $\NO{X}$ is simple and purely infinite.

% \begin{rmk}
% It follows that $\NO{X}$ is the unique quotient of Yeend's
% Toeplitz algebra which satisfies a gauge-invariant uniqueness
% theorem, and should therefore be regarded as the Cuntz-Krieger
% algebra of the topological higher-rank graph. That is, the
% co-universal property of $\NO{X}$ means that it is
% automatically the ``right" Cuntz-Krieger algebra (this is not
% established for $C^*(G_\Lambda)$ in~\cite{y}).
% \end{rmk}

\appendix

\section{Coactions, quotients and normality}

\begin{prop}\label{prp:deltaI existence}
Let $A$ be a $C^*$-algebra carrying a coaction $\delta$ of a
discrete group $G$. Let $I$ be an ideal of $A$ which is
generated as an ideal by $I^\delta_e:=I \cap A^\delta_e$. Let
$q_I \colon A \to A/I$ be the quotient map. Then there is a coaction
$\delta^I$ of $G$ on $A/I$ such that
\begin{equation}\label{eq:delta^I formula}
\delta^I \circ q_I = (q_I \otimes \id_{C^*(G)}) \circ \delta.
\end{equation}
\end{prop}

\begin{proof}
The proposition is trivially true if $A=\{0\}$, so we may, and
will, assume that $A$ contains a non-zero element.

To define the homomorphism $\delta^I$, observe that for $a \in
A^\delta_e$, we have $(q_I \otimes \id_{C^*(G)}) \circ \delta(a) =
q_I(a) \otimes i_G(e)$, which is equal to zero if and only
if $a \in I$; that is $\ker((q_I \otimes \id_{C^*(G)}) \circ
\delta) \cap A^\delta_e = I^\delta_e$. Since $I$ is generated as an
ideal by $I^\delta_e$, it follows that $I \subset \ker(q_I
\otimes \id_{C^*(G)}) \circ \delta)$, and hence $(q_I \otimes
\id_{C^*(G)}) \circ \delta$ descends to a homomorphism $\delta^I
\colon A/I \to (A/I) \otimes C^*(G)$ satisfying~\eqref{eq:delta^I
formula}.

We will show that $\delta^I$ is a coaction. Since $G$ is
discrete, it suffices to show that $\delta^I$ is nondegenerate
and injective and satisfies the coaction identity (see
\cite[Section~1]{qui:discrete coactions}). It follows from
\cite[Corollary 1.6]{qui:discrete coactions} that $A^\delta_e$
contains an approximate identity $(u_\lambda)_{\lambda \in
\Lambda}$ for $A$. To see that $\delta^I$ is nondegenerate,
note that $(q_I(u_\lambda))_{\lambda \in \Lambda}$ is an
approximate identity for $A/I$. Since each $u_\lambda \in
A^\delta_e$, we have $\delta^I(q_I(u_\lambda)) = q_I(u_\lambda)
\otimes \id_{C^*(G)}$ for all $\lambda$, and hence
$(\delta^I(u_\lambda))_{\lambda \in \Lambda}$ is an approximate
identity for $(A/I) \otimes C^*(G)$. So $\delta^I$ is
nondegenerate.

To see that $\delta$ is injective, let $\epsilon \colon C^*(G)
\to \CC$ be the augmentation representation $i_G(g) \mapsto 1$
for all $g \in G$. For $g \in G$ and $a \in A^\delta_g$, we have
\[
(\id_{A/I} \otimes \epsilon) \circ \delta^I(q_I(a))
    = (\id_{A/I} \otimes \epsilon)(q_I(a) \otimes i_G(g))
    = q_I(a) \otimes 1.
\]
Lemma~1.5 of \cite{qui:discrete coactions} shows that $A =
\clsp(\bigcup_{g \in G} A^\delta_g)$, so the preceding calculation
together with linearity and continuity of the homomorphism
$(\id_{A/I} \otimes \epsilon) \circ \delta$ show that
\[
(\id_{A/I} \otimes \epsilon) \circ \delta^I (x) = x \otimes 1
\]
for all $x \in A/I$. Hence $(\id_{A/I} \otimes \epsilon) \circ
\delta^I$ is injective, and in particular, $\delta^I$ is injective.

To see that $\delta^I$ satisfies the coaction identity, let
$\delta_G$ be the comultiplication on $C^*(G)$, fix $g \in G$
and $a \in A^\delta_g$, and calculate
\[
(\id_{A/I} \otimes \delta_G)\circ\delta^I(q_I(a))
    = q_I(a) \otimes i_G(g) \otimes i_G(g)
    = (\delta^I \otimes \id_{C^*(G)})\circ \delta^I(q_I(a)).
\]
It then follows from
linearity and continuity that $(\id_{A/I} \otimes
\delta_G)\circ\delta^I= (\delta^I \otimes \id_{C^*(G)})\circ
\delta^I$. We have now established that
$\delta^I$ is a coaction.
\end{proof}

\begin{rmk}
One could also prove Proposition \ref{prp:deltaI existence} by using
the duality between coactions of $G$ and topological $G$-gradings
(cf. \cite{qui:discrete coactions}) and \cite[Proposition 3.11]{Exel}.
\end{rmk}

\begin{notation}
Given a discrete group $G$, we will write $\lambda_G$ for the
left regular representation of $G$ on $\ell^2(G)$, and also for
the resulting homomorphism of $C^*(G)$ onto $C^*_\reduced(G)$ obtained
from the universal property of $C^*(G)$ applied to $\lambda_G$.
\end{notation}

\begin{lemma}\label{lem:when deltaI normal}
Resume the hypotheses of Proposition~\ref{prp:deltaI
existence}. The following are equivalent.
\begin{enumerate}
\item\label{it:deltaI normal} $\delta^I$ is normal.
\item\label{it:PhiI faithful} $\Phi^{\delta^I}$ is faithful
    on positive elements.
\item\label{it:deltar injective} $(\id_{A/I} \otimes \lambda_G)
    \circ \delta^I$ is injective.
\item\label{it:cond Q} $\ker((q_I \otimes
    \lambda_G)\circ\delta) = I$.
\end{enumerate}
\end{lemma}
\begin{proof}
The equivalence of \eqref{it:deltaI normal}~and~\eqref{it:PhiI
faithful} is an application of \cite[Lemma~1.4]{qui:discrete
coactions}. The equivalence of \eqref{it:deltaI
normal}~and~\eqref{it:deltar injective} is by definition of
normality; see \cite[Definitions A.39~and~A.50]{EKQR}. To
establish the equivalence of \eqref{it:cond
Q}~and~\eqref{it:deltar injective}, just observe that
\[
(q_I \otimes \lambda_G)\circ\delta
= (\id_{A/I} \otimes \lambda_G) \circ (q_I \otimes \id_{C^*(G)}) \circ \delta
= (\id_{A/I} \otimes \lambda_G) \circ \delta^I \circ q_I.\qedhere
\]
\end{proof}

Recall that a discrete group $G$ is called \emph{exact} if its
reduced $C^*$-algebra $C^*_\reduced(G)$ is exact.

\begin{prop} \label{prop:exact}
Let $G$ be a discrete group. Then the following are equivalent:
\begin{enumerate}
\item $G$ is exact.
\item For every normal coaction $\delta$ of $G$ on a $C^*$-algebra
  $A$, and every ideal $I$ of $A$ which is generated by its
  intersection with $A^\delta_e$, the induced coaction $\delta^I$ of $G$ on
  $A/I$ is normal.
\end{enumerate}
\end{prop}

\begin{proof}
  Assume that $G$ is not exact. Then there exists a $C^*$-algebra $A$
  which has an ideal $I$ such that $I\otimes C^*_\reduced(G)\subsetneq
  \ker (q_I\otimes \id_{C^*_\reduced(G)})$ where $q_I \colon A\to A/I$ is the
  quotient map. Let $\delta_{G}$ denote the coaction of $G$ on $C^*_\reduced(G)$
  given by $\delta_{G}(\lambda_G(g)) = \lambda_G(g) \otimes i_G(g)$
  for all $g \in G$ (see \cite[Example 1.15]{qui:full reduced} or
  \cite[Proposition 2.4]{qui:full duality}). The $*$-homomorphism $\delta:=\id_A\otimes\delta_G$ is
  then a coaction of $G$ on $A\otimes C^*_\reduced(G)$. Let
  $\Delta \colon C^*_\reduced(G)\to C^*_\reduced(G)\otimes C^*_\reduced(G)$ be
  given by $\Delta(x)=x\otimes x$ for all $x\in C^*_\reduced(G)$. We
  then have that
  $(\id_{A\otimes
    C^*_\reduced(G)}\otimes\lambda_G)\circ\delta=\id_A\otimes\Delta$,
  from which it follows that $(\id_{A\otimes
    C^*_\reduced(G)}\otimes\lambda_G)\circ\delta$ is
  injective, and thus that $\delta$ is normal. It is easy to
  check that $(A\otimes C^*_\reduced(G))^\delta_e=A\otimes
  1_{C^*_\reduced(G)}$, and that the ideal $I\otimes C^*_\reduced(G)$ of
  $A\otimes C^*_\reduced(G)$ is generated by its intersection with $A\otimes
  \id_{C^*_\reduced(G)}$. We have that $(q_{I\otimes
    C^*_\reduced(G)}\otimes\lambda_G)\circ\delta=q_I\otimes\Delta$
  from which it follows that
  \[I\otimes C^*_\reduced(G)\subsetneq
  \ker (q_I\otimes \id_{C^*_\reduced(G)})\subset
  \ker(q_I\otimes\Delta) = \ker((q_{I\otimes
    C^*_\reduced(G)}\otimes\lambda_G)\circ\delta),\]
  so it follows from Lemma \ref{lem:when deltaI normal} that
  $\delta^{I\otimes C^*_\reduced(G)}$ is not normal.

  Assume now that $G$ is exact and let $\delta$ be a coaction of $G$
  on a $C^*$-algebra $A$, and $I$ an ideal of $A$ which is generated by its
  intersection with $A^\delta_e$. If $x\in\ker((q_I \otimes
    \lambda_G)\circ\delta)$, then
    $(\id_A\otimes\lambda_G)(\delta(x))\ker (q_I\otimes
    \id_{C^*_\reduced(G)}) = I\otimes C^*_\reduced(G)$, from which it
    follows that $x\in I$. Thus $\ker((q_I \otimes
    \lambda_G)\circ\delta)= I$, so it follows from Lemma \ref{lem:when
      deltaI normal} that $\delta^I$ is normal.
\end{proof}

\begin{rmk}
  The first half of the proof is essentially taken from \cite[page 61]{Exel}, and
  is adapted to our coaction framework.
\end{rmk}


\begin{thebibliography}{00}
\bibitem{ALNR} S. Adji, M. Laca, M. Nilsen, and I. Raeburn,
    \emph{Crossed products by semigroups of endomorphisms and the
    Toeplitz algebras of ordered groups}, Proc.  Amer.  Math.  Soc.
    {\bf 122} (1994), 1133--1141.

\bibitem{Arveson1989} W. Arveson, \emph{Continuous analogues of
    {F}ock space}, Mem. Amer. Math. Soc. \textbf{80} (1989),
    no.~409, iv+66.

\bibitem{Arveson2008} W. Arveson, \emph{The noncommutative {C}hoquet boundary},
  J. Amer. Math. Soc. \textbf{21} (2008), no.~4, 1065--1084.

\bibitem{Black} B. Blackadar, \emph{Operator algebras. Theory
    of
    $C\sp *$-algebras and von Neumann algebras}. Encyclopaedia of
    Mathematical Sciences, 122. Operator Algebras and Non-commutative
    Geometry, III. Springer-Verlag, Berlin, 2006.

%\bibitem{CN2005} S. Campbell and G.A. Niblo, \emph{Hilbert
%    space compression and exactness of discrete groups}, J.
 %   Funct. Anal. {\bf222} (2005), 292--305.

\bibitem{CL2002} J. Crisp and M. Laca, \emph{On the Toeplitz
    algebras of right-angled and finite-type Artin groups},
 J. Austral. Math. Soc. {\bf 72} (2002), 223--245.

\bibitem{CL} J. Crisp and M. Laca, \emph{Boundary quotients and
    ideals of Toeplitz $C\sp*$-algebras of Artin groups}, J.
    Funct. Anal. \textbf{242} (2007), 127--156.

\bibitem{Cuntz1977} J. Cuntz, \emph{Simple {$C\sp*$}-algebras
    generated by isometries}, Comm. Math. Phys. \textbf{57}
    (1977), no.~2, 173--185.

\bibitem{CK1980} J. Cuntz  and W. Krieger, \emph{A class of
    {$C\sp{\ast} $}-algebras and topological {M}arkov chains},
    Invent. Math. \textbf{56} (1980), no.~3, 251--268.

\bibitem{Dinh1991} H.T. Dinh, \emph{Discrete product systems
    and their {$C\sp *$}-algebras}, J. Funct. Anal.
    \textbf{102} (1991), no.~1, 1--34.

\bibitem{Dix:C*-algebras} J. Dixmier, \emph{$C\sp*$-algebras.}
    Translated from the French by Francis Jellett. North-Holland
    Mathematical Library, Vol. 15. North-Holland Publishing Co.,
    Amsterdam-New York-Oxford, 1977.

\bibitem{EKQR} S. Echterhoff, S. Kaliszewski, J. Quigg and I.
    Raeburn, \emph{A categorical approach to imprimitivity theorems
    for $C\sp *$-dynamical systems},  Mem. Amer. Math. Soc. {\bf180}
    (2006), viii+169 pp.

\bibitem{EKQ} S. Echterhoff, S. Kaliszewski and J.  Quigg, \emph{Maximal coactions},
Internat. J. Math. {\bf 15} (2004), 47--61.

\bibitem{EQ}
  S. Echterhoff and J. Quigg, \emph{Induced coactions of discrete groups
  on {$C\sp *$}-algebras}, Canad. J. Math. \textbf{51} (1999), no.~4, 745--770.

\bibitem{Exel} R. Exel, \emph{Amenability for {F}ell bundles},
  J. reine angew. Math. \textbf{492} (1997), 41--73.


\bibitem{ELQ} R. Exel, M. Laca and J. Quigg, \emph{Partial
    dynamical systems and $C^*$-algebras generated by partial
    isometries}, J. Operator Theory {\bf47} (2002), 169--186.

\bibitem{Fell} J.M.G. Fell and R.S. Doran,
    \emph{Representations of {$\sp *$}-algebras, locally
    compact groups, and {B}anach {$\sp *$}-algebraic bundles.
    {V}ol. 2}. Pure and Applied Mathematics, Vol. 126. Academic
    Press Inc., Boston, MA, 1988.

\bibitem{Fowler1999} N.J. Fowler, \emph{Compactly-aligned
    discrete product systems, and generalization of {$\mathcal{O}_\infty$}}, Internat. J. Math. \textbf{10} (1999), no.~6,
    721--738.

\bibitem{F99} N.J. Fowler, \emph{Discrete product systems of
    Hilbert bimodules}, Pacific J. Math.  {\bf 204} (2002), 335--375.


\bibitem{FMR} N.J. Fowler, P.S. Muhly and I. Raeburn, \emph{Representations of Cuntz-Pimsner algebras}, Indiana
Univ. Math. J. \textbf{52} (2003), 569--605.
%\bibitem{Ham} M. Hamana, \emph{Injective envelopes of operator systems}, Publ. Res.
%  Inst. Math. Sci. \textbf{15} (1979), no.~3, 773--785.

\bibitem{HR1997} A. an Huef and I. Raeburn,
    \emph{The ideal structure of Cuntz-Krieger algebras},
     Ergod. Th. \& Dynam. Sys. \textbf{17} (1997), 611--624

\bibitem{ka1} T. Katsura, \emph{A class of $C^{*}$-algebras
    generalizing both graph algebras and homeomorphism
    $C^{*}$-algebras I, fundamental results}, Trans. Amer.
    Math. Soc. \textbf{356} (2004), 4287--4322.


\bibitem{ka2} T. Katsura, \emph{On $C^*$-algebras associated with $C^*$-correspondences},
J. Funct. Anal. \textbf{217}, 366-401.

\bibitem{ka3} T. Katsura, \emph{Ideal structure of $C^*$-algebras associated with $C^*$-correspondences},
Pacific J. Math. \textbf{230} (2007), 107--145.

\bibitem{L} M. Laca, \emph{Purely infinite simple
    {T}oeplitz algebras}, J. Operator Theory \textbf{41}
    (1999), no.~2, 421--435.

%\bibitem{LR} M. Laca and I. Raeburn, \emph{Semigroup crossed
 %   products and the Toeplitz algebras of nonabelian groups},
  %  J. Funct. Anal. {\bf 139} (1996), 415--440.

\bibitem{Lan} E.C. Lance, \emph{Hilbert $C^*$-modules: A
    toolkit for operator algebraists}, London Math. Soc. Lecture Note Series,
    vol. 210, Cambridge Univ. Press, Cambridge, 1994.

%\bibitem{land} M. Landstad, \emph{Duality theory for
 %   covariant systems}, Trans. Amer. Math. Soc. \textbf{248}
  %  (1979), 223--267.

%\bibitem{LPRS}
 %   M. Landstad, J. Phillips, I. Raeburn and C. Sutherland,
  %  \emph{Representations of crossed products by coactions and
   % principal bundles}, Trans. Amer. Math. Soc. \textbf{299}
    %(1987), 747--784.

\bibitem{N} A. Nica, \emph{$C^*$-algebras generated by isometries and
Wiener-Hopf operators}, J. Operator Theory, \textbf{27} (1992),
17--52.

\bibitem{Pimsner1997} M.V. Pimsner, A class of {$C\sp
    *$}-algebras generalizing both {C}untz-{K}rieger algebras
    and crossed products by {${\bf Z}$}, \emph{Free probability
    theory (Waterloo, ON, 1995)}, Amer. Math. Soc., Providence,
    RI, 1997, 189--212.

\bibitem{qui:full duality} J. Quigg, \emph{Full $C^*$-crossed
    product duality}, J. Austral. Math. Soc. \textbf{50}
    (1991), 34--52.

\bibitem{qui:discrete coactions}
    J. Quigg, \emph{Discrete coactions and $C^*$-algebraic
    bundles}, J. Austral. Math. Soc. \textbf{60} (1996),
    204--221.

\bibitem{qui:full reduced} J. Quigg, \emph{Full and reduced
    $C^*$-coactions}, Math. Proc. Camb. Phil. Soc. {\bf 116}
    (1994), 435--450.

\bibitem{QuiggRa} J. Quigg and I. Raeburn, \emph{Characterisations of crossed products by partial actions}, J. Operator
Theory {\bf 37} (1997), 311--340.

\bibitem{RS} I. Raeburn and A. Sims, \emph{Product systems of graphs and the {T}oeplitz
  algebras of higher-rank graphs}, J. Operator Theory \textbf{53} (2005),
  no.~2, 399--429.

\bibitem{TFB} I. Raeburn and D.P. Williams, \emph{Morita
    equivalence
    and continuous-trace $C^*$-algebras}, Math. Surveys and
    Monographs, vol. 60, Amer. Math. Soc., Providence, 1998.

\bibitem{r} J. Renault, \emph{A Groupoid Approach to $C^{*}$-algebras}, Lecture Notes in Mathematics, Vol. 793, Springer-Verlag, Berlin, 1980.

\bibitem{SY} A. Sims and T. Yeend, \emph{Cuntz-Nica-Pimsner
    algebras
    associated to product systems of Hilbert bimodules}, J. Operator
    Theory, to appear.

\bibitem{Yam09}
  S. Yamashita, \emph{Cuntz's $ax+b$-semigroup $C^*$-algebra over $\mathbb{N}$
  and product system $C^*$-algebras}, arXiv:0906.1994v1 (2009), 18 pp.

\bibitem{y} T. Yeend, \emph{Groupoid models for the $C^{*}$-algebras of topological higher-rank graphs}, J. Operator Theory \textbf{57} (2007), 95--120.

\end{thebibliography}
\end{document}